%% file: arxiv.tex
\documentclass[12pt]{article}

\input{header}
\usepackage[margin=3cm]{geometry}

\renewcommand{\bigvee}{\max}
\begin{document}

\title{\bf High-level methods for homotopy\\construction in associative $n$-categories}

\author{
\vc{David Reutter${}^{\,1}$\\\texttt{david.reutter@cs.ox.ac.uk}}
\hspace{10pt}
\vc{Jamie Vicary${}^{\,1\,2}$\\\texttt{vicaryjo@cs.bham.ac.uk}}
\vspace{-30pt}
}

\author{
\hspace{-5cm}
\def\authorboxwidth{6cm}
\makebox[\authorboxwidth]{\vc{\,\,\,David Reutter${}^{\,1}$\\\small\texttt{david.reutter@cs.ox.ac.uk}}}
\makebox[\authorboxwidth]{\vc{\,\,\,Jamie Vicary${}^{\,1\,2}$\\\small\texttt{vicaryjo@cs.bham.ac.uk}}}
\hspace{-5cm}
\vspace{-0pt}
}

\date{16 January 2019}

\maketitle

\begin{abstract}
A combinatorial theory of associative $n$-categories has recently been proposed, with strictly associative and unital composition in all dimensions, and the weak structure arising as a notion of `homotopy' with a natural geometrical interpretation. Such a theory has the potential to serve as an attractive foundation for a computer proof assistant for higher category theory, since it allows composites to be uniquely described, and relieves proofs from the bureaucracy of associators, unitors and their coherence. However, this basic theory lacks a high-level way to construct homotopies, which would be intractable to build directly in complex situations; it is not therefore immediately amenable to implementation.

We tackle this problem by describing a `contraction' operation, which algorithmically constructs complex homotopies that reduce the lengths of composite terms. This contraction procedure allows building of nontrivial proofs by repeatedly contracting subterms, and also allows the contraction of those proofs themselves, yielding in some cases single-step witnesses for complex homotopies.  We prove correctness of this procedure by showing that it lifts connected colimits from a base category to a category of zigzags, a procedure which is then iterated to yield a contraction mechanism in any dimension. We also present \homotopyio, an online proof assistant that implements the theory of associative $n$\-categories, and use it to construct a range of examples that illustrate this new contraction mechanism.
{%
\footnote{Department of Computer Science, University of Oxford, UK}%
\footnote{School of Computer Science, University of Birmingham, UK}%
}
$\begin{tikzpicture}\path [use as bounding box] (0,0) to (0,0);\draw [fill=white, draw=none] (-.5,0) rectangle (0,.4);\end{tikzpicture}$
\end{abstract}

\maketitle

\input{section1-introduction}

\input{section2-zigzag}

\input{section3-contraction}

\input{section4-homotopies}


\input{sectionA-moremethods}

\input{sectionB-moreproofs}

{
\bibliographystyle{IEEEtran} 
\footnotesize
\small
\bibliography{references}
}

\end{document}

%% file: header.tex
\usepackage{tikz}
\usetikzlibrary{cd}
\usepackage{xspace}
\usepackage{amsthm}
\usepackage{amsmath}
\usepackage{amssymb}       
\usepackage{mathabx}
\usepackage[compact]{titlesec}
\usepackage[numbers, sort, compress]{natbib}
\usepackage[hyperfootnotes=false]{hyperref}
\usepackage{docmute}
\usepackage{enumitem}
\setlist[itemize]{leftmargin=13pt}
\usepackage[labelfont=it,textfont=it]{caption}
\usepackage{etoolbox}

\def\homotopyio{{\emph{homotopy.io}}\xspace}
\pgfdeclarelayer{background}
\pgfdeclarelayer{foreground}
\pgfsetlayers{background,main,foreground}
\newcommand\alignedgraphics[2]{\begin{aligned}\includegraphics[#1]{#2}\end{aligned}}

\theoremstyle{plain}
\newtheorem{theorem}{Theorem}
\newtheorem{lemma}[theorem]{Lemma}
\newtheorem{proposition}[theorem]{Proposition}
\newtheorem{corollary}[theorem]{Corollary}

\theoremstyle{definition}
\newtheorem{remark}[theorem]{Remark}
\newtheorem{definition}[theorem]{Definition}

\newtheorem{example}[theorem]{Example}

\newcommand\vc[1]{\begin{tabular}{@{}c}#1\end{tabular}}
\newcommand\cat[1]{\ensuremath{\mathbf{#1}}\xspace}
\newcommand{\Cat}{{\cat{Cat}}}
\newcommand{\Hom}{\mathrm{Hom}}
\newcommand{ \ignore}[1]{}

\newcommand{\IC}{\mathrm{IC}}
\newcommand{\ob}{\mathrm{ob}}
\newcommand{\op}{\mathrm{op}}
\newcommand{\sing}{\mathrm{sing}}
\newcommand{\reg}{\mathrm{reg}}

\newcommand{\Set}{\mathrm{Set}}

\newcommand{\svdots}{\smash{\vdots}\vphantom{R_1}}
\tikzcdset{arrow style=tikz, diagrams={>={Straight Barb[scale=0.6]}}}

\newcommand{\id}{\mathrm{id}}
\renewcommand{\to}[1][]{\ensuremath{\xrightarrow{\smash{#1}}}}
\newcommand{\ot}[1][]{\ensuremath{\xleftarrow{\smash{#1}}}}

\def\Ob{\mathrm{Ob}}
\def\Mor{\mathrm{Mor}}

\usepackage{stmaryrd}

\allowdisplaybreaks
\renewcommand\paragraph[1]{

\noindent
\textbf{#1.} }

\usepackage[]{footmisc}

\makeatletter
\def\calign@preamble{%
   &\hfil\strut@
    \setboxz@h{\@lign$\m@th\displaystyle{##}$}%
    \ifmeasuring@\savefieldlength@\fi
    \set@field
    \hfil
    \tabskip\alignsep@
}
\let\cmeasure@\measure@
\patchcmd\cmeasure@{\divide\@tempcntb\tw@}{}{}{}
\patchcmd\cmeasure@{\divide\@tempcntb\tw@}{}{}{}
\patchcmd\cmeasure@{\ifodd\maxfields@
  \global\advance\maxfields@\@ne
  \fi}{}{}{}    
\newenvironment{calign}
{%
  \let\align@preamble\calign@preamble
  \let\measure@\cmeasure@
  \align
}
{%
  \endalign
}  
\makeatother

\definecolor{DRcolor}{rgb}{0.0,0.8,0.25}        
\definecolor{CDcolor}{rgb}{0.8,0.0,0.2}         
\definecolor{JVcolor}{rgb}{0.2,0.1,0.85}        
\newcommand{\DR}[1]{{\color{DRcolor}*}\marginpar{\vspace*{-20pt}\tiny\color{DRcolor}{* #1}\vspace*{20pt}}}

\newcommand{\DRs}[1]{}
\newcommand{\CDs}[1]{}
\newcommand{\DRcomm}[1]{{\color{DRcolor}{#1}}}

\newenvironment{tz}[1][]{\begin{tikzpicture}[baseline={([yshift=-.8ex]current bounding box.center)},#1]}{\end{tikzpicture}}

\def\cA{\cat A}
\def\cB{\cat B}
\def\cC{\cat{C}}
\def\cD{\cat D}
\def\cE{\cat E}

\def\cJ{\cat J}

\newcommand{\colim}{\mathrm{colim}}
\newcommand\N{\mathbb{N}\xspace}
\def\Simp{\Delta}

%% file: section1-introduction.tex
\begin{figure}[!b]
\[
\alignedgraphics{width=5.25cm}{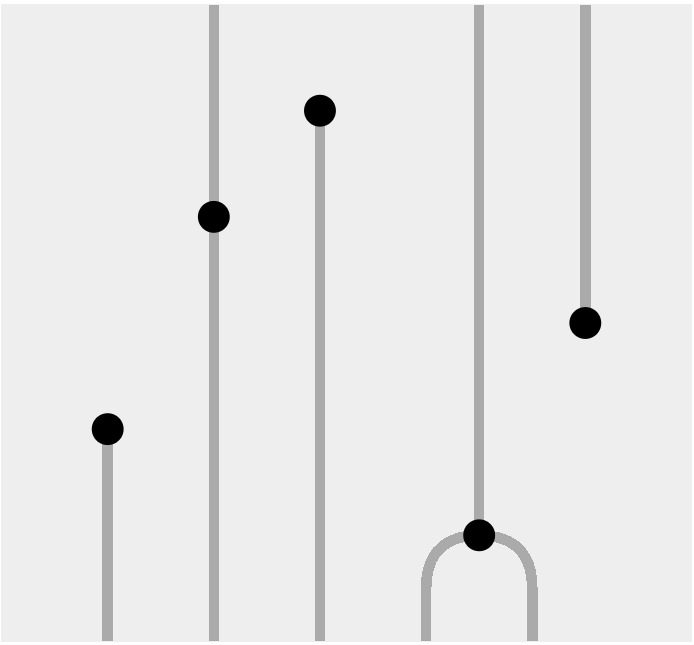} 
\quad\leadsto\quad
\alignedgraphics{width=5.25cm}{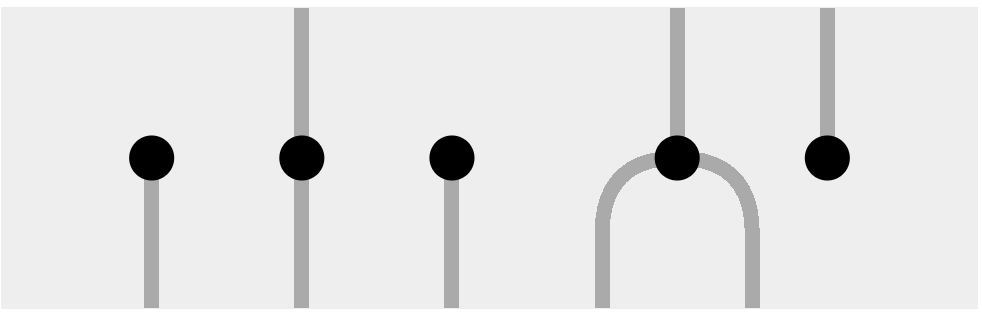} 
\]
\caption{\label{fig:algebraiccontraction}A 3-dimensional homotopy of 2-dimensional diagrams.
\href{http://www.cs.bham.ac.uk/~vicaryjo/homotopy.io/lics2019/figure_1.html}{(Link to online proof)}}
\end{figure}

\section{Introduction}

\noindent
The theory of associative $n$\-categories (ANCs) has recently been proposed~\cite{Dorn2018, Douglas2019}. As with the theory of strict $n$\-categories~\cite{Leinster2003}, composition in this theory is strictly associative and unital in all dimensions. However, unlike the strict theory, ANCs retain a significant amount of weak structure---in the form of homotopies\footnote{\textit{Homotopy} is a standard notion from algebraic topology, which can be understood informally to mean the  continuous deformation of one topological structure into another. We use the term only in an informal sense, basing our formal mathematical development on the theory of associative $n$-categories, which have a combinatorial foundation.}, with a natural geometric interpretation---making it reasonable to conjecture that every weak $n$\-category is equivalent to an ANC~\cite[Conjecture I.5.0.4]{Dorn2018}.\footnote{In dimension~3 the theory of ANCs agrees with the theory of Gray categories, a well-known model of 3-categories which is equivalent to the fully-weak theory~\cite{Gordon1995}, but which has strict associators and unitors; for \mbox{$n>3$}, ANCs are not expected to be equivalent to any previously-described theory.}

\begin{figure}[!b]
\[
\begin{aligned}
\includegraphics[width=3cm]{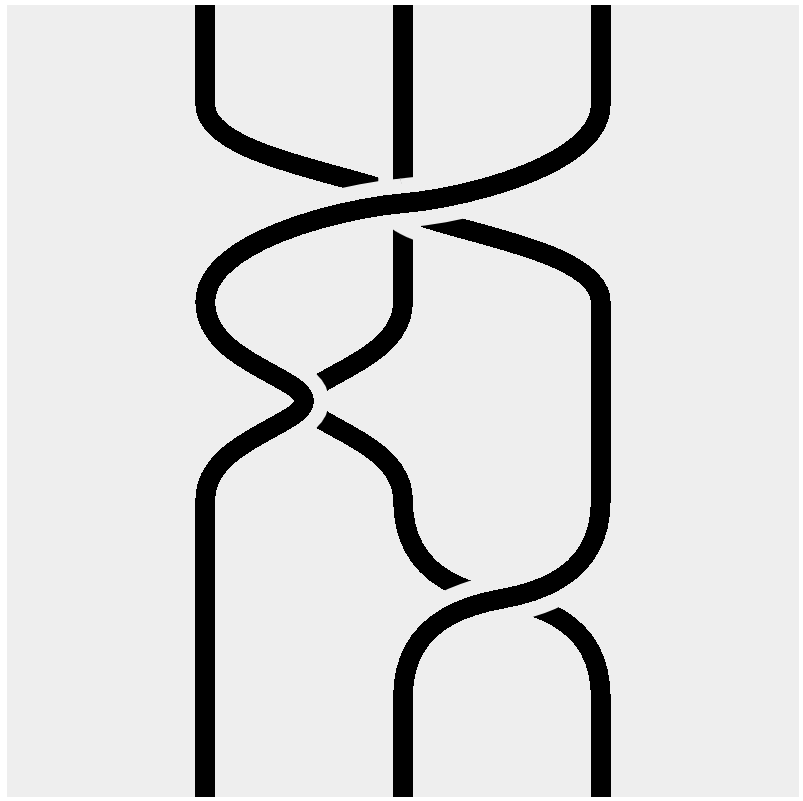}
\end{aligned}
\quad\leadsto\quad
\begin{aligned}
\includegraphics[width=4cm]{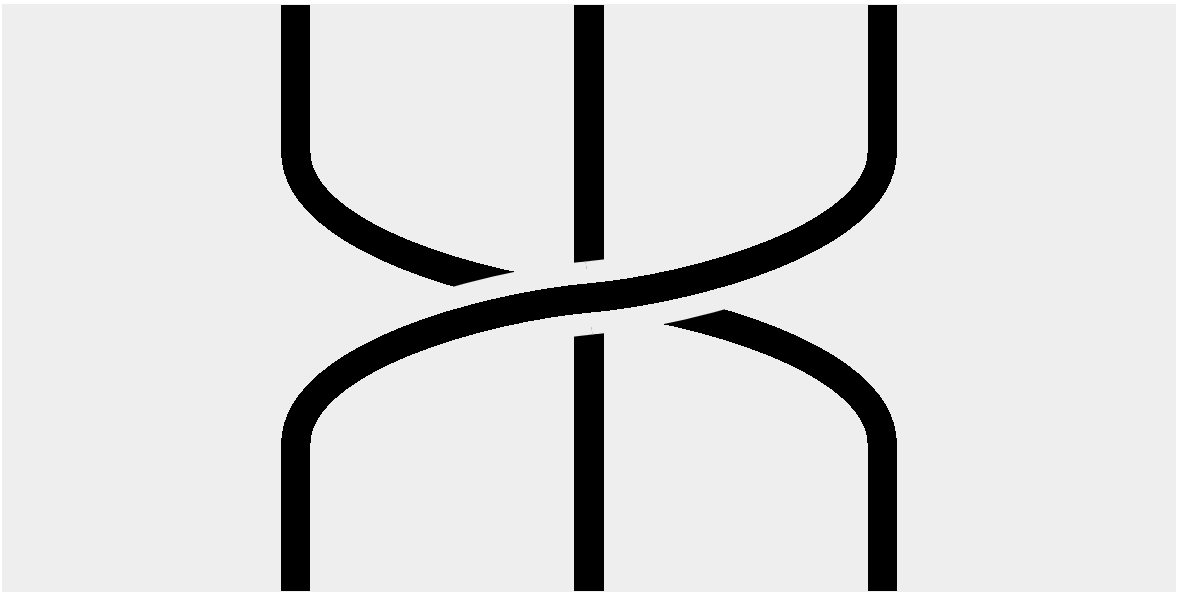}
\end{aligned}
\]
\caption{\label{fig:tanglecontraction}A 4-dimensional homotopy of 3-dimensional diagrams.
\href{http://www.cs.bham.ac.uk/~vicaryjo/homotopy.io/lics2019/figure_2.html}{(Link to online proof)}}
\end{figure}

It may therefore be possible for  ANCs to serve as an attractive general language for calculations in higher category theory, if suitably encoded into a computer proof assistant. Strict associators and units would make composites unique, eliminating some of the bureaucracy of coherence, while the remaining weak structure---while still potentially of high complexity---could be reasoned about geometrically.

The major obstacle to realizing this goal is the difficulty of constructing nontrivial terms of the theory. In principle these can encode complex data, including not only the composites of generating types in arbitrary dimension, but also arbitrary homotopies of these composites. Each term has a dimension, and the $n$\-dimensional terms are called \textit{$n$\-diagrams}.

As examples of such terms, consider Figures~\ref{fig:algebraiccontraction} and~\ref{fig:tanglecontraction}. Figure~\ref{fig:algebraiccontraction} shows two 2\-diagrams, which can be interpreted as expressions in the string diagram calculus for a monoidal category~\cite{Selinger2011}; the arrow represents a homotopy of these 2\-diagrams, and forms part of the data of a 3\-diagram. Figure~\ref{fig:tanglecontraction} shows two 3\-diagrams, which can be interpreted as tangles~\cite{Baez1995} in the string diagram calculus for a braided monoidal category; the arrow represents a homotopy of these 3\-diagrams, and forms part of the data of a 4\-diagram.

The main mathematical contribution of this paper is the description of a \textit{contraction} algorithm, which builds homotopies that reduce a given portion of the diagram in size, along with a full mathematical theory that demonstrates correctness of the procedure. This gives a  high-level method for building nontrivial homotopies in an associative $n$\-category. Figures~\ref{fig:algebraiccontraction} and~\ref{fig:tanglecontraction} both give examples; in each case, the second diagram was obtained by executing the contraction procedure on the first diagram.

Contraction can serve as the main workhorse for the construction of a range of nontrivial proofs in the theory. Given an initial composite $n$\-diagram, we produce our $(n+1)$-dimensional proof object by contracting various  $k$\-dimensional subdiagrams for $k \leq n$ to produce the content of the $(n+1)$\-dimensional proof object, as well as applying algebraic moves from the signature, and extending these recursively to the diagram as a whole using some further techniques described in Section~\ref{sec:additionalmethods}. Once our proof is complete, we can then contract that proof term itself, to yield a shorter proof of the same logical statement. In Section~\ref{sec:homotopies} we give two fully-worked examples of this entire proof construction workflow. Indeed, we conjecture that contractions, together the associated recursive techniques that we discuss in Section~\ref{sec:additionalmethods}, yield a universal toolkit which can in principle construct any homotopy in the theory.
\input{figures/figure_motivation}

We also present \homotopyio~\cite{nlab:homotopyio}, an online proof assistant that implements the theory of associative $n$\-categories. This proof assistant is enabled by our new theory of contractions, which serves as the main engine for homotopy construction, and is applied by clicking and dragging with the mouse on the graphical representation rendered by the tool. We present the tool  as an accompaniment to the claims of the paper, demonstrating that the theory of contractions that we build here is useful and practical.

We keep the focus of this paper on logical foundations, and do not give further details on the implementation. Nonetheless, we accompany many of our examples with direct hyperlinks to their online formalization in the tool, which we invite the reader to investigate. To explore these workspaces, change the parameters of the ``Slice'' control at the top-right of the window; to manipulate the diagrams homotopically, use the mouse to drag vertices (or crossings) up or down, or drag wires left or right.


\subsection{Related work}
\label{sec:relatedwork}

\noindent
This work builds on the existing theory of ANCs due to Dorn, Douglas and Vicary~\cite{Dorn2018, Douglas2019}. That theory defines \textit{signatures} that give families of admissible types,  \emph{diagrams} that encode composites and homotopies of these types, \emph{term normalization} which reduces a diagram to a standard form, and \emph{type checking} which verifies whether a normalized diagram is valid.  The tool \homotopyio implements all these core aspects of the theory, about which we give no further details in this paper. However,  that existing theory does not include the concept of contraction, or yield any other high-level method for homotopy construction, motivating our results here.

The theory of ANCs can be seen as a development into arbitrary dimension of the theory of quasistrict 4\-categories of Bar and Vicary~\cite{Bar2017}, implemented as the proof assistant \emph{Globular} by Bar, Kissinger and Vicary~\cite{Bar2018}. That proof assistant had a restricted notion of homotopy construction, limited fundamentally to dimension 4, and could not even in principle be generalized to arbitrary dimension, where our results apply.

Having in hand a high-level method for homotopy construction in arbitrary dimension, it is interesting to ask for an algorithm which, given a pair of $n$\-diagrams, either constructs a correct homotopy between them, or correctly reports that no such homotopy exists. Such an algorithm was recently given for the case of 2\-diagrams~\cite{Delpeuch2018}, running in quadratic time. The general case is known to be decidable by work of Makkai~\cite{Makkai2005}.

Contractions, as we present them in this paper,  are colimit constructions for sequences of cospans. Spans and cospans have seen wide application in the theory of higher categories, in particular by Baez and collaborators~\cite{Baez2009}, Grandis~\cite{Grandis2008}, Morton~\cite{Morton2009} and Stay~\cite{Stay2013}; however, in these approaches, a colimit construction usually yields cospan \textit{composition} is often given as a colimit construction. This highlights a key difference: in our work, we compose cospans just by arranging them side-by-side, with the colimits---which do not always exist---instead giving us the high-level contraction structure.


\subsection{Overview of the paper}

\noindent
Our contribution is structured as follows. In Section~\ref{sec:zigzag} we introduce zigzag categories and explore their properties, culminating in simple definitions of untyped and typed $n$\-diagrams. In Section~\ref{sec:contraction} we give a construction procedure for colimits on a zigzag category in terms of colimits in the base category, prove its correctness, and show how this gives rise to a contraction procedure for diagrams. In Section~\ref{sec:homotopies}, we show how contraction works together with some other simple mechanisms to give a general toolkit for homotopy construction, and we give a wide range of examples.

\subsection{Notation}

\noindent
For a natural number $n \in \N$, we write $[n]$ for the totally-ordered set $\{0, 1, \ldots, n-1\}$. We use boldface capital letters $\cat A, \cat B, \cat C, \ldots$ for categories. We write \cat 1 for the terminal category, with unique object $\bullet$ and only the identity morphism, and  $\Simp$ for the category of (possibly empty) finite totally-ordered sets and order-preserving functions.

\subsection{Acknowledgements}

\noindent
We are grateful to Christoph Dorn and Christopher Douglas for many useful discussions about associative $n$\-categories. The second author acknowledges funding from the Royal Society.

%% file: figures/figure_motivation.tex
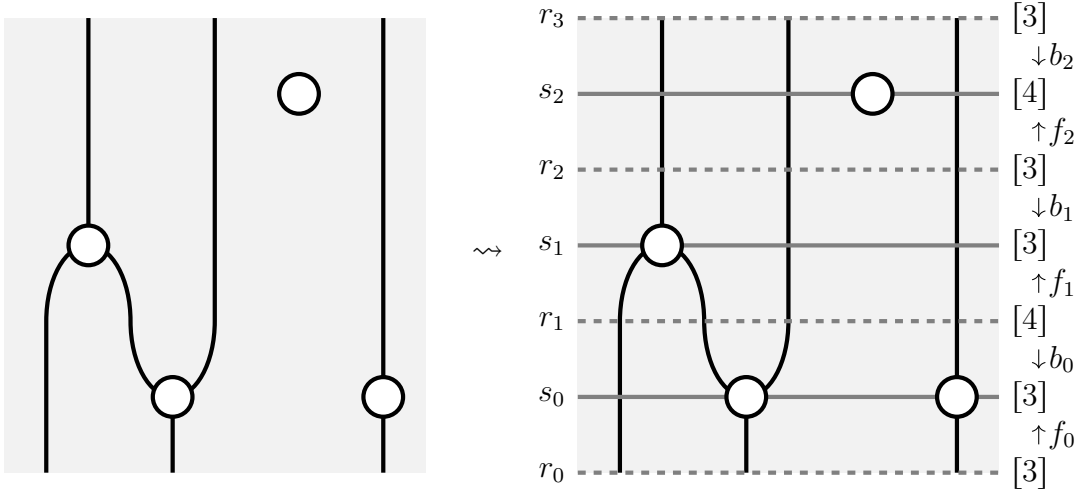
\begin{figure*}
\tikzset{dot/.style={draw=black, circle, fill=white, inner sep=0pt, font=\tiny, minimum width=10pt, fill=white, text=black}}
\tikzset{bluedot/.style={dot, fill=blue!50, text=black, fill=white, text=black}}
\tikzset{yellowdot/.style={dot, fill=yellow!50, text=black}}
\tikzset{bluedot/.style={dot, fill=blue!50, text=black}}
\tikzset{2cell/.style={draw=black, circle, fill=red!50, inner sep=0pt, minimum width=15pt, fill=white, text=black, ultra thick}}
\tikzset{scaffold/.style={->, shorten <=3pt, shorten >=3pt, black, thick}}
\tikzset{equality/.style={shorten <=5pt, shorten >=5pt, double distance=2pt}}
\tikzset{string/.style={black, ultra thick}}
\tikzset{vertical/.style={}}
\tikzset{regular/.style={opacity=0}}
\tikzset{equality/.style={opacity=0}}
\tikzset{horizontal/.style={opacity=0}}
\newcommand\only[1]{}
\[
\begin{aligned}
\begin{tikzpicture}[scale=0.7,yscale=1.2, xscale=0.8, yscale=1.2]
\path [red, ultra thick, use as bounding box] (0,0) rectangle +(10,6);
\tikzset{dot/.style={}}
\begin{pgfonlayer}{foreground}
\node [bluedot, regular] (00) at (0,0) {};
\node [dot] (10) at (1,0) {};
\node [dot, regular] (20) at (2.5,0) {};
\node [dot] (30) at (4,0) {};
\node [dot, regular] (40) at (6.5,0) {};
\node [dot] (50) at (9,0) {};
\node [dot, regular] (60) at (10,0) {};
\node [bluedot, regular] (01) at (0,1) {};
\node [dot] (11) at (1,1) {};
\node [dot, regular] (21) at (2.5,1) {};
\node [2cell] (31) at (4,1) {};
\node [dot, regular] (41) at (6.5,1) {};
\node [2cell] (51) at (9,1) {};
\node [dot, regular] (61) at (10,1) {};
\node [bluedot, regular] (02) at (0,2) {};
\node [dot] (12) at (1,2) {};
\node [dot, regular] (22) at (2,2) {};
\node [dot] (32) at (3,2) {};
\node [dot, regular] (42) at (4,2) {};
\node [dot] (52) at (5,2) {};
\node [dot, regular] (62) at (7,2) {};
\node [dot] (72) at (9,2) {};
\node [dot, regular] (82) at (10,2) {};
\node [bluedot, regular] (03) at (0,3) {};
\node [2cell] (13) at (2,3) {};
\node [dot, regular] (23) at (3.5,3) {};
\node [dot] (33) at (5,3) {};
\node [dot, regular] (43) at (7,3) {};
\node [dot] (53) at (9,3) {};
\node [dot, regular] (63) at (10,3) {};
\node [bluedot, regular] (04) at (0,4) {};
\node [dot] (14) at (2,4) {};
\node [dot, regular] (24) at (3.5,4) {};
\node [dot] (34) at (5,4) {};
\node [dot, regular] (44) at (7,4) {};
\node [dot] (54) at (9,4) {};
\node [dot, regular] (64) at (10,4) {};
\node [bluedot, regular] (05) at (0,5) {};
\node [dot] (15) at (2,5) {};
\node [dot, regular] (25) at (3.5,5) {};
\node [dot] (35) at (5,5) {};
\node [dot, regular] (45) at (6,5) {};
\node [2cell] (55) at (7,5) {};
\node [dot, regular] (65) at (8,5) {};
\node [dot] (75) at (9,5) {};
\node [dot, regular] (85) at (10,5) {};
\node [bluedot, regular] (06) at (0,6) {};
\node (16) at (2,6) {};
\node (36) at (5,6) {};
\node (56) at (9,6) {};
\node (66) at (10,6) {};
\end{pgfonlayer}
\draw [string] (30.center) to (31.center);
\draw [string] (31.center) to [out=left, in=down] (32.center) to [out=up, in=right] (13.center);
\draw [string] (10.center) to (12.center) to [out=up, in=left] (13.center) to (14.center) to (15.center) to (16.center);
\draw [string] (31.center) to [out=right, in=down] (52.center) to (33.center) to (34.center) to (35.center) to (36.center);
\draw [string] (50.center) to (51.center) to (72.center) to (53.center) to (54.center) to (75.center) to (56.center);
\begin{pgfonlayer}{background}
\draw [fill=black!5, draw=none] (00.center) rectangle (66.center);
\end{pgfonlayer}
\end{tikzpicture}
\end{aligned}
\hspace{.5cm}\leadsto\hspace{.8cm}
\begin{aligned}
\begin{tikzpicture}[scale=0.7, yscale=1.2, xscale=0.8, yscale=1.2] 
\path [red, ultra thick, use as bounding box] rectangle +(11,6);
\tikzset{dot/.style={}}
\begin{pgfonlayer}{foreground}
\node [bluedot, regular] (00) at (0,0) {};
\node [dot] (10) at (1,0) {};
\node [dot, regular] (20) at (2.5,0) {};
\node [dot] (30) at (4,0) {};
\node [dot, regular] (40) at (6.5,0) {};
\node [dot] (50) at (9,0) {};
\node [dot, regular] (60) at (10,0) {};
\node [bluedot, regular] (01) at (0,1) {};
\node [dot] (11) at (1,1) {};
\node [dot, regular] (21) at (2.5,1) {};
\node [2cell] (31) at (4,1) {};
\node [dot, regular] (41) at (6.5,1) {};
\node [2cell] (51) at (9,1) {};
\node [dot, regular] (61) at (10,1) {};
\node [bluedot, regular] (02) at (0,2) {};
\node [dot] (12) at (1,2) {};
\node [dot, regular] (22) at (2,2) {};
\node [dot] (32) at (3,2) {};
\node [dot, regular] (42) at (4,2) {};
\node [dot] (52) at (5,2) {};
\node [dot, regular] (62) at (7,2) {};
\node [dot] (72) at (9,2) {};
\node [dot, regular] (82) at (10,2) {};
\node [bluedot, regular] (03) at (0,3) {};
\node [2cell] (13) at (2,3) {};
\node [dot, regular] (23) at (3.5,3) {};
\node [dot] (33) at (5,3) {};
\node [dot, regular] (43) at (7,3) {};
\node [dot] (53) at (9,3) {};
\node [dot, regular] (63) at (10,3) {};
\node [bluedot, regular] (04) at (0,4) {};
\node [dot] (14) at (2,4) {};
\node [dot, regular] (24) at (3.5,4) {};
\node [dot] (34) at (5,4) {};
\node [dot, regular] (44) at (7,4) {};
\node [dot] (54) at (9,4) {};
\node [dot, regular] (64) at (10,4) {};
\node [bluedot, regular] (05) at (0,5) {};
\node [dot] (15) at (2,5) {};
\node [dot, regular] (25) at (3.5,5) {};
\node [dot] (35) at (5,5) {};
\node [dot, regular] (45) at (6,5) {};
\node [2cell] (55) at (7,5) {};
\node [dot, regular] (65) at (8,5) {};
\node [dot] (75) at (9,5) {};
\node [dot, regular] (85) at (10,5) {};
\node [bluedot, regular] (06) at (0,6) {};
\node (16) at (2,6) {};
\node (36) at (5,6) {};
\node (56) at (9,6) {};
\node (66) at (10,6) {};
\end{pgfonlayer}
\draw [string] (30.center) to (31.center);
\draw [string] (31.center) to [out=left, in=down] (32.center) to [out=up, in=right] (13.center);
\draw [string] (10.center) to (12.center) to [out=up, in=left] (13.center) to (14.center) to (15.center) to (16.center);
\draw [string] (31.center) to [out=right, in=down] (52.center) to (33.center) to (34.center) to (35.center) to (36.center);
\draw [string] (50.center) to (51.center) to (72.center) to (53.center) to (54.center) to (75.center) to (56.center);
\begin{pgfonlayer}{background}
\draw [fill=black!5, draw=none] (00.center) rectangle (66.center);
\end{pgfonlayer}
\tikzset{string/.style={opacity=0}}
\tikzset{2cell/.style={minimum width=15pt, inner sep=-4pt, circle, fill=black, draw, fill opacity=1, text opacity=0, font=\tiny, scale=0.5}}
\tikzset{dot/.style={2cell}}
\draw [gray, ultra thick, dashed] (0,0) node [left, black] {$r_0$} to +(10,0) node [right, black] (m0) {$[3]$};
\draw [gray, ultra thick] (0,1) node [left, black] {$s_0$} to +(10,0) node [right, black] (m1) {$[3]$};
\draw [gray, ultra thick, dashed] (0,2) node [left, black] {$r_1$} to +(10,0) node [right, black] (m2) {$[4]$};
\draw [gray, ultra thick] (0,3) node [left, black] {$s_1$} to +(10,0) node [right, black] (m3) {$[3]$};
\draw [gray, ultra thick, dashed] (0,4) node [left, black] {$r_2$} to +(10,0) node [right, black] (m4) {$[3]$};
\draw [gray, ultra thick] (0,5) node [left, black] {$s_2$} to +(10,0) node [right, black] (m5) {$[4]$};
\draw [gray, ultra thick, dashed] (0,6) node [left, black] {$r_3$} to +(10,0) node [right, black] (m6) {$[3]$};
\draw [->] ([xshift=-15pt] m0.north east) to ([xshift=-15pt] m1.south east);
\draw [<-] ([xshift=-15pt] m1.north east) to ([xshift=-15pt] m2.south east);
\draw [->] ([xshift=-15pt] m2.north east) to ([xshift=-15pt] m3.south east);
\draw [<-] ([xshift=-15pt] m3.north east) to ([xshift=-15pt] m4.south east);
\draw [->] ([xshift=-15pt] m4.north east) to ([xshift=-15pt] m5.south east);
\draw [<-] ([xshift=-15pt] m5.north east) to ([xshift=-15pt] m6.south east);
\node at (11.5,.5) {$f_0$};
\node at (11.5,1.5) {$b_0$};
\node at (11.5,2.5) {$f_1$};
\node at (11.5,3.5) {$b_1$};
\node at (11.5,4.5) {$f_2$};
\node at (11.5,5.5) {$b_2$};
\end{tikzpicture}
\end{aligned}
\]
\caption{Decomposing a 2-dimensional string diagram as a sequence of monotone functions.\label{fig:zigzagmotivation}}
\end{figure*}

%% file: section2-zigzag.tex
\section{Zigzag categories}
\label{sec:zigzag}

\noindent
Our theory is based on the notion of zigzag, a reworking and simplification of the notion of \textit{singular interval} from~\cite{Dorn2018, Douglas2019}, and zigzag maps, corresponding to the notion of \emph{limit} in that reference. In this section we develop the theory of zigzags and their maps, and show how they can be used to give definitions of untyped and typed $n$\-dimensional diagrams.

\subsection{Motivation}

\noindent
We motivate the theory of zigzags by examining a 2\-dimensional string diagram, as illustrated on the left of Figure~\ref{fig:zigzagmotivation}, drawn in the standard Joyal-Street graphical calculus for monoidal categories~\cite{Selinger2011, Joyal1993}, and considering how we could represent it combinatorially. At 3 distinct heights, the diagram contains vertices (which we imagine to be pointlike); we call these the \emph{singular heights}, and label them $s_0$, $s_1$ and $s_2$.

If we formally remove these heights, we disconnect the diagram into 4 sections, none of which contain any vertices. For any such section, the geometrical content of any two heights will be equivalent; in particular, the number of wires present at two such heights must be the same, since wires are only created or destroyed by vertices. So we arbitrarily choose one height in each of these sections, called the \emph{regular heights}, and label them $r_0$, $r_1$, $r_2$ and $r_3$.

We now have 7 chosen heights, and at each of them we count how many geometrical entities (either vertices or wires) are present at that height. For example, $r_0$ intersects 3 entities (all wires), and $s_2$ intersects 4 entities (3 wires and 1 vertex). These entities form a totally-ordered set in a natural way, in their order of appearance from left-to-right within each height, and we write the corresponding totally-ordered set $[n]$ at the right of the diagram.

We then choose a regular height $r_i$, and imagine it converging to one of its adjacent singular heights $s_j$. This process will induce an order-preserving function from the entities at height $r_i$ to the entities at height $s_j$, and we write these functions as $f_0$, $b_0$, $f_1$, $b_1$, etc, to the side of the diagram.

From our original diagram, we have therefore obtained a family of totally-ordered sets, and monotone functions between them, with an alternating pattern of directions. This is an instance of the general theory of \emph{singular intervals}~\cite{Dorn2018, Douglas2019}, and directly motivates the more elementary theory of \emph{zigzags}, which we now explore.

\subsection{Zigzags and zigzag maps}

\begin{definition}
In a category \cat C, a \emph{zigzag} $Z$ is a finite diagram of the following sort:
\begin{equation}
\label{eq:zigzag}
\begin{tz}[xscale=0.9, yscale=1, scale=1.5] 
\node (1) at (0,0) {$r_0$};
\node (2) at (1,1) {$s_0$};
\node (3) at (2,0) {$r_1$};
\node (4) at (3,1) {$s_1$};
\node (5) at (4,0) {\phantom{$s_1$}};
\node at (5,0.5) {$\cdots$};
\node (6) at (6,0) {\phantom{$s_1$}};
\node (7) at (7,1) {$s_{n-1}$};
\node (8) at (8,0) {$r_{n}$};
\draw [->] (1) to node [above left, inner sep=0pt] {$f_0$} (2);
\draw [->] (3) to node [above right, inner sep=0pt] {$b_0$} (2);
\draw [->] (3) to node [above left, inner sep=0pt] {$f_1$} (4);
\draw [->, dashed] (5) to node [above right, inner sep=0pt] {$b_1$} (4);
\draw [->, dashed] (6) to node [above left, inner sep=0pt] {$f_{n-1}$} (7);
\draw [->] (8) to node [above right, inner sep=0pt] {$b_{n-1}$} (7);
\end{tz}
\end{equation}
We write $Z_{\sing}=[n]$ for the ordered set of \textit{singular heights}, and $Z_\reg = [n+1]$ for the ordered set of \textit{regular heights}. The objects $r_0, r_1, \ldots$ are called the \emph{regular objects}, and the objects $s_0, s_1, \ldots$ are called the \emph{singular objects}. Such a zigzag has \emph{length $n$}, given by the number of singular heights. Zigzags of length $0$ are allowed, and consist of a single regular height only, and no morphisms. We write $f_i: r_i \to s_i$ and $b_i: r_{i+1} \to s_i$ for the \textit{forward and backward morphisms} in the diagram as indicated, for all $i \in Z_\sing$. Where the zigzag $Z$ is ambiguous, we will write $r_i^Z, s_i^Z, f_i^Z, b_i^Z$ instead of $r_i, s_i, f_i, b_i$.
\end{definition}
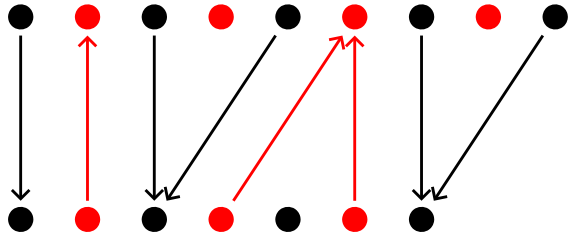
\begin{figure}[b]
\[
\def\regheight{\bullet}
\def\singheight{\textcolor{red}{\bullet}}
\begin{tz}
\node[scale=1.5, rotate=90, yscale=-1, every node/.style={inner sep=0.0pt}, scale=1.5] 
at (0,0){
\begin{tikzcd}[row sep=0.5em, column sep=0.15em]
\\ &&&&&&&&&&&&&&&& \regheight \arrow[to=R1]
\\ &&&&&&&&&&&&&&&& \singheight
\\ |[alias=R1]|\regheight &&&&&&&&&&&&&&&& \arrow[to=R1] \regheight\\
\singheight \arrow[to=S1,red]&&&&&&&&&&&&&&&&  |[alias=S1]|\singheight\\
|[alias=R2]|\regheight &&&&&&&&&&&&&&&&  \arrow[to=R3]\regheight\\
\singheight \arrow[to=S1, red] &&&&&&&&&&&&&&&&  |[alias=S2]|\singheight\\
|[alias=R3]|\regheight  &&&&&&&&&&&&&&&& \arrow[to=R3] \regheight\\
\singheight \arrow[to=S3, red]&&&&&&&&&&&&&&&&  |[alias=S3]|\singheight\\
 |[alias=R4]| \regheight&&&&&&&&&&&&&&&& \arrow[to=R4] \regheight\\
\end{tikzcd}
};
\end{tz}
\]
\caption{An interleaved illustration of a monotone map \mbox{$f:[3] \to{} [4]$} in $\Delta$ going up the page in red, and $f': [5] \to{} [4]$ in $\Delta_=$ going down the page in black.\label{fig:interleaved}}
\end{figure}
\input{figures/figure_zigzagmap}

Before we can define maps of zigzags, we need a short formal development.

\begin{definition}
Let $(-)^\mathrm T:\Simp \to \Simp$ be the functor that adds an element ``on top'' of the total orders, acting on objects as $n \mapsto n+1$, and on a monotone map $f: [n] \to{} [m]$  as $f ^\mathrm T (n) = m$ and $f ^\mathrm T (i) = f(i)$ for $0 \leq i < n$.
\end{definition}

\begin{definition}
\hspace{-3.3pt}Given a monotone map \mbox{$f: [n] \to {} [m]$}, for any element $j \in [m]$ in the target, define \mbox{$f_{\geq j} = \{ i \in [n] | f(i) \geq j \}$} as the elements in the source whose image is above $j$.
\end{definition}

\noindent
For any monotone map $f:[n] \to{} [m]$, for any $j \in [m]$, note that that $(f ^\mathrm T) _{\geq j}$ is always nonempty, due to the additional top element, and thus $\bigvee \big( (f ^\mathrm T) _{\geq j} \big)$ is well-defined.

\begin{definition}\label{def:deltaequal}
The category $\Simp _=$ has nonzero natural numbers as objects, and as morphisms $n \to m$, monotone maps \mbox{$[n] \to{} [m]$} preserving the first and last elements. 
\end{definition}

\begin{definition}
\label{def:k}
The functor $(-)': \Simp \to \Simp _= ^\op$ acts on objects as $n \mapsto n+1$, and acts on morphisms $f: [n] \to{} [m]$ as $f'(j) = \bigvee \big((f^\mathrm T)_{\geq j} \big).$ 
\end{definition}

\noindent
We illustrate this in Figure~\ref{fig:interleaved}, which shows a monotone map $f$ in red, and its ``reversal'' $f'$ in black.

\begin{lemma}
The functor $(-)'$ is an equivalence.
\end{lemma}

\begin{proof}
The functor is clearly surjective on objects. That it is fully faithful can be seen by inspection of Figure~\ref{fig:interleaved}: from any black right-to-left monotone preserving first and last elements, a red left-to-right  monotone can be constructed by ``filling in the gaps'', and vice-versa.
\end{proof}

\noindent
Abusing notation, we also denote the inverse of this functor $\Simp_= \to \Simp$ by $(-)'$, restricting its use in this way to situations where there is no ambiguity.

We now use this technology to define maps of zigzags.

\begin{definition}
\label{def:zigzagmap}
In a category \cat C, a \emph{zigzag map} $f: Z \to Z'$ comprises a monotone function $f_\sing:Z_\sing \to Z'_\sing$, and for each $i \in Z_\sing$ a morphism $f_i: s_i \to s'_{f_\sing(i)}$. Defining $f_\reg = (f_{\sing})':Z_{\reg}' \to Z_{\reg}$, we then require that the diagram constructed as follows, which can always be laid out in a planar way, is commutative: 
\begin{enumerate}[label=(\arabic*)]
\item Take the disjoint union of $Z$ and $Z'$ as diagrams in \cat C.
\item For every $i \in Z_\sing$, add the arrow $f_i$ to the diagram, going from $s \in Z_\sing$ to $f_\sing(s) \in Z' _\sing$.
\item For every $j \in Z'_\reg$, add an identity arrow to  the diagram, between $j \in Z' _\reg$ and $f _\reg(j) \in Z _\reg$.
\end{enumerate}
\end{definition}

\noindent
This construction is quite simple to use in practice. We illustrate it in Figure~\ref{fig:zigzagmap}. Informally, it amounts to the following. (1)~Draw the zigzags $Z$ and $Z'$ one above the other. (2)~For each singular object of $Z$, add an arrow to some singular object of $Z'$, such that the implied function $Z_\sing \to Z'_\sing$ is monotone. (3)~In the spaces between these arrows, add all possible equalities between regular objects of $Z$ and $Z'$.

In the example of Figure~\ref{fig:zigzagmap}, $Z$ has length 4 and $Z'$ has length 5, with $g: Z \to Z'$ running bottom-to-top. The monotone $g_\sing : [4] \to{} [5]$ acts as $0 \mapsto 0$, $1 \mapsto 1$, $2 \mapsto 1$ and $3 \mapsto 4$, with the morphisms $g_0$, $g_1$, $g_2$ and $g_3$ having source and target objects as indicated. The equalities between regular heights  force equalities of regular objects $r_0 = r_0'$, $r_1 = r_1'$, $r_3 = r_2' = r_3' = r_4'$, and $r_4 = r_5'$. The diagram is formed from 9 squares, all of which must commute, leading in this case to the requirements $f_0' = g_0 \circ f_0$, $f_0' = g_0 \circ b_0$, $f_1' = g_1 \circ f_1$, $f_1 \circ b_1 = g_2 \circ f_2$, $b_1' = g_2 \circ b_2$, $f_2' = b_2'$, $f_3' = b_3'$, $b_4' = g_3 \circ f_3$, and $f_4' = g_3 \circ b_3$.

Zigzags and their maps form a category in the obvious way.

\begin{definition}
Given a category \cat C, the \textit{zigzag category} $Z_ \cat C$ is defined to have zigzags as objects and zigzag maps as morphisms.
\end{definition}

\noindent
Composition, associativity and units are clear. We will often be interested in iterating this construction, as follows.
\begin{definition}
Given a category \cat C, the \emph{iterated zigzag category} $Z^n_\cat C$ is the category obtained by starting with the category $\cat C$, and taking the zigzag category $n$ times.
\end{definition}

Every zigzag category has forgetful functors to $\Simp$ and $\Simp_= ^\op$.

\begin{definition}[Regular and singular monotone functors]
For a category \cat C, the \textit{singular monotone functor} $S _\cat C: Z_\cat C \to \Simp$ acts as $S_\cat C(Z) = Z_\sing$ and $S_\cat C(f) = f_\sing$, and the \textit{regular monotone functor} $R_\cat C: Z_\cat C \to \Simp_= ^\op$ acts as $R_\cat C = (-)' \circ S_\cat C$.
\end{definition}

\begin{example}
Starting with the terminal category \cat 1, we see that $Z _\cat 1 = \Simp$, the  singular monotone functor \mbox{$S_\cat 1 : \Simp \to \Simp$} is the identity, and the regular monotone functor \mbox{$R_\cat 1 = (-)'$}.
\end{example}

\input{figures/figure_untyped}

\subsection{Untyped and typed diagrams}
\label{sec:untypedandtyped}

\noindent
We can use zigzag categories to give a straightforward notion of untyped $n$\-diagram, yielding an elementary untyped version of the diagrams which form the terms of the theory of associative $n$\-categories~\cite{Dorn2018}.

\begin{definition}
\label{def:untypedndiagram}
An \textit{untyped $n$\-diagram} is an object of the iterated zigzag category $Z^n _\cat 1$.
\end{definition}

\noindent
We explore this definition through the following examples, illustrated in Figure~\ref{fig:untyped}.

\begin{example}
The only untyped 0\-diagram is $\bullet$, the point.
\end{example}
\begin{example}
An untyped 1\-diagram is an object of $Z_\cat 1=\Simp$, a finite ordinal. So the only parameter is the length of the composite.
\end{example}
\begin{example}
\hspace{-2.9pt}An untyped 2\-diagram is an object of \mbox{$Z_{Z_\cat 1} = Z_\cat \Simp$}.
\end{example}

\noindent
Note in particular that the untyped 2\-diagram illustrated in Figure~\ref{fig:untyped}(c) corresponds exactly to the example of Figure~\ref{fig:zigzagmotivation}, which motivated our construction in the first place.

To develop the theory of typed diagrams, suppose that $L$ is a set of labels equipped with an arbitrary \textit{dimension  function} \mbox{$\mathrm d:L \to \N$}. We can build a poset $\cat L$, whose objects are elements of the set $L$, and where for any $l,l' \in L$, there is a morphism $l \to l'$ in $\cat L$ just when $l=l'$ or $\dim(l) < \dim(l')$. We can use this to give a generalization of Definition~\ref{def:untypedndiagram} appropriate for the typed setting.
\begin{definition}
\label{def:typedndiagram}
For a set of types $L$ equipped with a dimension function $\mathrm d:L \to \N$, an \emph{$L$-typed $n$\-diagram} is an object of the iterated zigzag category $Z^n _{\cat L}$.
\end{definition}

An $L$\-typed $n$\-diagram is a similar structure to an untyped $n$\-diagram, except that every ``bottom-level point'' is assigned an element of the label set $L$, in a way which is arbitrary, except that as we pass from one type to another along a zigzag map, the types must either stay the same, or increase in dimension. Indeed, it is clear that if we choose $L=\{ \bullet \}$ with dimension function $\bullet \mapsto 1$, we we recover the theory of untyped $n$\-diagrams as a special case.

The full theory of associative $n$\-categories~\cite{Dorn2018, Douglas2019} has a notion of \emph{type signature} $\Sigma$, which defines a set of \emph{type labels} $|\Sigma|$, and for each label $l \in |\Sigma|$ a \textit{canonical neighbourhood}. There is then a \textit{type checking scheme}, which takes as input a {$|\Sigma|$\-typed $n$\-diagram}, and returns a boolean, indicating whether or not it is well-typed with respect to $\Sigma$; that is, whether for every instance of every type label in the diagram, its neighbourhood in the diagram normalizes to its canonical neighbourhood. We do not discuss this further, as our contraction procedure operates just at the level of the categories  $Z^n_{\cat L}$.

In the remainder of the paper, the $n$\-diagrams we will draw will be typed, and therefore objects of $Z^n_\cat L$ for some label set $L$, about which we will not give details. (It can be assumed that $L$ is sufficiently large to label all the distinct types of regions, wires and vertices that appear in the diagram.) We will generally use the more attractive ``type notation'' of Figure~\ref{fig:zigzagmotivation} for these diagrams, rather than the bare ``untyped notation'' of Figure~\ref{fig:untyped}. In these $n$\-diagrams, which we typically draw in a 2\-dimensional projection, vertices correspond to labels of dimension $n$, wires correspond to labels of dimension $n-1$, and regions correspond to labels of dimension at most $n-2$.

\subsection{Further zigzag constructions}

\noindent
Here we collect some further technical results on zigzags and zigzag maps, which will be used later.

There is an obvious way in which zigzags can be concatenated, by gluing their diagrams horizontally.
\begin{definition}[Zigzag concatenation]
\label{def:concatenate}
In a category \cat C, given zigzags $Z,Z'$ such that the last regular object of $Z$ equals the first regular object of $Z'$, their \emph{concatenation} is the zigzag $Z \circ Z'$ of length $n^Z + n^{Z'}$, obtained by drawing $Z$ to the left of $Z'$ such that their last and first regular level respectively coincides. For any such $Z,Z'$, given zigzag maps \mbox{$f:Z \to Y$} and $f':Z' \to Y'$, we can also concatenate $f$ and $g$ in a precisely analogous way, yielding $f \circ g : Z \circ Z' \to Y \circ Y'$.
\end{definition}

\noindent
These compositional properties are perhaps unsurprising, given that we will use zigzags as the foundation of our approach to associative $n$\-categories. Also note that zigzag concatenation is strictly associative, a property that is inherited by the theory of associative $n$\-categories for composition in all dimensions.

Functors on base categories extend to zigzag categories, and the zigzag construction as a whole extends to \Cat. We omit the proofs, which are straightforward.
\begin{lemma}[Zigzag functors]
\label{lem:zigzagfunctor}
A functor $F:\cat C \to \cat D$ extends to a \textit{zigzag functor} $Z_F : Z_\cat C \to Z _ \cat D$, acting on objects by direct application to all objects and morphisms in the zigzag diagram, and on morphisms by direct application to the entire commutative diagram defining the zigzag map. Furthermore, if $F$ is fully faithful, so is $Z_F$.
\end{lemma}%
\begin{figure}[b!]
\[
\begin{tz}
\node (1) at (0.75,-.75) {$Z$};
\node (2) at (2,-1) {$Z''$};
\node (3) at (1,-2) {$Z'$};
\draw [->] (1) to node [above] {$h$} (2);
\draw [->] (1) to node [right, pos=0.4] {$g$} (3);
\end{tz}
\leadsto
\begin{tz}
\node[scale=0.7] at (0,0){
\begin{tikzcd}[row sep=0.5em, column sep=0.5em]
&&&& r_2 \arrow[to=R2RT, equal] \arrow[to= R2RB, equal] \arrow[to=R2B,equal] \arrow[dl, "b_1"'] &&&&&&&&&&&
\\&&&   s_1 \arrow[to=S1',"g_1", pos=0.6] \arrow[to=S1'', crossing over, "h_1"] &&&&&&&&&&& |[alias=R2RT]| r_2'' \arrow[dl, "b_1''"] 
\\&& r_1 \arrow[ur, "f_1"] \arrow[dl, "b_1"'] &&&&&&&&&&&  |[alias=S2'']| s_1''
\\& s_0 \arrow[to=S1',"g_0", pos=0.6] \arrow[to=S1'',crossing over, "h_0"] &&&&&&&&&&& |[alias=R2RB]|  r_1''\arrow[ur, "f_1''"'] \arrow[dl,"b_0''"]
\\ r_0 \arrow[ur, "f_0"] \arrow[to= R0B,equal]\arrow[to=R0R, equal,crossing over] &&&&&&&&&&& |[alias = S1'']| s_0''
\\ &&&&&&&&&& |[alias= R0R]| r_0''\arrow[ur, "f_0''"'] 
\\ \phantom{a}
\\
\\ &&& |[alias=R2B]|r_1'\arrow[dl, "b_0'"]
\\ && |[alias=S1']| s_0'
\\ & |[alias = R0B]|r_0' \arrow[ur, "f_0'"']
\end{tikzcd}
};
\end{tz}
\]
\caption{The deconstruction of a diagram in $Z_\cat C$, given as a diagram in \cat C.\label{fig:deconstructedzigzagdiagram}}
\end{figure}\begin{figure}[b]
\def\colar{gray!50} 
\def\collab{gray!70} 
\[
 \begin{tikzcd}[row sep=5em, column sep=1.0em]
|[alias=newr]| \textcolor{\collab}{r_0'} \arrow[r,\colar] & |[alias=news]| \textcolor{\collab}{s_0'} &|[alias=r0t]|r'_1\arrow[r] \arrow[l,\colar] &|[alias=s1t]|s'_1 & |[alias=r2t]|r'_2 \arrow[l] \arrow[r,\colar] & |[alias=s2t]|\textcolor{\collab}{s'_2} & \arrow[l,\colar]\textcolor{\collab}{r'_3}\arrow[r,\colar] & \textcolor{\collab}{s_3'} & \textcolor{\collab}{r'_4}\arrow[l,\colar]
\\
\textcolor{\collab}{r_0} \arrow[to=newr, equal,\colar] \arrow[r,\colar] & \textcolor{\collab}{s_0} \arrow[to=news,\colar] & r_1 \arrow[to=r0t, equal] \arrow[r] \arrow[l,\colar] & s_1 \arrow[to=s1t]& r_2 \arrow[l] \arrow[r] & s_2 \arrow[to=s1t] & \arrow[l]r_3 \arrow[to=r2t,equal]  \arrow[r,\colar]& \textcolor{\collab}{s_3}\arrow[ull,\colar] &\arrow[l,\colar] \textcolor{\collab}{r_4}\arrow[ull,equal,\colar]\arrow[u, equal,\colar]
 \end{tikzcd}
\]
\caption{\label{fig:zigzagmaprestriction}The restriction of a zigzag map $f:Z \to Z'$ (the entire diagram) to the regular heights $1,2 \in Z'_\reg$, yielding $f _{(1,2)} :Z^{}_{(1,3)} \to Z'_{(1,2)}$, (drawn in black.)}
\end{figure}%
\begin{lemma}
The zigzag construction extends to a functor \mbox{$Z: \Cat \to \Cat$}, mapping categories to zigzag categories, and functors to zigzag functors.
\end{lemma}

The equalities in step (3) of Definition \ref{def:zigzagmap} give a strong restriction on which zigzags can have maps between them. In particular, if $f: Z \to Z'$ is a zigzag map, then the first and last regular objects of $Z$ and $Z'$ must be the same. This yields a natural partition of $Z_\cat C$ into a disjoint union of full subcategories, as follows.

\begin{definition}[Local zigzag category]\label{def:local}
Given a category \cat C with chosen objects $A, B$, the \textit{local zigzag category} $Z_\cat C(A, B)$ is the full subcategory of $Z_ \cat C$ containing zigzags whose first regular object is $A$, and whose last regular object is $B$.
\end{definition}

\begin{lemma}[Decomposition]
$Z_ \cat C = \coprod _{A, B \in \Ob(\cat C)} Z_ \cat C(A, B)$.
\end{lemma}

Zigzag maps are defined in terms of the construction of a commutative diagram. This construction is important, and we make it completely explicit with the following definition.

\begin{definition}
\label{def:zigzagmapdiagram}
Given a zigzag map $f: Z \to Z'$, its \emph{zigzag map diagram} is the corresponding diagram (as presented in Figure~\ref{fig:zigzagmap}) with which it was defined.
\end{definition}

\noindent
This gives us a formal way to project diagrams in $Z_\cat C$ to give diagrams in $\cat C$, which we illustrate in Figure~\ref{fig:deconstructedzigzagdiagram}.

\begin{definition}
\label{def:deconstruction}
For a diagram $D: \cat J \to Z_\cat C$, its \textit{deconstruction} $D^*: \cat J^* \to \cat C$ is the diagram obtained by taking the union of the diagrams of the zigzag maps $D(f)$ for all $f \in \Mor(\cat J)$.
\end{definition}

\noindent
More precisely, the objects of the deconstructed diagram category $\cat J^*$ are given by a choice of $j \in \Ob(\cat J)$, and a choice of a regular or singular height of $D(j)$; we write such an object as $\big(j, r_i ^{D(j)} \big)$ or $\big(j, s_i ^{D(j)} \big)$, where $r_i^{D(j)} \in D(j)_{\reg}$ and $s_i^{D(j)} \in D(j)_{\sing}$. The morphisms of $\cat J ^*$ are given by adding for all $j \in \Ob(\cat J)$ and $i \in D(j)_\sing$ morphisms $(j,r_{i} ^{D(j)}) \to (j, s_i ^{D(j)}) \ot (j,r_{i+1} ^{D(j)})$, and for all $f \in \Mor(\cat J)$ with $f:j \to j'$, additional morphisms between singular heights  $\smash{\big( j,s_i ^{D(j)} \big) \to \big( j', s ^{D(j')} _{D(f)_\sing(i)} \big)}$ and regular heights $\big( j',r_i ^{D(j')} \big) \to \big( j,r ^{D(j)} _{D(f)_\reg(i)} \big)$.

Given a zigzag map $f:Z \to Z'$, we can restrict it to some contiguous subset of $Z'_\reg$. We illustrate this idea in Figure~\ref{fig:zigzagmaprestriction}, and develop it formally as follows. Here and throughout, the function $f_\reg: Z'_\reg \to Z_\reg$ is understood to act on pairs $a,b \in Z'_\reg$ elementwise.

\begin{definition}[Zigzag restriction]
For a zigzag $Z$ and a pair $a,b \in Z_\reg$ with $a \leq b$, the \textit{restricted zigzag} $Z_{(a,b)}$ is that part of the zigzag diagram for $Z$ that includes the regular objects $r_a$ and $r_b$, and everything in between.
\end{definition}

\begin{definition}[Zigzag map restriction]
For a zigzag map \mbox{$f: Z \to Z'$} and $a,b \in Z'_\reg$ with $a<b$, the \textit{restricted zigzag map} $f_{(a,b)} : Z_{f_\reg(a,b)} \to Z' _{(a,b)}$ is that part of the zigzag map diagram for $f$ that includes the zigzag diagrams for $Z_{f_\reg (a,b)}$ and $Z'_{(a,b)}$, and the morphisms going between these parts.
\end{definition}

\noindent

%% file: figures/figure_zigzagmap.tex
\begin{figure*}[h]
\vspace{1.5cm}
\[
\def\yoff{0.0cm}
\begin{tz}[xscale=1.5, yscale=.5, scale=.8, scale=1.2]
\node (1) at (0,0) {$r_0'$};
\node (2) at (1,\yoff) {$s_0'$};
\node (3) at (2,0) {$r_1'$};
\node (4) at (3,\yoff) {$s_1'$};
\node (5) at (4,0) {$r_2'$};
\node (6) at (5,\yoff) {$s_2'$};
\node (7) at (6,0) {$r_3'$};
\node (8) at (7,\yoff) {$s_3'$};
\node (9) at (8,0) {$r_4'$};
\node (10) at (9,\yoff) {$s_4'$};
\node (11) at (10,0) {$r_5'$};
\draw [->] (1) to node [above] {$f_0'$} (2);
\draw [->] (3) to node [above] {$b_0'$} (2);
\draw [->] (3) to node [above] {$f_1'$} (4);
\draw [->] (5) to node [above] {$b_1'$} (4);
\draw [->] (5) to node [above] {$f_2'$} (6);
\draw [->] (7) to node [above] {$b_2'$} (6);
\draw [->] (7) to node [above] {$f_3'$} (8);
\draw [->] (9) to node [above] {$b_3'$} (8);
\draw [->] (9) to node [above] {$f_4'$} (10);
\draw [->] (11) to node [above] {$b_4'$} (10);

\begin{scope}[xshift=0cm,yshift=-5cm]
\node (1b) at (0,0) {$r_0$};
\node (2b) at (1,\yoff) {$s_0$};
\node (3b) at (2,0) {$r_1$};
\node (4b) at (3,\yoff) {$s_1$};
\node (5b) at (4,0) {$r_2$};
\node (6b) at (5,\yoff) {$s_2$};
\node (7b) at (6,0) {$r_3$};
\node (8b) at (7,\yoff) {$s_3$};
\node (9b) at (8,0) {$r_4$};
\draw [->] (1b) to node [below] {$f_0$} (2b);
\draw [->] (3b) to node [below] {$b_0$} (2b);
\draw [->] (3b) to node [below] {$f_1$} (4b);
\draw [->] (5b) to node [below] {$b_1$} (4b);
\draw [->] (5b) to node [below] {$f_2$} (6b);
\draw [->] (7b) to node [below] {$b_2$} (6b);
\draw [->] (7b) to node [below] {$f_3$} (8b);
\draw [->] (9b) to node [below] {$b_3$} (8b);
\end{scope}

\tikzset{eq/.style={double equal sign distance}}
\draw [eq] (1b) to (1);
\draw [eq] (3b) to (3);
\draw [eq] (7b) to (5);
\draw [eq] (7b) to (7);
\draw [eq] (7b) to (9);
\draw [eq] (9b) to (11);

\draw [->] (2b) to node [left, pos=0.5] {$g_0$} (2);
\draw [->] (4b) to node [left, pos=0.5] {$g_1$} (4);
\draw [->] (6b) to node [left, pos=0.5] {$g_2$} (4);
\draw [->] (8b) to node [left, pos=0.5] {$g_3$} (10);

\end{tz}
\]
\caption{A zigzag map diagram $g:Z \to Z'$, running bottom-to-top.\label{fig:zigzagmap}}
\end{figure*}

%% file: figures/figure_untyped.tex
\begin{figure*}
\tikzset{vertical/.style={}}
\tikzset{dot/.style={draw=black, circle, fill=white, inner sep=0pt, font=\tiny, minimum width=10pt, fill=white, text=black}}
\tikzset{bluedot/.style={dot, fill=blue!50, text=black, fill=white, text=black}}
\tikzset{yellowdot/.style={dot, fill=yellow!50, text=black}}
\tikzset{bluedot/.style={dot, fill=blue!50, text=black}}
\tikzset{2cell/.style={draw=black, circle, fill=red!50, inner sep=0pt, minimum width=15pt, fill=white, text=black, ultra thick}}
\tikzset{scaffold/.style={->, shorten <=3pt, shorten >=3pt, black, thick}}
\tikzset{equality/.style={shorten <=5pt, shorten >=5pt, double equal sign distance}}
\tikzset{string/.style={black, ultra thick}}
\tikzset{regular/.style={opacity=1}}
\tikzset{horizontal/.style={opacity=1}}
\tikzset{dot/.style={minimum width=0pt, inner sep=-1pt, circle, fill=black, draw, fill opacity=1, text opacity=0}}
\begin{calign}
\nonumber
\begin{aligned}
\begin{tikzpicture}[xscale=0.6]
\node (1) [dot] at (0,0) {};
\end{tikzpicture}
\end{aligned}
&
\begin{aligned}
\begin{tikzpicture}[xscale=0.6]
\node (1) [dot] at (0,0) {};
\node (2) [dot] at (1,0) {};
\node (3) [dot] at (2,0) {};
\node (4) [dot] at (3,0) {};
\node (5) [dot] at (4,0) {};
\node (6) [dot] at (5,0) {};
\node (7) [dot] at (6,0) {};
\draw [scaffold, horizontal, ->] (1) to (2);
\draw [scaffold, horizontal, ->] (3) to (2);
\draw [scaffold, horizontal, ->] (3) to (4);
\draw [scaffold, horizontal, ->] (5) to (4);
\draw [scaffold, horizontal, ->] (5) to (6);
\draw [scaffold, horizontal, ->] (7) to (6);
\end{tikzpicture}
\end{aligned}
&
\begin{aligned}
\begin{tikzpicture}[xscale=0.7, scale=0.8, scale=0.8]
\path [red, ultra thick, use as bounding box] (0,0) rectangle +(10,6);
{
\tikzset{vertical/.style={opacity=0}}
\tikzset{dot/.style={opacity=0}}
}
\begin{pgfonlayer}{foreground}
\node [dot, regular] (00) at (0,0) {};
\node [dot] (10) at (1,0) {};
\node [dot, regular] (20) at (2.5,0) {};
\node [dot] (30) at (4,0) {};
\node [dot, regular] (40) at (6.5,0) {};
\node [dot] (50) at (9,0) {};
\node [dot, regular] (60) at (10,0) {};
\node [dot, regular] (01) at (0,1) {};
\node [dot] (11) at (1,1) {};
\node [dot, regular] (21) at (2.5,1) {};
\node [dot] (31) at (4,1) {};
\node [dot, regular] (41) at (6.5,1) {};
\node [dot] (51) at (9,1) {};
\node [dot, regular] (61) at (10,1) {};
\node [dot, regular] (02) at (0,2) {};
\node [dot] (12) at (1,2) {};
\node [dot, regular] (22) at (2,2) {};
\node [dot] (32) at (3,2) {};
\node [dot, regular] (42) at (4,2) {};
\node [dot] (52) at (5,2) {};
\node [dot, regular] (62) at (7,2) {};
\node [dot] (72) at (9,2) {};
\node [dot, regular] (82) at (10,2) {};
\node [dot, regular] (03) at (0,3) {};
\node [dot] (13) at (2,3) {};
\node [dot, regular] (23) at (3.5,3) {};
\node [dot] (33) at (5,3) {};
\node [dot, regular] (43) at (7,3) {};
\node [dot] (53) at (9,3) {};
\node [dot, regular] (63) at (10,3) {};
\node [dot, regular] (04) at (0,4) {};
\node [dot] (14) at (2,4) {};
\node [dot, regular] (24) at (3.5,4) {};
\node [dot] (34) at (5,4) {};
\node [dot, regular] (44) at (7,4) {};
\node [dot] (54) at (9,4) {};
\node [dot, regular] (64) at (10,4) {};
\node [dot, regular] (05) at (0,5) {};
\node [dot] (15) at (2,5) {};
\node [dot, regular] (25) at (3.5,5) {};
\node [dot] (35) at (5,5) {};
\node [dot, regular] (45) at (6,5) {};
\node [dot] (55) at (7,5) {};
\node [dot, regular] (65) at (8,5) {};
\node [dot] (75) at (9,5) {};
\node [dot, regular] (85) at (10,5) {};
\node [dot, regular] (06) at (0,6) {};
\node [dot] (16) at (2,6) {};
\node [dot, regular] (26) at (3.5,6) {};
\node [dot] (36) at (5,6) {};
\node [dot, regular] (46) at (7,6) {};
\node [dot] (56) at (9,6) {};
\node [dot, regular] (66) at (10,6) {};
\end{pgfonlayer}
\draw [scaffold, horizontal] (00) to (10);
\draw [scaffold, horizontal] (20) to (10);
\draw [scaffold, horizontal] (20) to (30);
\draw [scaffold, horizontal] (40) to (50);
\draw [scaffold, horizontal] (40) to (30);
\draw [scaffold, horizontal] (60) to (50);
\draw [scaffold, horizontal] (01) to (11);
\draw [scaffold, horizontal] (21) to (11);
\draw [scaffold, horizontal] (21) to (31);
\draw [scaffold, horizontal] (41) to (31);
\draw [scaffold, horizontal] (41) to (51);
\draw [scaffold, horizontal] (61) to (51);
\draw [scaffold, horizontal] (02) to (12);
\draw [scaffold, horizontal] (22) to (12);
\draw [scaffold, horizontal] (22) to (32);
\draw [scaffold, horizontal] (42) to (32);
\draw [scaffold, horizontal] (42) to (52);
\draw [scaffold, horizontal] (62) to (52);
\draw [scaffold, horizontal] (62) to (72);
\draw [scaffold, horizontal] (82) to (72);
\draw [scaffold, horizontal] (03) to (13);
\draw [scaffold, horizontal] (23) to (13);
\draw [scaffold, horizontal] (23) to (33);
\draw [scaffold, horizontal] (43) to (33);
\draw [scaffold, horizontal] (43) to (53);
\draw [scaffold, horizontal] (63) to (53);
\draw [scaffold, horizontal] (04) to (14);
\draw [scaffold, horizontal] (24) to (14);
\draw [scaffold, horizontal] (24) to (34);
\draw [scaffold, horizontal] (44) to (34);
\draw [scaffold, horizontal] (44) to (54);
\draw [scaffold, horizontal] (64) to (54);
\draw [scaffold, horizontal] (05) to (15);
\draw [scaffold, horizontal] (25) to (15);
\draw [scaffold, horizontal] (25) to (35);
\draw [scaffold, horizontal] (45) to (35);
\draw [scaffold, horizontal] (45) to (55);
\draw [scaffold, horizontal] (65) to (55);
\draw [scaffold, horizontal] (65) to (75);
\draw [scaffold, horizontal] (85) to (75);
\draw [scaffold, horizontal] (06) to (16);
\draw [scaffold, horizontal] (26) to (16);
\draw [scaffold, horizontal] (26) to (36);
\draw [scaffold, horizontal] (46) to (36);
\draw [scaffold, horizontal] (46) to (56);
\draw [scaffold, horizontal] (66) to (56);
\draw [equality, vertical] (00) to (01);
\draw [scaffold, vertical] (10) to (11);
\draw [equality, vertical] (20) to (21);
\draw [scaffold, vertical] (30) to (31);
\draw [equality, vertical] (40) to (41);
\draw [scaffold, vertical] (50) to (51);
\draw [equality, vertical] (60) to (61);
\draw [equality, vertical] (02) to (01);
\draw [scaffold, vertical] (12) to (11);
\draw [equality, vertical] (22) to (21);
\draw [scaffold, vertical] (32) to (31);
\draw [scaffold, vertical] (52) to (31);
\draw [equality, vertical] (62) to (41);
\draw [scaffold, vertical] (72) to (51);
\draw [equality, vertical] (82) to (61);
\draw [equality, vertical] (02) to (03);
\draw [scaffold, vertical] (12) to (13);
\draw [scaffold, vertical] (32) to (13);
\draw [equality, vertical] (42) to (23);
\draw [scaffold, vertical] (52) to (33);
\draw [equality, vertical] (62) to (43);
\draw [scaffold, vertical] (72) to (53);
\draw [equality, vertical] (82) to (63);
\draw [equality, vertical] (04) to (03);
\draw [scaffold, vertical] (14) to (13);
\draw [equality, vertical] (24) to (23);
\draw [scaffold, vertical] (34) to (33);
\draw [equality, vertical] (44) to (43);
\draw [scaffold, vertical] (54) to (53);
\draw [equality, vertical] (64) to (63);
\draw [equality, vertical] (04) to (05);
\draw [scaffold, vertical] (14) to (15);
\draw [equality, vertical] (24) to (25);
\draw [scaffold, vertical] (34) to (35);
\draw [equality, vertical] (44) to (45);
\draw [equality, vertical] (44) to (65);
\draw [scaffold, vertical] (54) to (75);
\draw [equality, vertical] (64) to (85);
\draw [equality, vertical] (06) to (05);
\draw [scaffold, vertical] (16) to (15);
\draw [equality, vertical] (26) to (25);
\draw [scaffold, vertical] (36) to (35);
\draw [equality, vertical] (46) to (45);
\draw [equality, vertical] (46) to (65);
\draw [scaffold, vertical] (56) to (75);
\draw [equality, vertical] (66) to (85);
\end{tikzpicture}
\end{aligned}
\\[10pt]
\nonumber
\text{\emph{(a) An untyped 0-diagram.}}
&
\text{\emph{(b) An untyped 1-diagram.}}
&
\text{\emph{(c) An untyped 2-diagram.}}
\end{calign}
\caption{\label{fig:untyped}Examples of untyped diagrams.}
\end{figure*}

%% file: section3-contraction.tex
\section{Contraction}

\label{sec:contraction}

\noindent
We define contraction as follows.
\begin{definition}
Given a zigzag in \cat C, we define its \emph{contraction} to be the zigzag of length 1 arising from the colimit in \cat C, if it exists, of its zigzag diagram:
\begin{equation}
\label{eq:zigzag}
\begin{tz}[xscale=0.9, yscale=1]
\node (1) at (0,0) {$r_0$};
\node (2) at (1,1) {$s_0$};
\node (3) at (2,0) {$r_1$};
\node (4) at (3,1) {$s_1$};
\node (5) at (4,0) {\phantom{$s_1$}};
\node at (5,0.5) {$\cdots$};
\node (6) at (6,0) {\phantom{$s_1$}};
\node (7) at (7,1) {$s_{n-1}$};
\node (8) at (8,0) {$r_{n}$};
\draw [->] (1) to node [above left, inner sep=0pt] {$f_0$} (2);
\draw [->] (3) to node [above right, inner sep=0pt] {$b_0$} (2);
\draw [->] (3) to node [above left, inner sep=0pt] {$f_1$} (4);
\draw [->, dashed] (5) to node [above right, inner sep=0pt] {$b_1$} (4);
\draw [->, dashed] (6) to node [above left, inner sep=0pt] {$f_{n-1}$} (7);
\draw [->] (8) to node [above right, inner sep=0pt] {$b_{n-1}$} (7);
\node (C) at (4,3) {$C$};
\draw [->] (2) to node [left, pos=0.52] {$c_0$} (C);
\draw [->] (4) to node [left, pos=0.53] {$c_1$} (C);
\draw [->] (7) to node [right] {$c_{n-1}$} (C);
\node at (4.5,2) {$\cdots$};
\end{tz}
\end{equation}
\end{definition}

\noindent
Given the structure of a zigzag diagram, we of course only need to define the cocone maps for the singular objects. If the colimit exists, then the contraction is defined to be the following zigzag in \cat C of length 1:
\begin{equation}
\begin{aligned}
\begin{tikzpicture}[yscale=0.5]
\node (1) at (0,0) {$r_0$};
\node (2) at (2,0) {$r_n$};
\node (3) at (1,2) {$C$};
\draw [->] (1) to node [above left] {$c_0 \circ f_0$} (3);
\draw [->] (2) to node [above right] {$c_{n-1} \circ b_{n-1}$} (3);
\end{tikzpicture}
\end{aligned}
\end{equation}

If the colimit does not exist, then the contraction is not defined. For our intended application this will frequently be the case, as the categories $Z^n _\cat L$ that we will be working with lack many colimits. We can interpret this as saying that contraction is nontrivial, and not always possible for an $n$\-diagram.

\begin{remark}
In such a zigzag colimit diagram, note that the first and last regular objects $r_0$ and $r_n$, and their associated morphisms $f_0$ and $b_{n-1}$, do not affect the colimit. When it simplifies the narrative to do so, we will ignore them in our formal developments below.
\end{remark}

\subsection{Constructing zigzag colimits}

\noindent
Given a connected diagram in $Z_\cat C$, we build its colimit, or detect that such a colimit does not exist, by the following scheme. This scheme, and its correctness proofs, are the main mathematical contributions of this paper. Note that we do not assume that \cat C itself has any particular colimits; but if \cat C has few colimits, then the same will be true for $Z_\cat C$.

\begin{definition}[Zigzag colimit]\label{def:zigzagcolimit}
For a category \cat C with a terminal object, given a non-empty connected diagram \mbox{$D: J \to Z_\cat C$}, we build its colimit, or fail, according to the following scheme. To fix notation, we write $C$ for the final colimit zigzag that we are trying to construct, and for each $j \in \Ob(\cat J)$, we write  $f^j: S_\cat C  D(j) \to C$ for the corresponding cocone zigzag map.
\begin{enumerate}[label=(\arabic*)]
\item\label{it:colim1}  Build the diagram $\cat J \to[D] Z _\cat C \to[S_\cat C] \Simp$, and obtain its colimit. If no colimit exists, fail.
\item Otherwise, we have a colimit object $c \in \Ob(\Simp)$, and cocone monotone functions $c_j : S_\cat C  D(j) \to c$ for every $j \in \Ob(\cat J)$.
\item We choose the zigzag $C$ to have length $c$, and we choose the monotone functions $(f^j)_\sing = c_j$.
\item\label{it:colim21} We now perform the following subconstruction  for each $k \in [c]$, as follows.
\label{enum:itemlabel1}
\begin{enumerate}[label=(\roman*)]
\item Restrict the diagram $D : \cat J \to Z_\cat C$ to a diagram \mbox{$D_k: \cat J \to Z_\cat C$}, by defining $D_k(j)$ on an object \mbox{$j \in \Ob(\cat J)$} as the restricted zigzag $D(j) _{c_j'(k,k+1)}$, and similarly on morphisms.\footnote{Recall Definition~\ref{def:k} of $(-)': \Simp \to \Simp ^\op_=$.}
\item\label{it:colim2} Build the deconstruction $(D_k)^*: {\cat J}^* \to \cat C$,\footnote{Recall Definition~\ref{def:deconstruction} of the deconstruction procedure.} and obtain its  colimit. If no colimit exists, fail.
\label{enum:itemlabel2}
\item Otherwise, we have a colimit object $p \in \Ob(\cat C)$, and for any $j \in \Ob(\cat J)$ and $i \in D(j)_\sing$, a cocone morphism of type $p_{i}^j : (D_k)^*(j,s_i ^{D(j)}) \to p$.
\item Build a zigzag $C_k$ of height 1 as follows. Choose some $j \in \Ob(\cat J)$ with  $D_k(j)_\sing = m> 0$.\footnote{Such a $j$ must exist, since a colimit in $\Simp$ of empty sets is empty.} Define the forward map as $f^{C_k}_0 = p_{0}^j \circ f^{D_k(j)}_0$, and the backward map as $b_0 ^{C_k} = p ^j _{m-1}\circ b_{m-1}^{D_k(j)}$. Hence obtain $r_0 ^{C_k}$ and $r_1^{C_k}$ as the sources of $f_0 ^{C_k}$ and $b_0 ^{C_k}$ respectively. We set~$s_0 ^{P_k} = p$.
\item For a fixed $j \in \Ob(\cat J)$, build a zigzag map of type $f^{j,k}:D_k(j) \to C_k$ by choosing the monotone map as the unique one of type $D_k(j)_\sing \to{} [1]$, and by choosing the singular morphisms at source singular height $i$ as~$p_i ^j$.
\end{enumerate}
\item Build the colimit zigzag $C$ as the concatenation of the length-1 zigzags $C_k$.\footnote{Recall Definition~\ref{def:concatenate} of concatenation of zigzags and their  maps.}
\item For each value of $j$, build the cocone zigzag map $f^j$ as the concatenation of the zigzag maps $f ^{j,k}$ for $k \in [c]$.
\end{enumerate}
\end{definition}

\noindent
This completes the description of the colimit construction scheme. The correctness proofs follow in Section~\ref{sec:colimitcorrectness}.

\begin{figure}[b!]
\[
\begin{tz}
\node[scale=0.7] at (0,0){
\begin{tikzcd}[row sep=0.5em, column sep=0.5em]
&&&& r_2 \arrow[to=R2RT, equal] \arrow[to= R2RB, equal] \arrow[to=R2B,equal] \arrow[dl] &&&&&&&&&&&&&&&&&&&&&&&&
\\&&&   s_2 \arrow[to=S1'] \arrow[to=S1'', crossing over] &&&&&&&&&&&&&&&&&&&&&&&& |[alias=R2RT]| r_2 \arrow[dl] \arrow[to= R2BRT, equal]
\\&& r_1 \arrow[ur] \arrow[dl] &&&&&&&&&&&&&&&&&&&&&&&&  |[alias=S2'']| s_2''
\\& s_1 \arrow[to=S1'] \arrow[to=S1'',crossing over] &&&&&&&&&&&&&&&&&&&&&&&& |[alias=R2RB]|  r_2\arrow[ur] \arrow[dl]
\\ r_0 \arrow[ur] \arrow[to= R0B,equal]\arrow[to=R0R, equal,crossing over] &&&&&&&&&&&&&&&&&&&&&&&& |[alias = S1'']| s_1''
\\ &&&&&&&&&&&&&&&&&&&&&&& |[alias= R0R]| r_0\arrow[ur] 
\\ \phantom{a}
\\ &&&&&&&&&&&&&&&&&&&&&&&&&&&&
\\ &&& |[alias=R2B]|r_2\arrow[dl] \arrow[to=R2BRB, equal] \arrow[to=R2BRT, equal]  &&&&&&&&&&&&&&&&&&&&&&&&  |[alias=R2BRT]| r_2 \arrow[dl]
\\ && |[alias=S1']| s_1' \arrow[to=S1t] &&&&&&&&&&&&&&&&&&&&&&&&  |[alias=S2t]| \widetilde{s}_2 \arrow[from =S2'',crossing over]
\\ & |[alias = R0B]|r_0 \arrow[ur]   \arrow[to = R0BR, equal] &&&&&&&&&&&&&&&&&&&&&&&&  |[alias=R2BRB]| r_2 \arrow[ur] \arrow[dl] \arrow[from=R2RB, equal, crossing over]
\\ &&&&&&&&&&&&&&&&&&&&&&&& |[alias=S1t]| \widetilde{s}_1 \arrow[from=S1'', crossing over]
\\ &&&&&&&&&&&&&&&&&&&&&&& |[alias = R0BR]| r_0 \arrow[ur] \arrow[from=R0R, crossing over, equal]
\end{tikzcd}
};
\end{tz}
\]
\caption{A pushout in $Z_\cat C$.\label{fig:pushout}}
\end{figure}%
We illustrate this procedure in Figure~\ref{fig:pushout}, which shows the computation of a pushout in $Z_\cat C$. The top-left, bottom-left and top-right zigzags are given, as well as the  maps between them.  The length of the bottom-right zigzag, and its incoming monotone maps, are determined by taking a pushout in $\Simp$. The regular objects of the bottom-right zigzag are completely determined by the incoming maps, and the singular objects are computed as colimits over the `incoming diagrams':
\[
\widetilde{s}_1 = \mathrm{colim}\left(
\begin{tz}
\node[scale=0.75] at (0,0){
\begin{tikzcd}[row sep=0.5em, column sep=0.5em]
&& |[alias=S2]|s_2 \arrow[to=S1'] \\
& r_1 \arrow[ur] \arrow[dl]&&&&&&&& \\
|[alias=S1]|s_1 \arrow[to=S1'] &&&&&&&& |[alias=S1'']|s_1'' \arrow[from=S1,crossing over] \arrow[from=S2, crossing over]\\
 \phantom{a}\\ \phantom{a}\\ 
 & |[alias=S1']|s_1'
\end{tikzcd}
};
\end{tz}
\right)
\hspace{.5cm}
\widetilde{s}_2 = s_2''
\]
The morphisms into the singular objects are given by the obvious morphisms into the colimits.

In the implementation \homotopyio, this colimit construction scheme provides the main recursive algorithm for performing contractions of typed diagrams, as objects of $Z^n _\cat L$. While we do not go into detail regarding the implementation, it is at least worth noting that termination is clear, since colimits in the base category \cat L can be trivially computed, and for a finite diagram, this colimit construction scheme involves only finitely many loops, with all recursion being to strictly lower-dimensional instances.

\subsection{Examples}

\noindent
We already encountered some nontrivial examples of contractions, in Figures~\ref{fig:algebraiccontraction} and~\ref{fig:tanglecontraction}. We give some further examples here. In the online versions of the proofs, you can view the contraction yourself by changing the setting of the ``Slice'' control in the top-right, or perform the contraction  yourself (where possible),  using the mouse to drag one of the vertices vertically towards the centre of the diagram.

\begin{example}
[\href{http://www.cs.bham.ac.uk/~vicaryjo/homotopy.io/lics2019/contraction_example_1.html}{Link to online proof}]
Here we perform a contraction in $Z^2 _\cat L$ of a zigzag of length 2,  containing 2 vertices. In the contracted diagram, these vertices are at the same height.
\begin{equation}
\alignedgraphics{width=2cm}{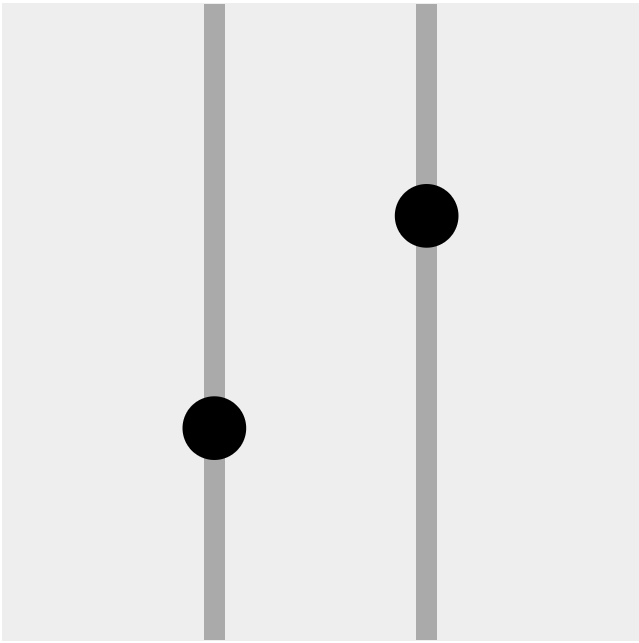}
\quad\leadsto\quad
\alignedgraphics{width=2cm}{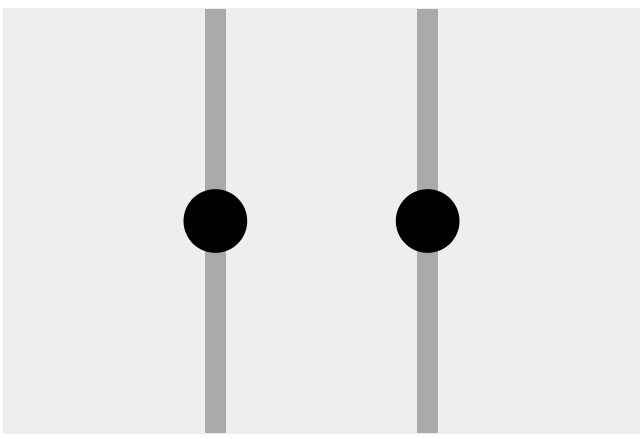}
\end{equation}
\end{example}

\begin{example}
[\href{http://www.cs.bham.ac.uk/~vicaryjo/homotopy.io/lics2019/contraction_example_2.html}{Link to online proof}]
\label{ex:opposingunitcounit}
In this non-example, again in $Z^2_\cat L$, the colimit construction procedure fails at step (1), since the diagram $\cat J \to[D] Z _\cat C \to[S_\cat C] \Simp$ has image  \mbox{$[1] \ot{} [0] \to{} [1]$}, which does not have a colimit:
\begin{equation}
\label{eq:opposingunitcounit}
\alignedgraphics{height=2cm}{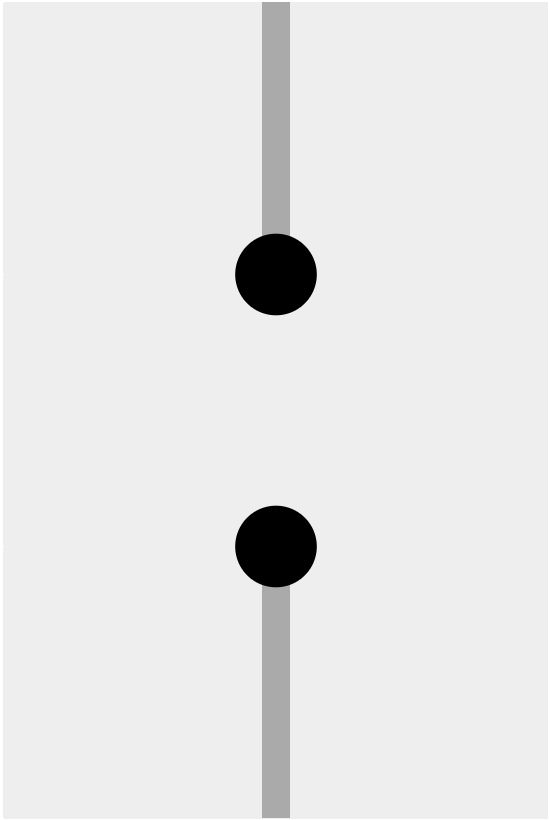}
\end{equation}
To understand why this contraction does not exist, consider that, if it could be constructed, the resulting unique singular height would have to contain 2 vertices, with one to the left of the other, as follows:
\begin{calign}
\label{eq:symmetrybroken}
\alignedgraphics{width=2cm}{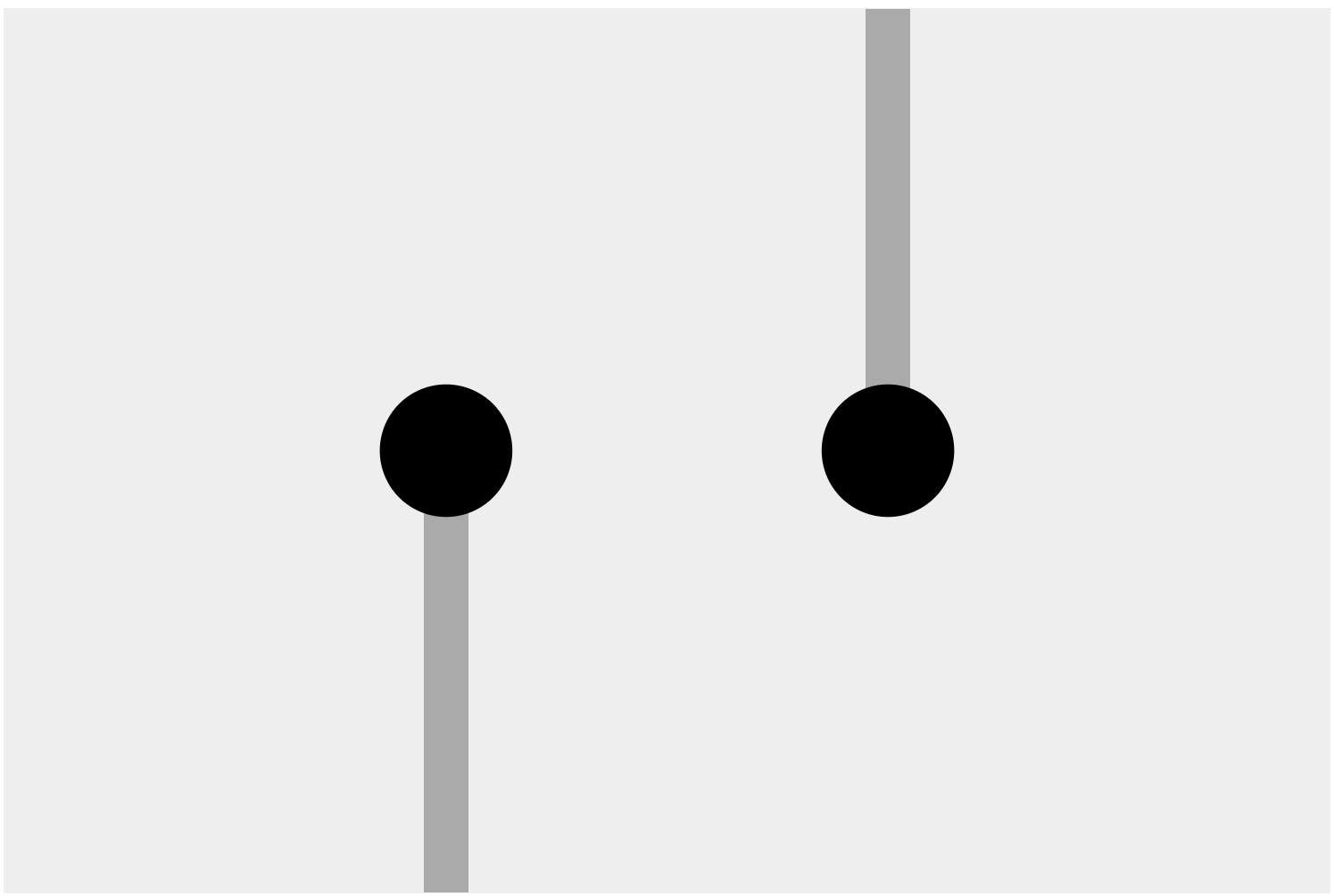}
&
\alignedgraphics{width=2cm}{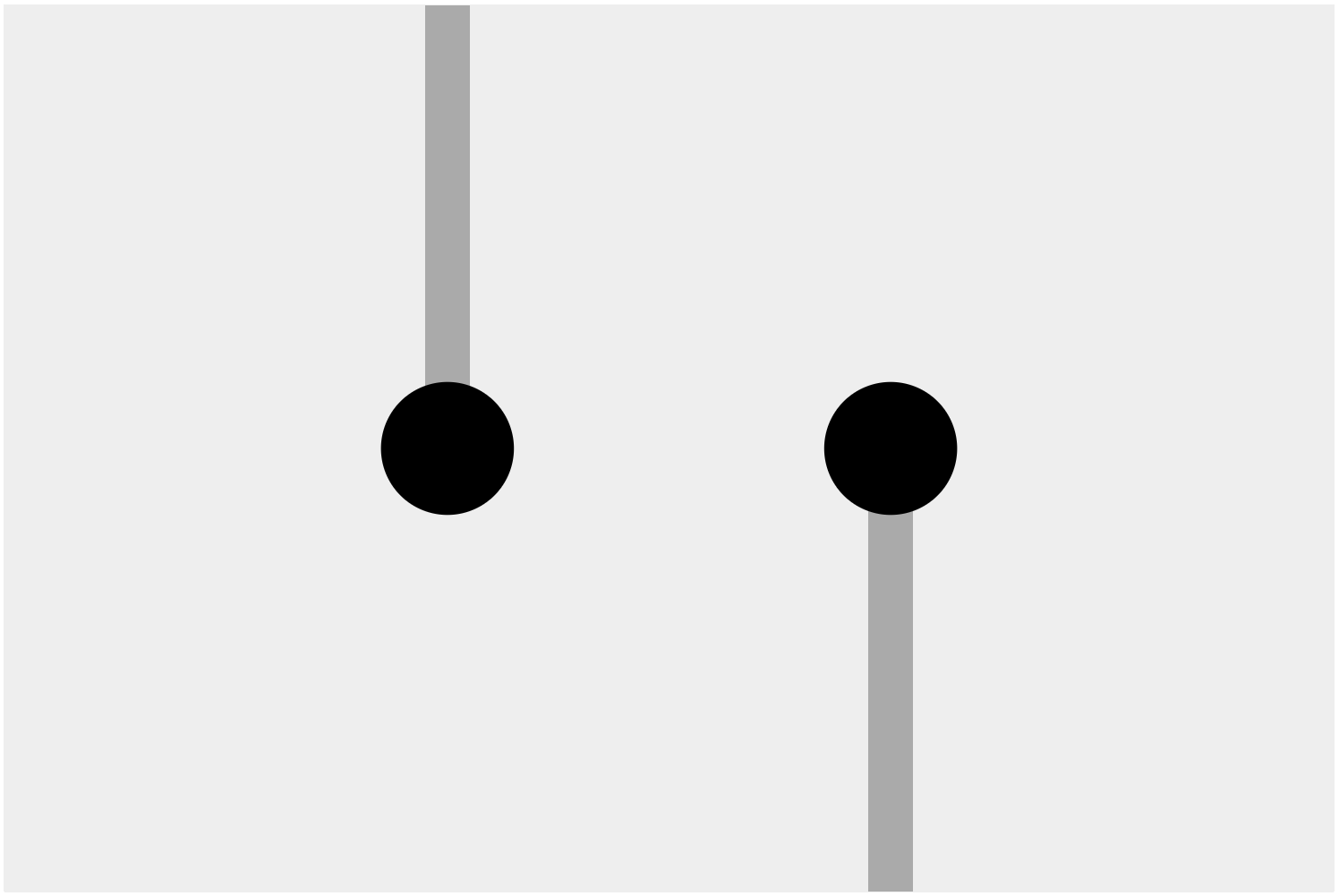}
\end{calign}
However, the colimit construction algorithm has no way to ``break the symmetry'', and cannot proceed. The implementation \homotopyio uses some additional techniques (see Section~\ref{sec:additionalmethods}) which allow us to break the symmetry here; in the online proof, we apply these techniques by dragging the upper vertex of~\eqref{eq:opposingunitcounit} in a south-east or south-west direction, to produce the two images given in~\eqref{eq:symmetrybroken}.
\end{example}

\begin{example}
[\href{http://www.cs.bham.ac.uk/~vicaryjo/homotopy.io/lics2019/contraction_example_3.html}{Link to online proof}]
If we modify the previous example by putting a wire in between the vertices, the diagram will now contract successfully:
\begin{equation}
\alignedgraphics{height=2cm}{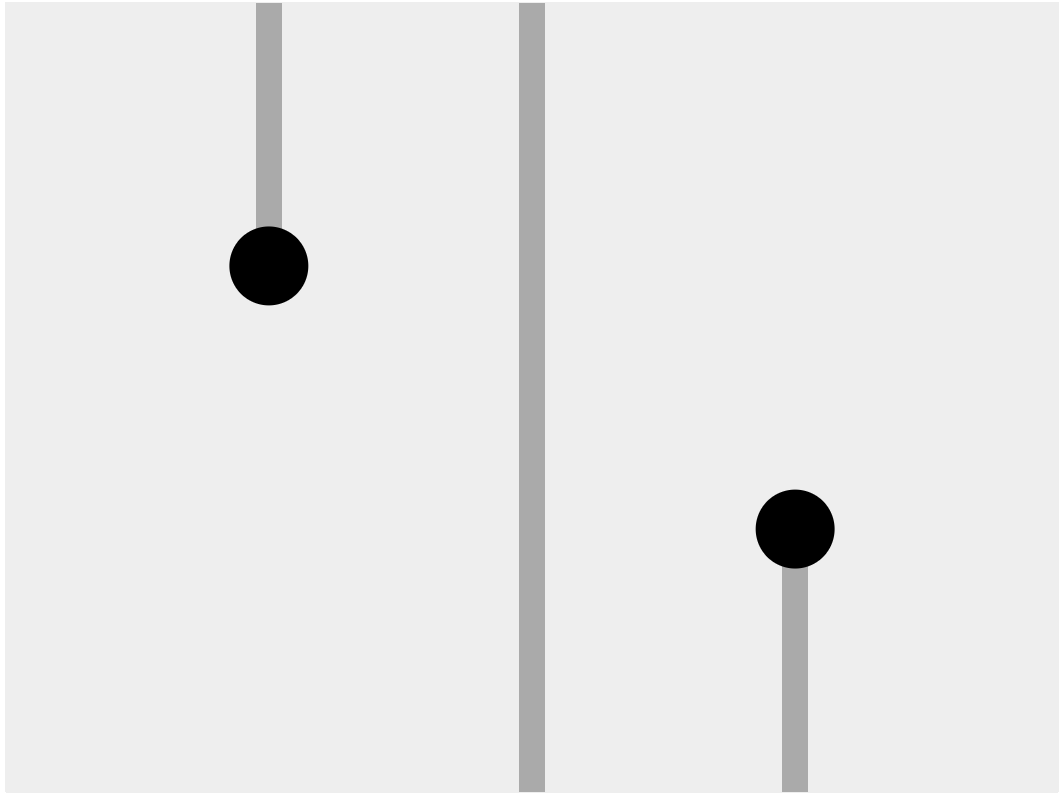}
\quad\leadsto\quad
\alignedgraphics{height=1.25cm}{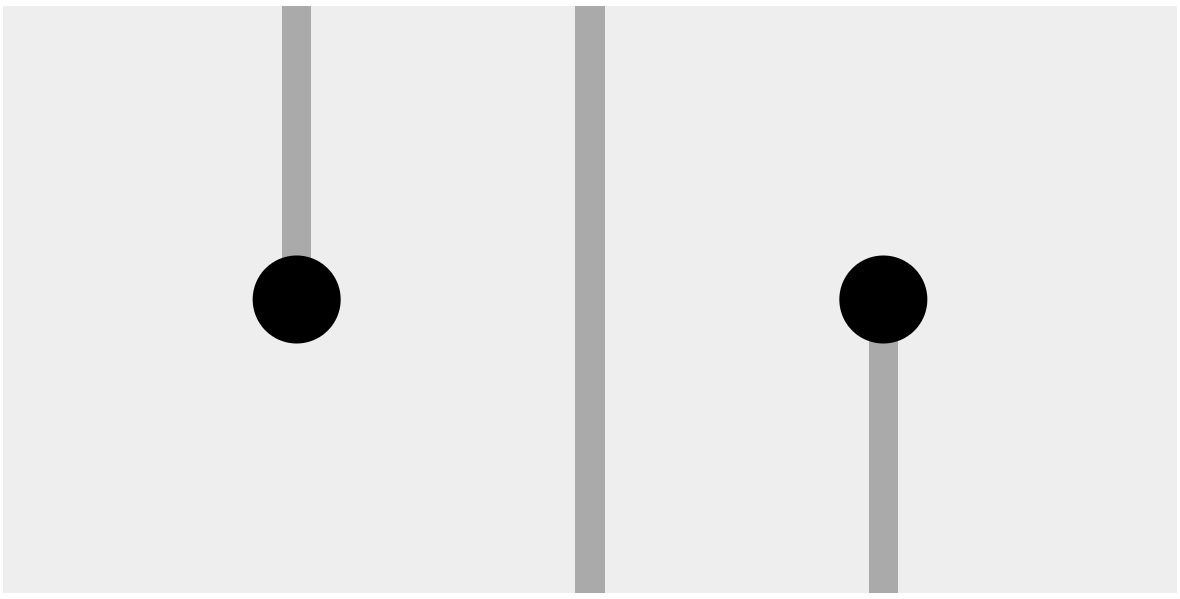}
\end{equation}
This is because the colimit diagram  in $\Delta$ now has the image $[2] \ot[0 \mapsto 1] [1] \to[0 \mapsto 0] [2]$, which does have a colimit.
\end{example}

Suppose that we are taking the contraction of a typed $n$\-diagram $D$---that is, an object of $Z^n _{\cat L}$---which is well-typed with respect to some signature $\Sigma$ (see Section~\ref{sec:untypedandtyped} for a brief discussion of type checking.) Even if the contraction of $D$ exists, yielding a new object $D'$ of $Z^n _\cat L$, it does not follow that $D'$ will again be well-typed with respect to $\Sigma$; the entire contraction $D \to D'$ must be passed through the type checker to verify this. We show such an example here.

\begin{example}
[\href{http://www.cs.bham.ac.uk/~vicaryjo/homotopy.io/lics2019/contraction_example_4.html}{Link to online proof}]
In this example, again in $Z_\cat L^2$, we contract a zigzag of length 2, as follows:
\begin{equation}
\alignedgraphics{width=1.3333cm}{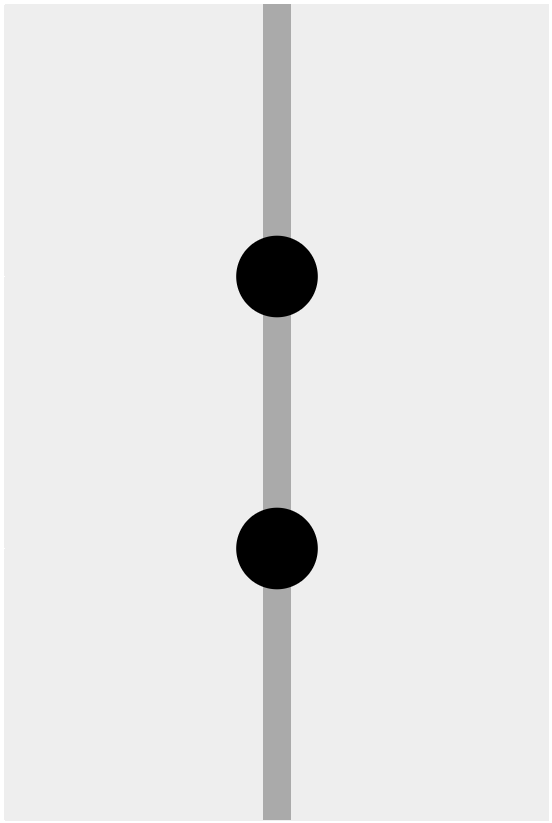}
\quad\leadsto\quad
\alignedgraphics{width=1.3333cm}{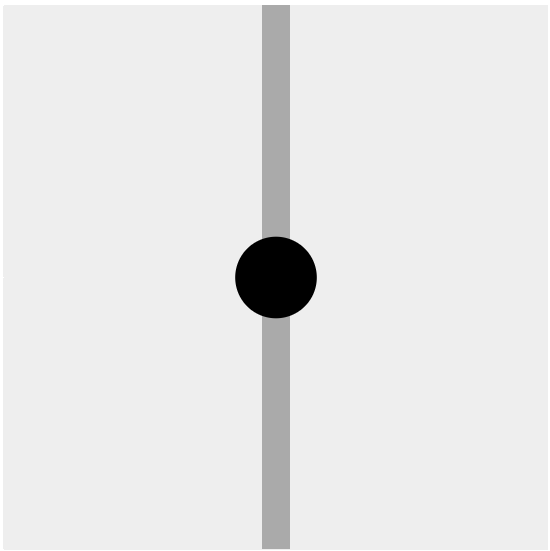}
\end{equation}
Here we ``fuse'' two endomorphisms on a wire into a single endomorphism, with the colimit construction procedure successfully returning the right-hand diagram. Both of these diagrams type check, but the contraction process as a whole does not, because homotopies may only ``move'' parts of the diagram around, not change the structure of individual labels. As a result, in the online proof, clicking and dragging either of the two vertices will have no effect, as the contraction above will be silently blocked by the type checker. This shows the way that contraction and type checking interact in the implementation.
\end{example}

\subsection{Correctness}
\label{sec:colimitcorrectness}

\begin{theorem}
\label{thm:zigzagcolimits}
Let $\cC$ be a category with a terminal object and let $D: \cJ \to Z_{\cC}$ be a non-empty connected diagram. Then, $D$ has a colimit if and only if the procedure in Definition~\ref{def:zigzagcolimit} succeeds (that is, if the colimits in step~\ref{it:colim1} and~(4.ii) exist), and the procedure constructs it. 
\end{theorem}

\begin{remark}
Since the category $Z_{\cC}$ is a disjoint union of local zigzag categories $Z_{\cC}(a,b)$ for objects $a,b $ in $\cC$ (see Definition~\ref{def:local}), and since Theorem~\ref{thm:zigzagcolimits} applies to connected non-empty diagrams, it is also true as stated for the categories $Z_{\cC}(a,b)$ replacing $Z_{\cC}$. Moreover, note that if $\cC$ is a category with a terminal object $*$, then $Z_{\cC}(a,b)$ has a terminal object (namely, the zigzag $a\to * \ot b$). Since all categories of diagrams may be obtained as iterated local zigzag categories, Theorem~\ref{thm:zigzagcolimits} holds for such categories. 
\end{remark}

We prove Theorem~\ref{thm:zigzagcolimits} in two steps. First, we show that if the colimits in step~\ref{it:colim1} and~(4.ii) exist, then the constructed cocone is indeed colimiting. Then, we prove that if a colimit of a diagram $D:J \to Z_{\cC}$ exist, then the colimits in step~\ref{it:colim1} and~(4.ii) must also exist.

\subsection{The procedure correctly computes colimits}

\noindent
We prove the first part of Theorem~\ref{thm:zigzagcolimits}: If the colimits in step~\ref{it:colim1} and~(4.ii) exist, then the constructed cocone is indeed a colimiting cocone of the diagram $D:J\to Z_{\cC}$.

The proof boils down to the following categorical fact: given a (Grothendieck) opfibration $F:\cA\to \cB$, then colimits in $\cA$ can be computed in terms of colimits in $\cB$ and in the fibre categories $F^{-1}(b)$ for objects $b\in \cB$.

\vspace{4pt}
\paragraph{Opfibrations and colimits}{}
We recall the following terminology.
Given a functor $F:\cA\to\cB$, a morphism $\phi: a \to a'$ in $\cA$ is called \emph{opcartesian} if for any morphism $\psi: a\to a''$ in $\cA$ and $g:F(a') \to F(a'')$ in $\cB$ such that $g\circ F(\phi) = F(\psi)$, there exists a unique $\chi:a' \to a''$ such that $\chi\circ \phi =\psi $ and $F(\chi)=g$. A functor $F:\cA\to \cB$ is an \emph{opfibration} if for any $a\in \cA$ and $h:F(a) \to b$ in $\cB$, there is a opcartesian morphism $\phi:a \to a'$ with $F(\phi) = h$. For an opfibration $F:\cA\to \cB$ and an object $b\in \cB$, the \emph{fibre category} $F^{-1}(b)$ is the subcategory of $\cA$ with objects and morphisms mapping to $b$ and $\id_b$, respectively. Given a morphism $\sigma: b\to b'$ in $\cB$, the \emph{base change functor} $\sigma_*: F^{-1}(b) \to F^{-1}(b')$ maps an object $a$ in the fibre over $b$ to the codomain of the opcartesian morphism lifting $\sigma:Fa\to b'$ and a morphism $f:a\to a'$ over $\id_b$ to the morphism $\sigma_*a \to \sigma_*a'$ obtained from opcartesianity of the lift of $\sigma:Fa\to b'$. 

We recall the following basic fact about opfibrations.

\ignore{Recall that a functor $F$ \emph{reflects colimits} if a cocone is colimiting whenever\DR{maybe don't say `whenever' here. I mean: given cocone C, such that $FC$ is colimiting. Then C colimiting}  its image under $F$ is.
\begin{proposition}
Let $F:\cA\to \cB$ be an opfibration and let $D:\cJ \to \cA$ be a functor such that $FD$ has a colimit. Suppose that the induced functor $\cJ \to F^{-1}\left(\colim FD \right)$ has a colimiting cocone $\mu_j: g_j \to g$. Then, $\mu_j\circ \phi_j:D_j \to g$ is a colimiting cocone of $D$. 
\end{proposition}
 }
\begin{proposition}
\label{prop:opfibration}
Let $F:\cA \to \cB$ be an opfibration and let $D:\cJ \to \cA$ be a diagram such that $F D$ has a colimit. If all fibres have $\cJ$-colimits and the base change functor $\sigma_*:F^{-1}(b) \to F^{-1}(b')$ preserves them for all $\sigma:b\to b'$ in $\cB$, then $D$ has a colimit and $F$ preserves it.
\end{proposition}

\noindent
This proposition is proven later as Proposition~\ref{prop:opfibrationapp}.

Explicitly, we can compute this colimit in terms of the colimit of $FD: \cJ \to \cB$ as follows: Lift the universal cocone morphisms $\lambda^j: FD_j \to \colim FD$ to opcartesian morphisms $\phi^j: D_j \to \lambda^j_*(D_j)$, where $F\left(\lambda^j_*(D_j)\right) = \colim FD$. Opcartesianity of $\phi^j$ gives rise to morphisms $\lambda_{\sigma}: \lambda^j_*(D_j) \to \lambda^{j'}_*(D_{j'})$ for $\sigma : j \to j'$ in $\cJ$, making this into a diagram $\cJ\to F^{-1}(\colim FD)$. A colimiting cocone $\mu^j: \lambda^j_*(D_j) \to X$ of this diagram $\cJ \to F^{-1}(\colim FD)$ induces a colimiting cocone $\mu^j \circ \phi^j: D_j \to X$ of $D$.

\vspace{5pt}
\paragraph{$S_{\cC}$ is an opfibration for cocomplete $\cC$}
 Given a zigzag $Z$ (drawn on the left) with a chosen regular object (here labelled $r$), we define a new zigzag $\widetilde{Z}$ (drawn on the right) in which the regular object is `expanded' into two regular objects, and a morphism of zigzags $Z\to \widetilde{Z}$ as follows:
 \[\begin{tz}
 \node[scale=0.75, yscale=.8] at (0,0){
 \begin{tikzcd}[row sep=1em, column sep=1em]
&&&&&&&&&&&&  \svdots  \arrow[d]\\
 \svdots\arrow[d]  &&&&&&&&&&&& |[alias=S1R]| s_1\\
 s_1 \arrow[to=S1R, "\id_{s_1}"] &&&&&&&&&&&&  |[alias=RRT]|r\vphantom{_1}\arrow[u] \arrow[d, "\id_r"]\\
 r\vphantom{_1}\arrow[u] \arrow[d] \arrow[to=RRT, equal] \arrow[to=RRB, equal]  &&&&&&&&&&&& |[alias=RRM]| r\vphantom{_1}\\
 s_2 \arrow[to=S2R, "\id_{s_2}" ']&&&&&&&&&&&& |[alias=RRB]| r\vphantom{_1}  \arrow[d]\arrow[u, "\id_r" ']\\
 \svdots \arrow[u] &&&&&&&&&&&& |[alias=S2R]| s_2\\
 &&&&&&&&&&&&  \svdots  \arrow[u]
 \end{tikzcd}
 };
 \end{tz}
 \]

If $\cC$ is cocomplete, and $Z$ is a zigzag with a chosen pair of adjacent singular objects (here labelled $s_1$ and $s_2$ on the left), we define a new zigzag $\widetilde{Z}$ in which the singular objects are `collapsed' into a single singular object, given by the pushout of $s_1$ and $s_2$ over the intermediate regular object, and a morphism $Z\to \widetilde{Z}$:
\[\begin{tz}
 \node[scale=0.75, yscale=.8] at (0,0){
 \begin{tikzcd}[row sep=1em, column sep=1em]
 \svdots &&&&&&&&&&&&\\
 r_1 \arrow[to=R1R, equal] \arrow[u] \arrow[d]&&&&&&&&&&&&  \svdots\\
 s_1 \arrow[to=colim]  &&&&&&&&&&&&  |[alias=R1R]|r_1\arrow[u]\arrow[d]\\
 r\arrow[u] \arrow[d]  &&&&&&&&&&&& |[alias=colim]| s_1 \sqcup_{r} s_2\\
 s_2 \arrow[to=colim]&&&&&&&&&&&& |[alias=R2R]| r_2 \arrow[d]\arrow[u]\\
r_2 \arrow[u] \arrow[d] \arrow[to=R2R, equal] &&&&&&&&&&&&\svdots \\
\svdots &&&&&&&&&&&&  
 \end{tikzcd}
 };
 \end{tz}
\]
Given a zigzag $Z$ and a monotone map $h:Z_{\sing} \to I$ into some finite totally ordered set $I$, we iterate these operations to produce a zigzag $\widetilde{Z}$ of length $|I|$ and a morphism of zigzags $\widetilde{h}:Z\to \widetilde{Z}$ with underlying monotone map $\widetilde{h}_{\sing} = h$ as illustrated in the following example lifting the constant monotone map \mbox{$\{1<2<3\}\to \{1<2<3\},~x\mapsto 1$}:
\[\begin{tz}
 \node[scale=0.75, yscale=.8] at (0,0){
 \begin{tikzcd}[row sep=1em, column sep=1em]
 r_0\arrow[d] \arrow[to=R0R, equal]                                     
 &&&&&&&&&&&& |[alias=R0R]| r_0\arrow[d]         \\
 s_1 \arrow[to=col] &&&&&&&&&&&& |[alias=col]| s_1\sqcup_{r_1} s_2 \sqcup_{r_2} s_3\\
 r_1 \arrow[u] \arrow[d]                
 &&&&&&&&&&&&   
|[alias=R3T]| r_3\arrow[d] \arrow[u]\\
 s_2 \arrow[to=col]
 &&&&&&&&&&&&   
r_3\\
 r_2\arrow[u]\arrow[d]  
 &&&&&&&&&&&&   
 |[alias=R3M]|r_3\arrow[u]\arrow[d]\\
 s_3 \arrow[to=col]
 &&&&&&&&&&&&
 r_3\\
 r_3 \arrow[u]\arrow[to=R3B, equal] \arrow[to=R3M, equal]\arrow[to=R3T, equal]
 &&&&&&&&&&&&
 |[alias=R3B]| r_3 \arrow[u]
 \end{tikzcd}
 };
 \end{tz}
 \]
Here, the left zigzag and the underlying monotone map are given; the right zigzag and the map of zigzags are produced by `expanding' and `collapsing'. This ability to `lift' monotone maps $h:Z_{\sing}\to I$ to maps of zigzags leads to the following proposition.

\begin{proposition} If $\cC$ is cocomplete, then the singular monotone functor $S_\cC: Z_{\cC} \to \Simp$ is an opfibration.
\end{proposition}
\begin{proof}
Given a zigzag $Z$ and a monotone map $h: Z_{\sing} \to I$ into some totally ordered set $I$, we lift $h$ to a map of zigzags $\widetilde{h}: Z \to \widetilde{Z}$ obtained by `collapsing' and `expanding' $Z$, as described above.
The fact that $\widetilde{h}$ is opcartesian corresponds precisely to the universal property of the colimits in the collapse operation.
\ignore{
Define $a'$ as the iterated cospan with $|I|$ singular heights as follows:
The labels of the regular heights are uniquely determined by the map $h:\sing(a) \to I$ and the regular height labels of $a$. Given a singular height $i \in I$ in the image of $h$, label the singular height by the colimit of the diagram $a_{h^{-1}(i)}$ \ignore{in the preimage $h^{-1}(i) \subseteq \sing(a)$ and take the colimit \DRcomm{of the diagram in $\cC$ between the top and bottom singular height of $f^{-1}(k)$ including all regular heights in-between}. } (which can be computed in terms of iterated pushouts).\DR{colimit in over-category} If $i$ is not in the image of $h$, then the regular heights adjacent to $i$ are labelled by the same object $x\in \cC$; in this case, we label the singular height $i$ by $x$ with identity incoming morphisms. 
}
\end{proof}

For a finite totally ordered set $I$, the fibre category $S_{\cC}^{-1}(I)$ is the category of zigzags of length $|I|$ with morphisms the maps of zigzags whose underlying monotone map is the identity. Explicitly, this category is the disjoint union
\[\bigsqcup_{r_i\in \mathrm{ob}\cC \text{ for }i \in I} ~~~~~\bigtimes_{0 \leq i \leq |I|-1} (r_i,r_{i+1})/\cC,
\]
where $(a, b)/\cC$ denotes the over-category whose objects are pairs of morphisms $(a\to x, b\to x)$ and morphisms $(a\to x, b\to x)\to (a\to y, b\to y)$ are morphisms $x\to y$ making the obvious triangles commute.

Given a monotone map $\lambda:I \to J$, the induced base change functor $S_\cC^{-1}(I) \to S_\cC^{-1}(J)$ maps a zigzag of length $|I|$ to a zigzag of length $|J|$ by expanding and collapsing according to the monotone map $\lambda$.

\begin{corollary}\label{cor:cocomplete} Let $\cC$  be cocomplete and let $D:\cJ \to Z_{\cC}$ be a connected, non-empty diagram such that $S_{\cC} D$ has a colimit. Then, $D$ has a colimit $C$, which is preserved by $S_{\cC}$, and which can be explicitly constructed as follows:
\begin{enumerate}
\item Construct a colimit $C_{\sing}$ of $S_{\cC} D$ with colimiting cocone $f^j_{\sing}:D^j_{\sing}\to C_{\sing}$.
\item For every $j \in \cJ$, `expand' and `collapse' the zigzag $D^j$ to a zigzag $\widetilde{D}^j$ of length $|C_{\sing}|$ according to the monotone map $f^j_{\sing}:D^j_{\sing} \to C_{\sing}$. This gives rise to a diagram $\widetilde{D}:\cJ \to Z_{\cC}$ in which every zigzag has the same length and every morphism of zigzags has underlying identity monotone map.
\item \label{it:colimcor} For every singular height $i \in C_{\sing}$, let $s_i$ be the colimit in $\cC$ over the diagram $\widetilde{D}|_i: \cJ \to \cC$ obtained by restricting the diagram $\widetilde{D}: \cJ \to Z_{\cC}$ to the singular objects at height $i$ and the morphisms between them (recall that all maps of zigzags in the image of $\widetilde{D}$ have underlying identity monotone map). 
\item For every regular height $i \in (C_{\sing})' = C_{\reg}$, define the regular object $r_i$ to be equal to the regular object of $\widetilde{D}^j$ at height $i$ for some (and hence any) $j\in \cJ$.
\item Define the forward and backward morphisms of $C$ and the singular morphisms of $f^j:D^j \to C$ as the obvious morphisms into the colimits $s_i$.
\end{enumerate}
\end{corollary}
\begin{proof} 
The fibre $S_{\cC}^{-1}(I)$ has all connected colimits since connected colimits in over-categories can be constructed as colimits in the original category $\cC$ in the obvious way. The base change functors can be factored into functors expanding a single regular object or collapsing a pair of adjacent singular objects. Explicitly, the corresponding base change functors are of the form 
\[\cdots \times \id_{(r_{i-1},r_i)/\cC}\times  r_i \times \id_{(r_{i}, r_{i+1})/\cC} \times \cdots, \]
where $r_i: * \to (r_i, r_i)/\cC$ picks out the object $r_i\to[\id_{r_i}] r_i \ot[\id_{r_i}] r_i$, and 
\[ \cdots \times \id_{(r_{i-2}, r_{i-1})/\cC} \times (-\sqcup_{r_i}-) \times \id_{(r_{i+1}, r_{i+2})/\cC} \times \cdots,
\]
where $-\sqcup_{r_i}-: (r_{i-1}, r_i)/\cC \times (r_i, r_{i+1})/\cC\to (r_{i-1}, r_{i+1})/\cC$ takes the pushout of the inner span. It is clear that both functors preserve connected, non-empty colimits.

It therefore follows from Proposition~\ref{prop:opfibration} that $D$ has a colimit which is preserved by $S_{\cC}$ and is constructed as described.
\end{proof}

\ignore{
Given a diagram $D:\cJ \to \IC(\cC)$, compute the colimit $I:=\colim (\sing\circ  D: \cJ \to \Simp_a)$ and consider the fiber $\sing^{-1}(I)$, that is, the category of iterated cospans of height $|I|$ and with underlying identity monotone maps. Since $\sing$ is an opfibration, the diagram $D:\cJ \to \IC(\cC)$ induces a diagram $\cJ\to \sing^{-1}(I)$. Explicitly, this induced diagram can be obtained from the original diagram $D$ by `collapsing' and `padding' the iterated cospans $D_j$ (taking the colimits of subsets of singular heights or adding new singular heights via identities), until all the underlying monotone maps are straigthened into identities.

Lastly, for every singular height in $I$ we take the colimit 
}

\paragraph{Colimits in $S_{\cC}$ if $\cC$ is not cocomplete}
Categories of typed or untyped diagrams---such as the category $\Simp = Z_{\cat 1}$, or iterated zigzag categories on $\Simp$---are far from cocomplete. In particular, Corollary~\ref{cor:cocomplete} does not hold in this setting.

Recall that we have `collapsed' singular objects by taking a colimit in $\cC$, and have later again taken colimits in $\cC$ to compute the colimit of the diagram in the fiber. In other words, we have computed the colimit of $\cJ \to Z_{\cC}$ by first computing the colimit $C_{\sing}$ in $\Simp$ and then, for every $i\in C_{\sing}$, taking several consecutive colimits in $\cC$. If $\cC$ is not cocomplete, it is possible that some of these intermediate colimits do not exist, even if the overall colimit does exist. We can avoid this issue by only `formally' taking intermediate colimits. This can be formalized by passing to the free completion of $\cC$, as follows.

Let $y: \cC \to \widehat{\cC} := [ \cC, \Set]^{\op}$ denote the `dual' Yoneda embedding of $\cC$. The functor $y$ has the convenient property that it preserves and reflects all colimits; in particular, a diagram $D:\cJ \to \cC$ has a colimit if and only if the colimit of the diagram $y D: \cJ \to \widehat{\cC}$ is representable (that is, is in the essential image of $y$.) Moreover, $y$ gives rise to a fully faithful functor $Z_y: Z_{\cC}\to Z_{\widehat{\cC}}$. 

\ignore{We may therefore use Proposition~\ref{} to compute colimits in $\IC(\cC)$ for general categories $\cC$: Let $D:\cJ \to \IC(\cC)$ be a connected diagram and suppose that $\sing \circ D = \sing \circ \IC(y) : \cJ \to \Simp_a$ has a colimit $I$ with monotone cocone  $\lambda_j : \sing D_j \to I$. By collapsing and padding with identities, replace all iterated cospans in the diagram $\IC(y) \circ D$ by iterated cospans of height $|I|$ turning all underlying monotone maps into identities. Then, for every fixed height $i\in I$ compute the colimits of the diagrams $\cJ\to \widehat{\cC}$ at fixed height $i$. Overall, this results in an interated cospan of height $I$ whose regular heights are labelled by objects of $\cC$ and whose singular heights are labelled by objects of $\widehat{\cC}$. Each of these singular height labels is obtained as the result of two successive colimits --- first collapsing all singular and regular heights in the preimage $\lambda_{j}^{-1}(i)$ and then taking the colimit over all these collapsed diagrams --- ; they are representable if and only if the corresponding total colimit exists in $\cC$. 
}

\begin{proposition}\label{prop:thmpart1}
Let $D:\cJ \to Z_{\cC}$ be a connected, non-empty diagram such that the colimits in step~\ref{it:colim1} and~(4.ii) of Definition~\ref{def:zigzagcolimit} exist. Then the cocone constructed in Definition~\ref{def:zigzagcolimit} is colimiting.
\end{proposition}
\begin{proof} It follows from the existence of the colimit in step~\ref{it:colim1} and Corollary~\ref{cor:cocomplete} that the composite $\cJ \to Z_{\cC} \hookrightarrow Z_{\widehat{\cC}}$ has a colimiting cocone, constructed as in Corollary~\ref{cor:cocomplete}. The existence of the colimits in step~(4.ii) of Theorem~\ref{thm:zigzagcolimits} imply that the singular objects of the constructed zigzag (constructed in step~\ref{it:colimcor} of Corollary~\ref{cor:cocomplete}) are representable. Hence, the constructed cocone is in the image of the fully faithful $Z_y: Z_{\cC} \to Z_{\widehat{\cC}}$, and is therefore a colimit of $\cJ \to  Z_{\cC}$.  
\end{proof}

\subsection{The procedure detects all colimits}

\noindent
We now prove the second part of Theorem~\ref{thm:zigzagcolimits}: if a connected, non-empty diagram $D:\cJ \to Z_{\cC}$ has a colimit, then the colimits in step~\ref{it:colim1} and~(4.ii) of Definition~\ref{def:zigzagcolimit} exist.

\begin{proposition}
\label{prop:singpreserves} Let $\cC$ be a category with a terminal object. The functor $S_{\cC}: Z_{\cC} \to \Simp$ preserves connected colimits.
\end{proposition}
\begin{proof} Given a set $X$, we define $\Simp_=(X)$ as the following generalization of the category $\Simp_=$ from Definition~\ref{def:deltaequal}: its objects are pairs $(O,f)$ of a non-empty totally ordered set $O$ and a function $f:O \to X$, and its morphisms $(O, f) \to (O', f')$ are regular monotone maps $\rho:O \to O'$ such that $f'\circ \rho = f$. Note that $\Simp_=(X)$ is the comma category $F/X$, where $F: \Simp_= \to \Set$ is the forgetful functor.

\def\L{\mathbb{L}}
\def\R{\mathbb{R}}
 The regular monotone functor $R_{\cC}: Z_{\cC} \to \left(\Simp_=\right)^\op$ factors through a functor $\L: Z_{\cC} \to \left(\Simp_=(\ob\cC)\right)^{\op} $ mapping a zigzag $Z$ to its totally ordered set of regular objects $Z_{\reg}$ together with the function $Z_{\reg}\to \ob \cC, ~i \mapsto r_i$. We construct a right adjoint $\R: \left( \Simp_=(\ob \cC)\right)^{op} \to Z_{\cC}$ as follows. The functor $\R$ maps an object $(O,f)$ to the zigzag of length $|O|$ with regular objects determined by $f$ and with singular objects given by the terminal object of $\cC$. It maps a morphism $\lambda: (O,f) \to (O', f')$ to the unique morphism of zigzag with underlying regular monotone map $\lambda$. The natural transformation \[\Hom_{Z_{\cC}} (Z, \R(O,f)) \to  \Hom_{\Simp_=(\ob \cC)}((O,f), \L A)  \]
 mapping a map of zigzags $Z\to \R(O,f)$ to its underlying regular map is a natural isomorphism. Hence, $\R$ is right adjoint to $\L$ and in particular, $\L: Z_{\cC} \to \left(\Simp_=(\ob \cC)\right)^{\op}$ preserves colimits.
 
 Therefore, to show that $S_{\cC}: Z_{\cC} \to \Simp$ preserves connected colimits, it suffices to show that the composite functor \[\left(\Simp_=(\ob \cC)\right)^{\op} \to \left(\Simp_=\right)^{\op} \to[\left(-\right)'] \Simp\] preserves connected colimits. Since $\left(-\right)'$ is an equivalence, it suffices to show that $\Simp_=(\ob \cC) \to \Simp_=$ preserves connected limits. Since $\Simp_=(\ob \cC) = F/X$ is a comma category, this follows from Proposition~\ref{prop:commalimit}.
\end{proof}

\begin{corollary}\label{cor:ICpreserve}
Let $\cC$ be a category with a terminal object. The functor $Z_{\cC} \to Z_{\widehat{\cC}}$ preserves connected, non-empty colimits.
\end{corollary}
\begin{proof} Let $D: \cJ \to Z_{\cC}$ be a connected non-empty diagram, and let $a^j:D^j \to C$ be a colimiting cocone. By Proposition~\ref{prop:singpreserves}, the cocone $c^j_{\sing}:D^j_{\sing} \to C_{\sing}$ is colimiting in~$\Simp$. By Corollary~\ref{cor:cocomplete}, the composite $\cJ \to Z_{\cC} \to Z_{\widehat{\cC}}$ has a colimiting cocone $\widehat{c}^j:D^j \to \widehat{C}$. In particular, there is a morphism of cocones $\mu:\widehat{C} \to C$ in $Z_{\widehat{C}}$. In the following, we show that $\mu$ is an isomorphism.

  Since $\mu_{\sing}: \widehat{C}_{\sing} \to C_{\sing}$ is a morphism of cocones $\big(D^j_{\sing} \to \widehat{C}_{\sing}\big) \to \big(D^j_{\sing}\to C_{\sing}\big)$, and since both cocones are colimiting, it follows that $\mu_{\sing}$ is the identity and that $c^j_{\sing}= \widehat{c}^j_{\sing}$. Denote the regular objects of $C$ by $r_0, \ldots, r_n$. Since there is a morphism $\mu: \widehat{C} \to C$ with $\mu_{\sing} = \id$, it follows that the regular objects of $\widehat{C}$ are also $r_0, \ldots, r_n$. In particular, the morphism $\mu$ can be understood as a morphism in the category $\widehat{\cE}:=\bigtimes_{i} (r_i, r_{i+1})/\widehat{\cC}$. Denoting $\cE := \bigtimes_{i} (r_i, r_{i+1})/\cC$, we observe that the obvious functor $\widehat{\cE} \to[] [\cE, \Set]^{op}$ is an equivalence. 
  
Let $E$ be an object of $\cE$---or equivalently, a zigzag with regular objects $r_0,\ldots, r_n$ and singular objects in $\cC$---and let $\lambda: \widehat{C} \to E$ be a morphism in $\widehat{\cE}$. Then, the composite $\lambda \circ \widehat{c}^j:D^j \to E$ is a cocone of $\cJ \to Z_{\cC}$. In particular, there is a unique morphism $\phi: C\to E$ in $Z_{\cC}$ such that $ \lambda \circ \widehat{c}^j = \phi\circ c^j  = \phi \circ \mu \circ \widehat{c}^j$, or equivalently such that $\lambda = \phi \circ \mu$. Applying the singular monotone functor $S_{\cC}$ and using that $\lambda_{\sing} = \mu_{\sing} = \id$, it follows that $\phi_{\sing} =\id$. We can therefore summarize the preceding paragraph as follows: given an object $E$ of $\cE$ and a morphism $\lambda: \widehat{C} \to E$ in $\widehat{\cE}$, there is a unique morphism $\phi:C\to E$ in $\cE$ such that $\lambda =  \phi \circ \mu$. Since $\widehat{\cE}$ is equivalent to the free completion $[\cE, \Set]^{\op}$ of $\cE$, this means that $\widehat{C}$ is in the essential image of $\cE\to \widehat{\cE}$ and hence isomorphic to $C$.\looseness=-2
\end{proof}

\begin{corollary}\label{cor:thmpart2} Let $\cC$ be a category with a terminal object and let $\cJ \to Z_{\cC}$ be a connected, non-empty diagram admitting a colimit. Then, the colimits in step~\ref{it:colim1} and~(4.ii). of Definition~\ref{def:zigzagcolimit} exist.
\end{corollary}
\begin{proof} The colimit in step~\ref{it:colim1} exists since the singular monotone functor $S_{\cC}:Z_{\cC}\to \Simp$ preserves connected, non-empty colimits (Proposition~\ref{prop:singpreserves}.)
The existence of the colimit in step~(4.ii) is equivalent to the representability of the singular objects in step 3 of Corollary~\ref{cor:cocomplete}. This follows since $Z_{\cC} \to Z_{\widehat{\cC}}$ preserves connected, non-empty colimits (Corollary~\ref{cor:ICpreserve}.)
\end{proof}

\noindent
We can now combine Proposition~\ref{prop:thmpart1} and Corollary~\ref{cor:thmpart2} into a proof of Theorem~\ref{thm:zigzagcolimits}.

\begin{proof}[Proof of Theorem~\ref{thm:zigzagcolimits}]
Proposition~\ref{prop:thmpart1} asserts that if the procedure succeeds, then the cocone constructed in Definition~\ref{def:zigzagcolimit} is colimiting. Conversely, Corollary~\ref{cor:thmpart2} shows that if $D:\cJ \to Z_{\cC}$ has a colimit, then the procedure succeeds. 
\end{proof}

%% file: section4-homotopies.tex
\section{Homotopy construction}
\label{sec:homotopies}

\def\gap{\hspace{0.1cm}}

\noindent
Here we show that contraction can be used as a general technique to construct  nontrivial homotopies. In particular, we analyze the 4\-dimensional ``naturality'' homotopy, and the 5\-dimensional ``naturality of naturality'' homotopy. We first introduce some simple additional techniques, which are used together with contraction in the tool to produce these examples.

The examples come with direct links to the formalized proofs in the online proof assistant \homotopyio, where the interested reader can investigate them. To explore them, change the parameters of the ``Slice'' control at the top-right. You can also manipulate them directly; for example to execute a homotopy, use the mouse to drag a vertex (or a crossing) up or down, or drag a wire to the left or right. Further guidance on using the tool is available on the $n$Lab~\cite{nlab:homotopyio}.

As well as contraction, the tool makes use of some simple additional homotopy construction methods, which we summarize in Section~\ref{sec:additionalmethods}.


\subsection{Naturality (\href{http://www.cs.bham.ac.uk/~vicaryjo/homotopy.io/lics2019/naturality.html}{Link to online proof})}

\noindent
Here we build the following ``naturality'' homotopy, during which a vertex moves through a braiding, as the following zigzag of length 1 in $Z_\cat L^3$ (or equivalently, as an object of $Z_\cat L^4$):
\begin{equation}
\label{eq:shortnaturalitymovie}
\alignedgraphics{width=2cm}{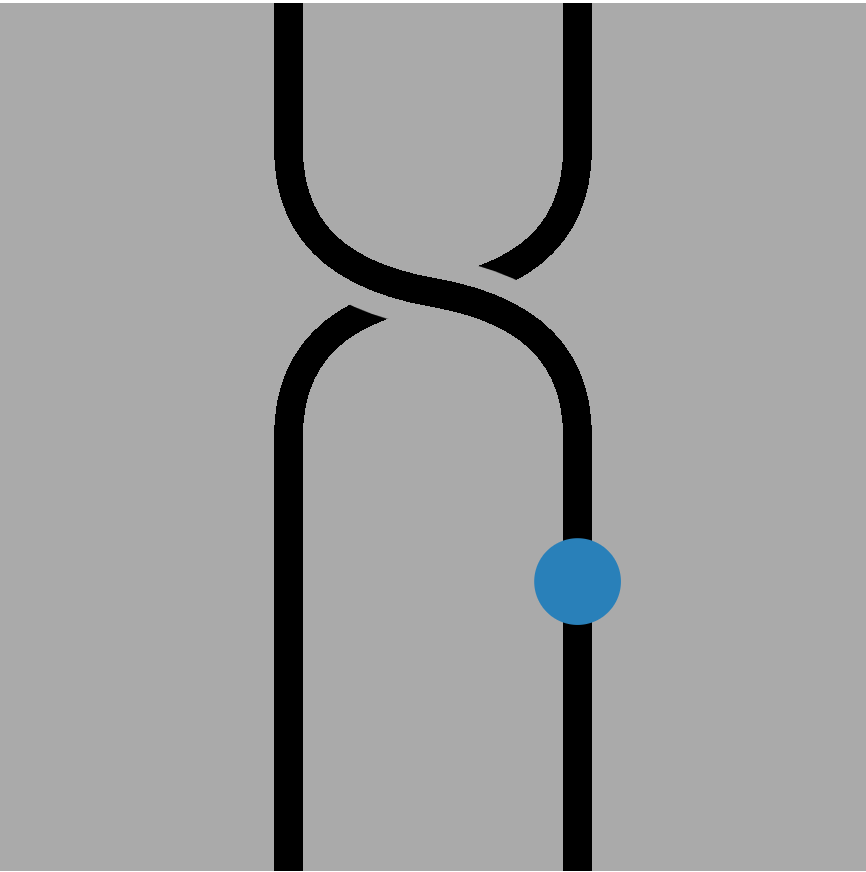}
\gap\to\gap
\alignedgraphics{width=2cm}{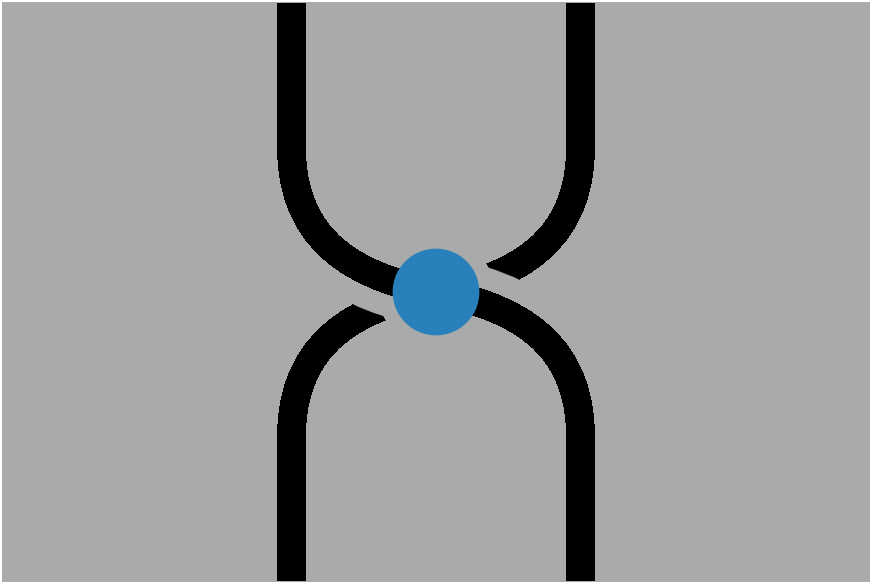}
\gap\ot\gap
\alignedgraphics{width=2cm}{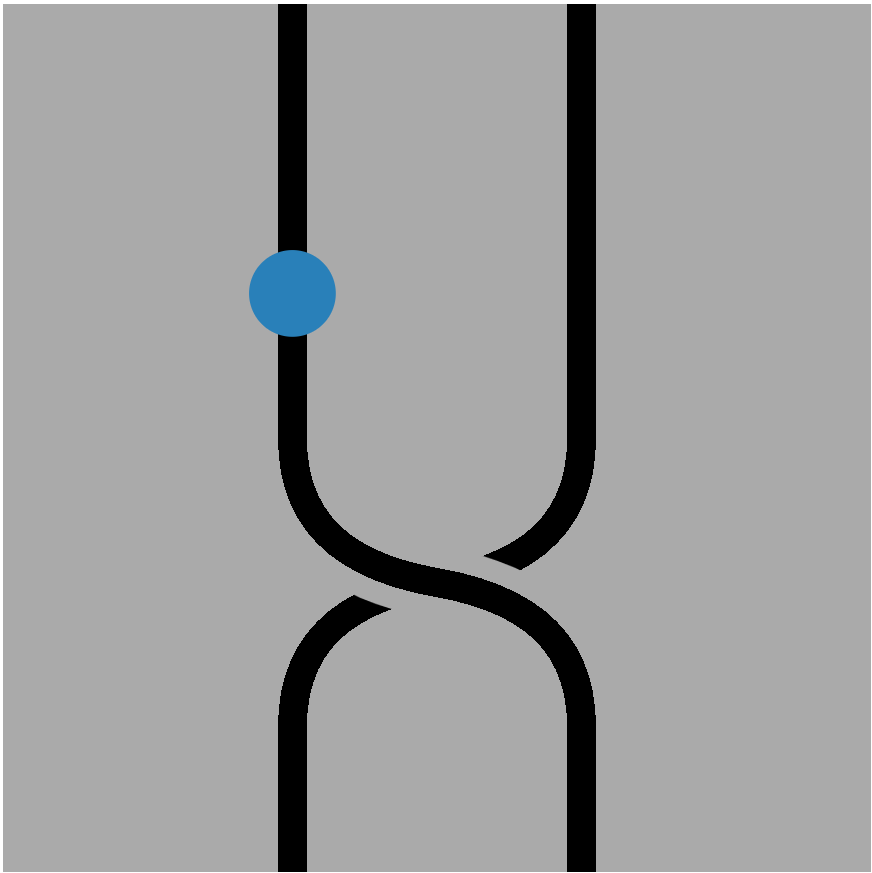}
\end{equation}

\begin{figure*}
\def\figwidth{1.5cm}
\def\figgap{\hspace{0.4cm}}
\[
\alignedgraphics{width=\figwidth}{images/naturality_1}
\figgap\to[*]\figgap
\alignedgraphics{width=\figwidth}{images/naturality_2}
\figgap\ot\figgap
\alignedgraphics{width=\figwidth}{images/naturality_2}
\figgap\to[\dag]\figgap
\alignedgraphics{width=\figwidth}{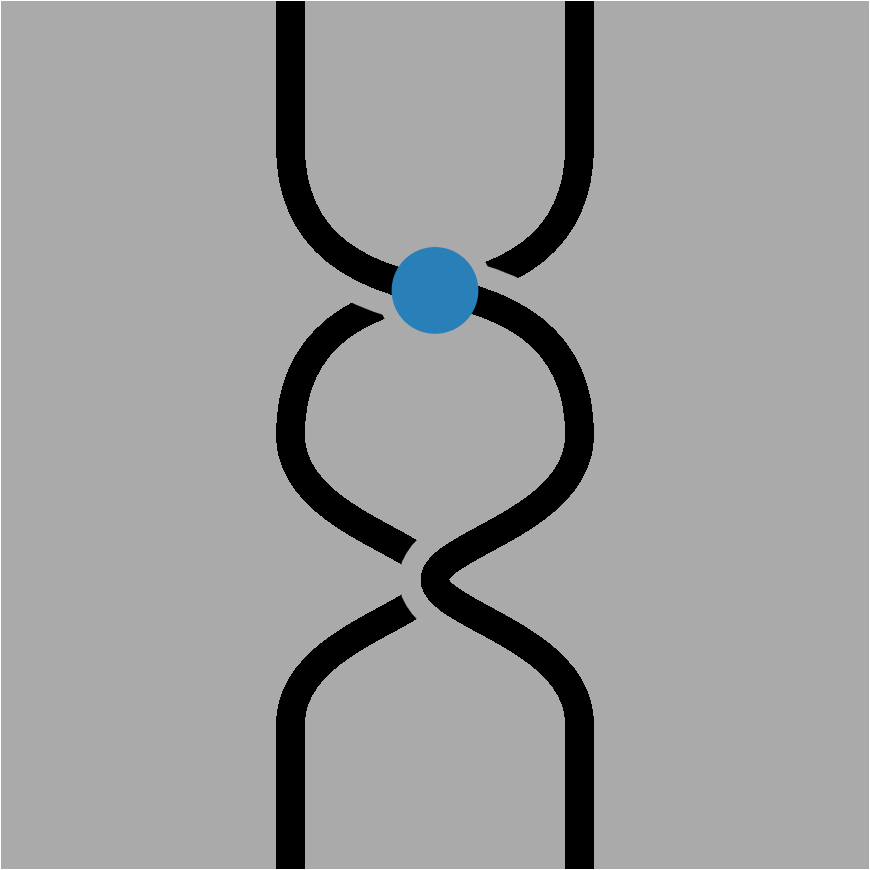}
\figgap\ot\figgap
\alignedgraphics{width=\figwidth}{images/naturality_3}
\figgap\to\figgap
\alignedgraphics{width=\figwidth}{images/naturality_3}\]
\[
\cdots 
\figgap\ot\figgap
\alignedgraphics{width=\figwidth}{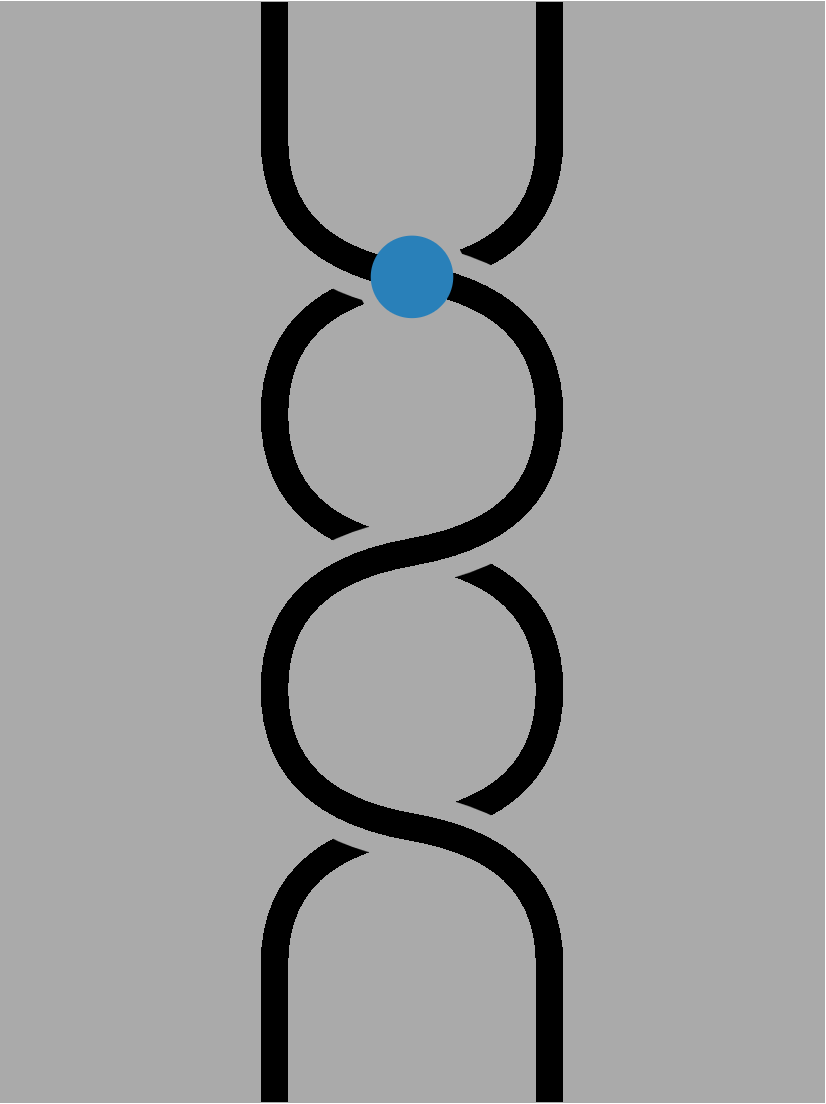}
\figgap\to\figgap
\alignedgraphics{width=\figwidth}{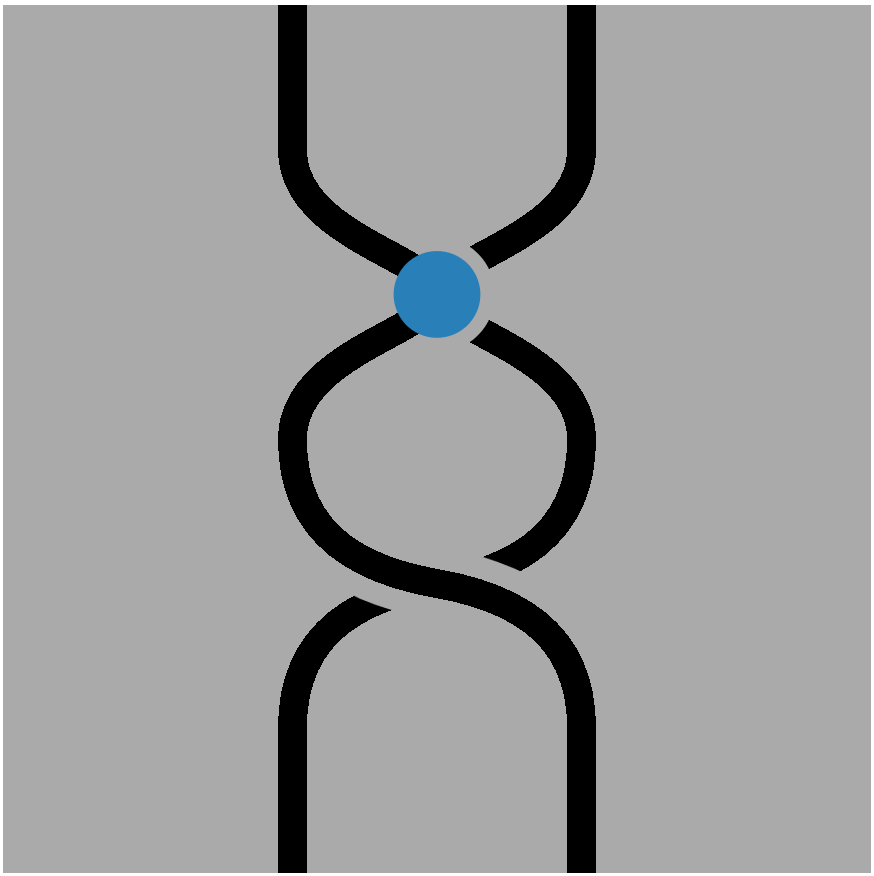}
\figgap\ot\figgap
\alignedgraphics{width=\figwidth}{images/naturality_5}
\figgap\to\figgap
\alignedgraphics{width=\figwidth}{images/naturality_5}
\figgap\ot\figgap
\alignedgraphics{width=\figwidth}{images/naturality_6}
\]
\caption{\label{fig:constructnaturality}Constructing the 4-dimensional naturality homotopy.}
\end{figure*}

To construct this homotopy, we begin by following the steps illustrated in Figure~\ref{fig:constructnaturality}, yielding a proof which is a length-5 zigzag in $Z^3_\cat L$. Each of these steps is obtained by contracting, or performing one of the recursive methods of Section~\ref{sec:additionalmethods} in a slice of the diagram; for example, in the arrow labelled $*$ we contract the entire diagram, and in the arrow labelled $\dag$ we perform a contraction within the first regular height of the diagram. By projecting out an extra dimension, we can view the entire proof that we have constructed as a 2\-dimensional graphic, giving information about the overall structure of our proof, as shown in the first image here:
\begin{equation}
\label{eq:naturalityproof}
\alignedgraphics{width=3.5cm}{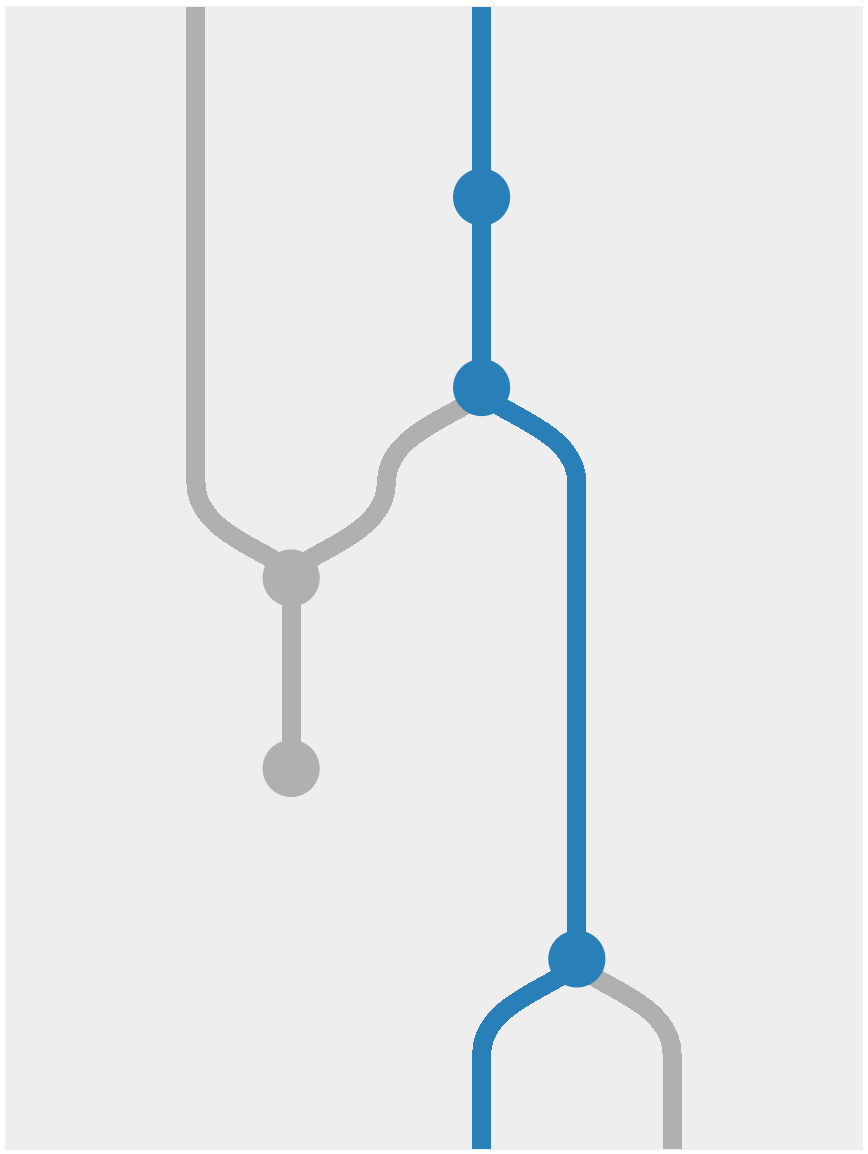}
\quad\leadsto\quad
\alignedgraphics{width=2.3cm}{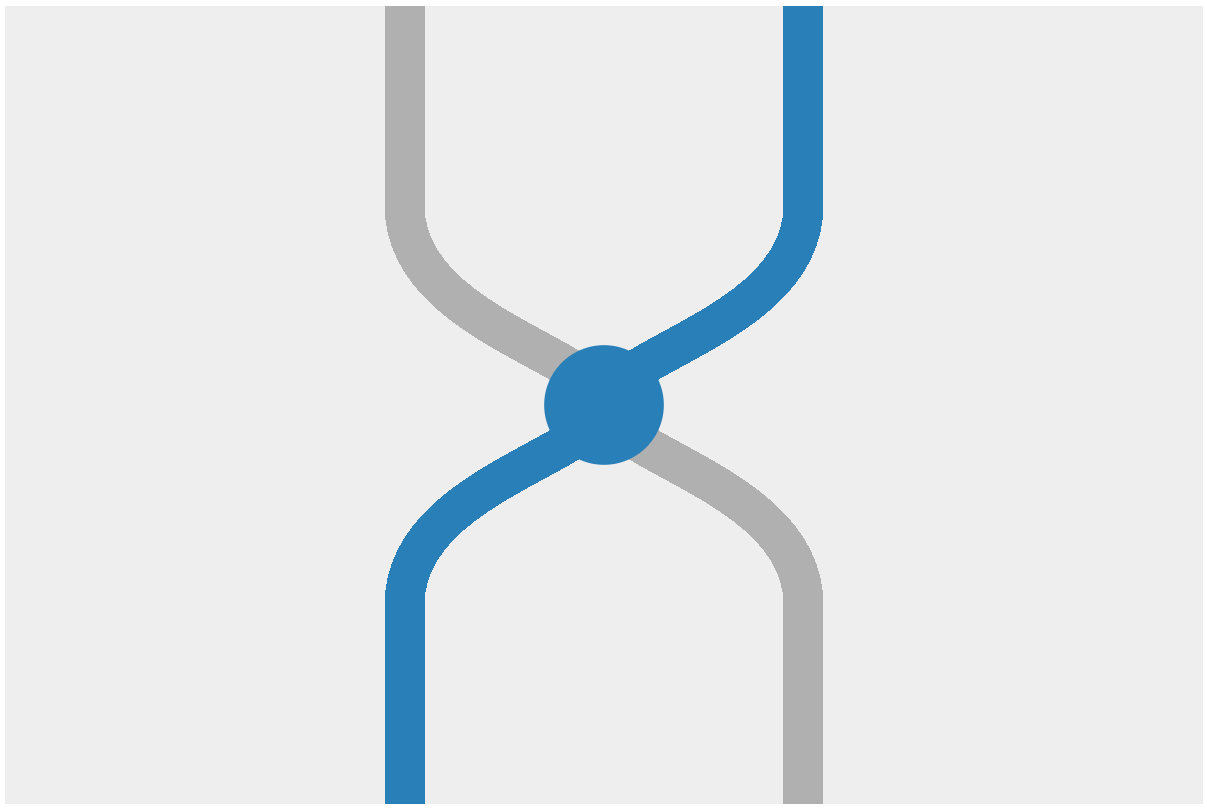}
\end{equation}
We then contract this, and this entire proof collapses to a zigzag of length 1, which performs the naturality move in a single step, shown in projection as the second diagram above. Viewing this as a ``movie'' gives back precisely the desired homotopy~\eqref{eq:shortnaturalitymovie} above.

\subsection{Naturality of naturality (\href{http://www.cs.bham.ac.uk/~vicaryjo/homotopy.io/lics2019/naturality_of_naturality.html}{Link to online proof})}

\noindent
This homotopy has the following 4\-dimensional diagrams as its source and target respectively:
\def\gap{\vspace{0.1cm}}
\begin{calign}
\label{eq:nnsource}
\alignedgraphics{width=1.1cm}{images/naturality_1}
\gap\to\gap\alignedgraphics{width=1.1cm}{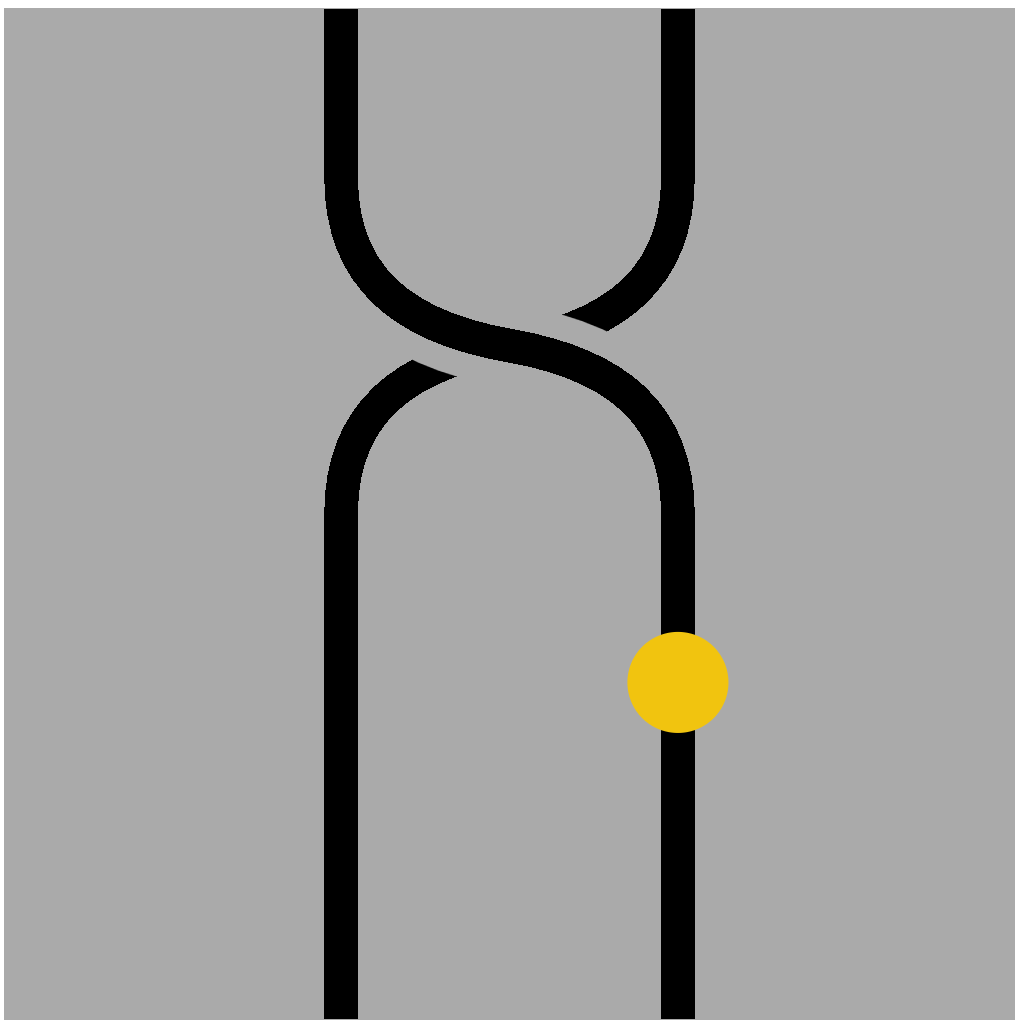}
\gap\ot\gap\alignedgraphics{width=1.1cm}{images/naturality_1}
\gap\to\gap\alignedgraphics{width=1.1cm}{images/naturality_2}
\gap\ot\gap\alignedgraphics{width=1.1cm}{images/naturality_6}
\\*
\label{eq:nntarget}
\alignedgraphics{width=1.1cm}{images/naturality_1}
\gap\ot\gap\alignedgraphics{width=1.1cm}{images/naturality_2}
\gap\to\gap\alignedgraphics{width=1.1cm}{images/naturality_6}
\gap\to\gap\alignedgraphics{width=1.1cm}{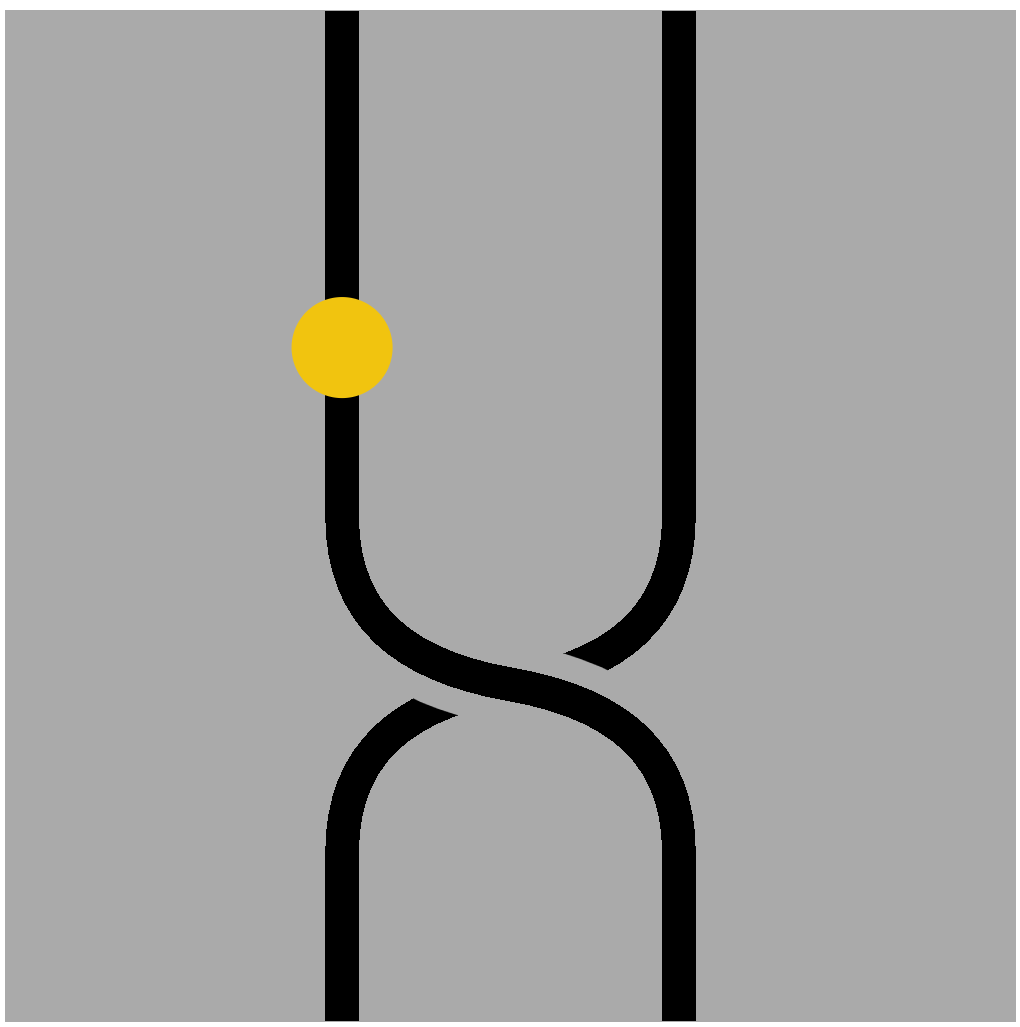}
\gap\ot\gap\alignedgraphics{width=1.1cm}{images/naturality_6}
\end{calign}
These diagrams feature a 3\-cell drawn in blue, a crossing, and a 4\-cell drawn in yellow, which acts as an endomorphism of the blue 3\-cell. The source~\eqref{eq:nnsource} describes a composite process that applies the 4\-cell to the 3\-cell, then pulls the 3\-cell through the braiding; in the target~\eqref{eq:nntarget}, we instead first pull the 3\-cell through the braiding, and then apply the 4\-cell. These source and target 4\-diagrams are not homotopies, since they involve the yellow 4\-cell, which is an algebraic move.

\begin{figure*}
\def\figwidth{1.3cm}
\def\figgap{\hspace{0.0cm}}
$
\alignedgraphics{width=\figwidth}{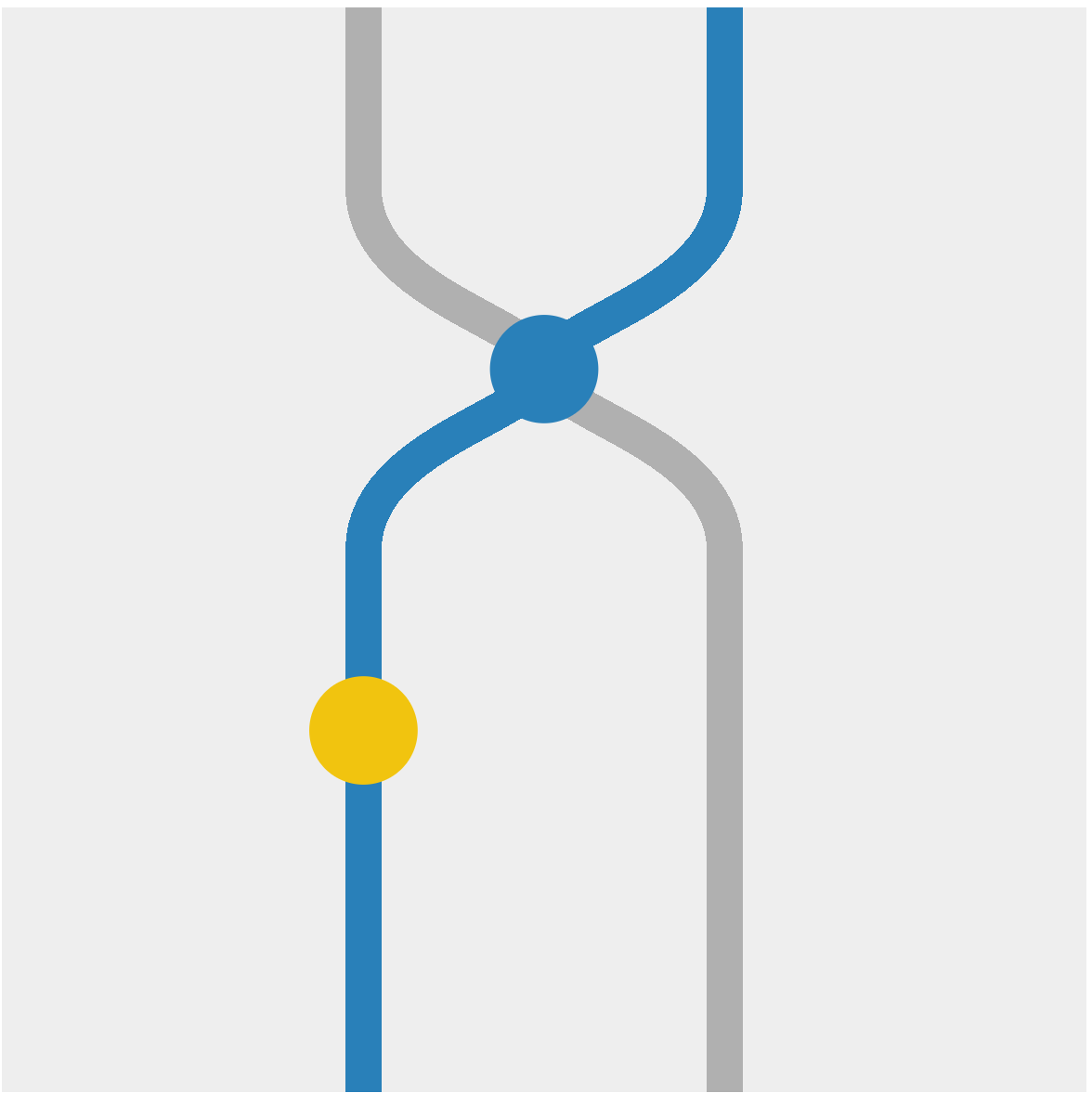}
\figgap\to\figgap\alignedgraphics{width=\figwidth}{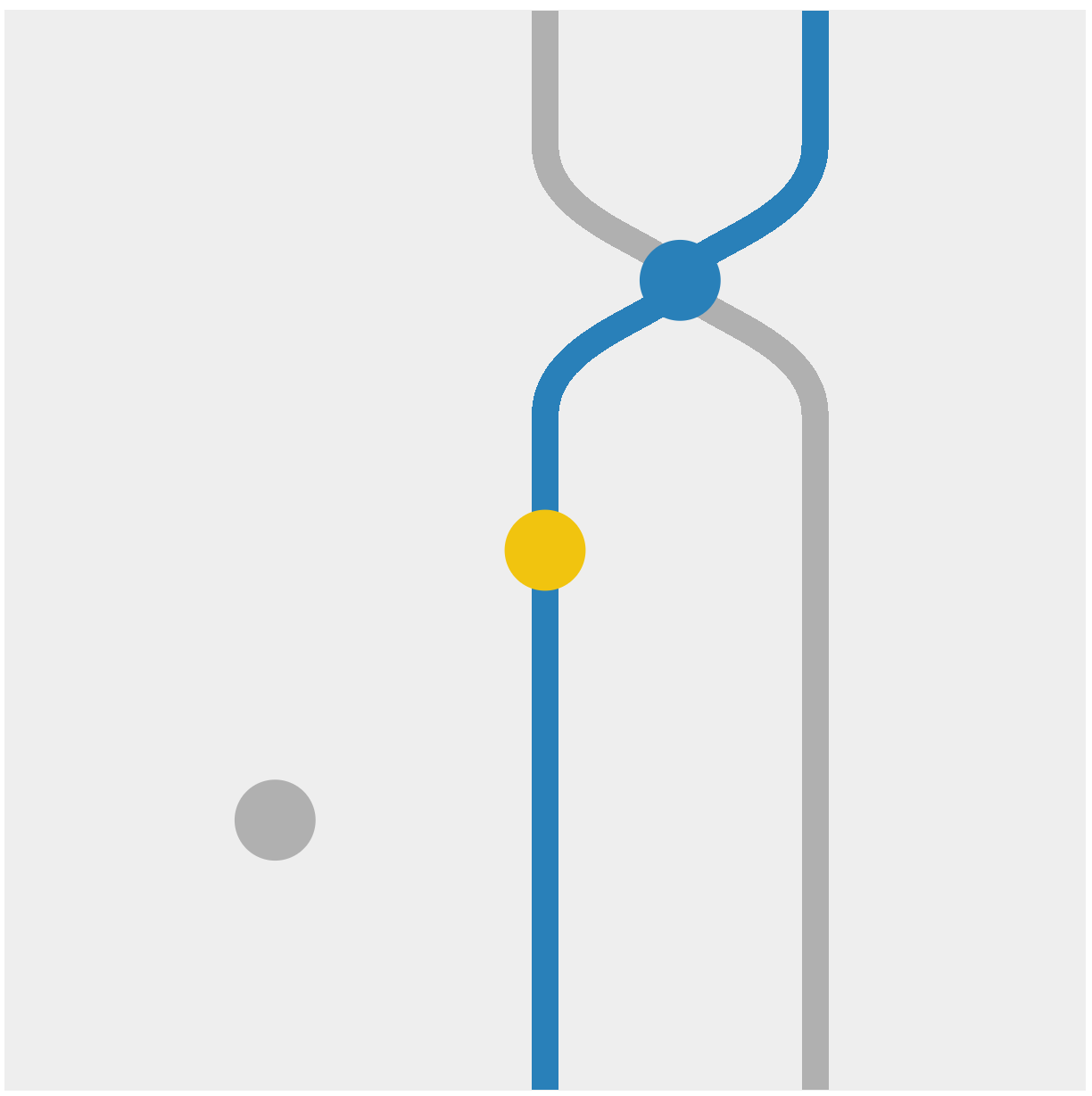}
\figgap\ot\figgap\alignedgraphics{width=\figwidth}{images/nn_2}
\figgap\to\figgap\alignedgraphics{width=\figwidth}{images/nn_2}
\figgap\ot\figgap\alignedgraphics{width=\figwidth}{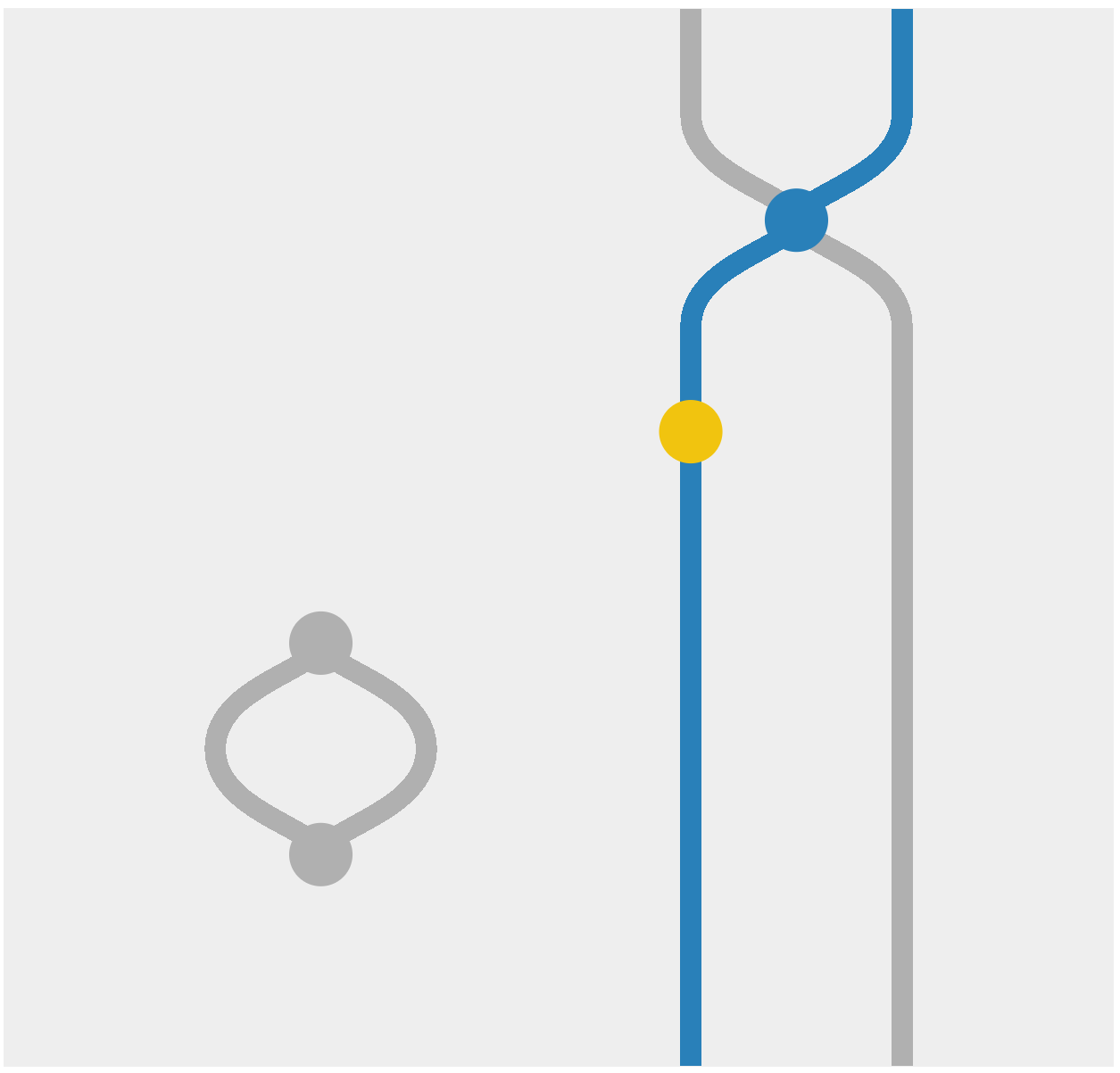}
\figgap\to\figgap\alignedgraphics{width=\figwidth}{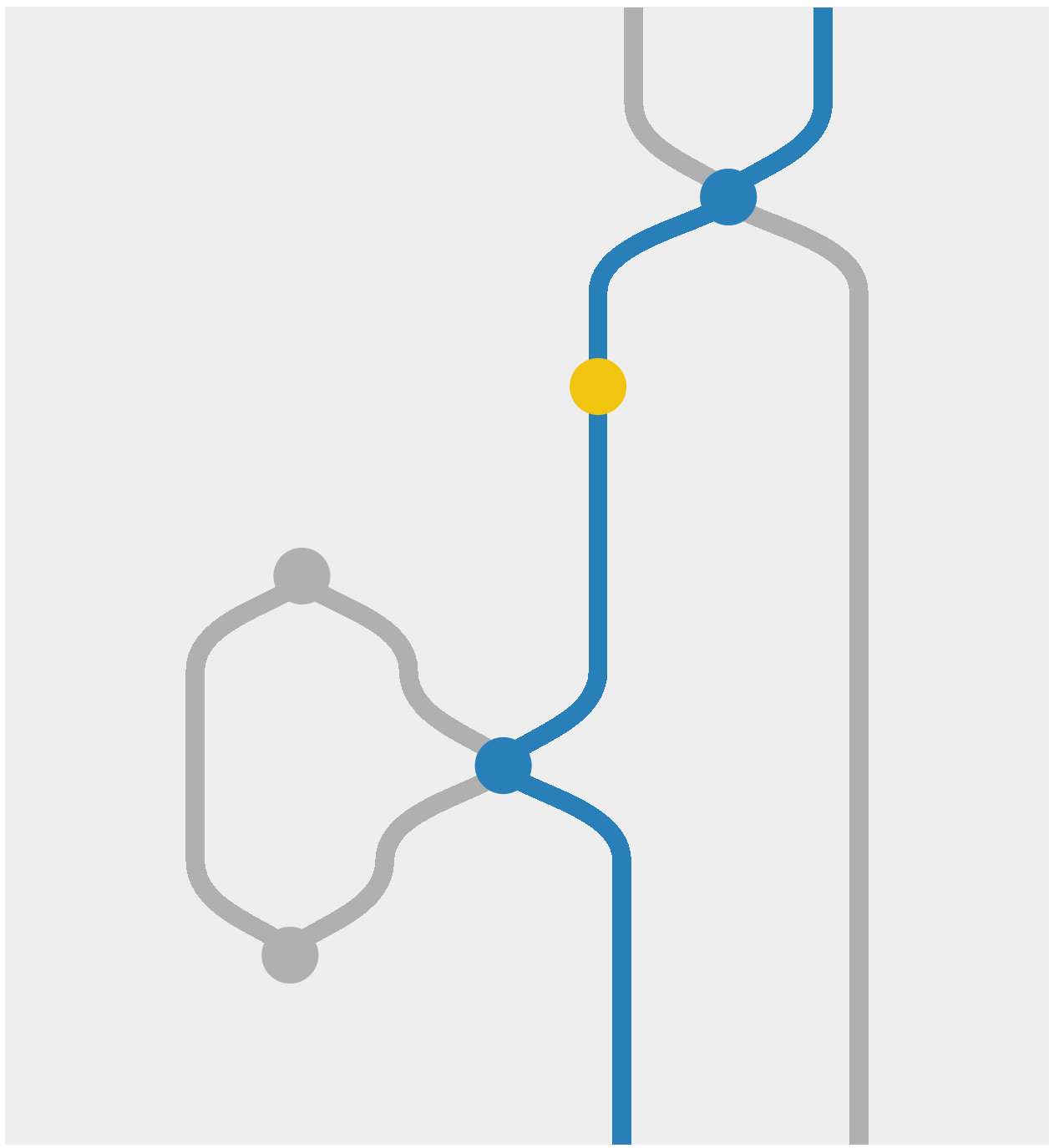}
\figgap\ot\figgap\alignedgraphics{width=\figwidth}{images/nn_4}
\figgap\to\figgap\alignedgraphics{width=\figwidth}{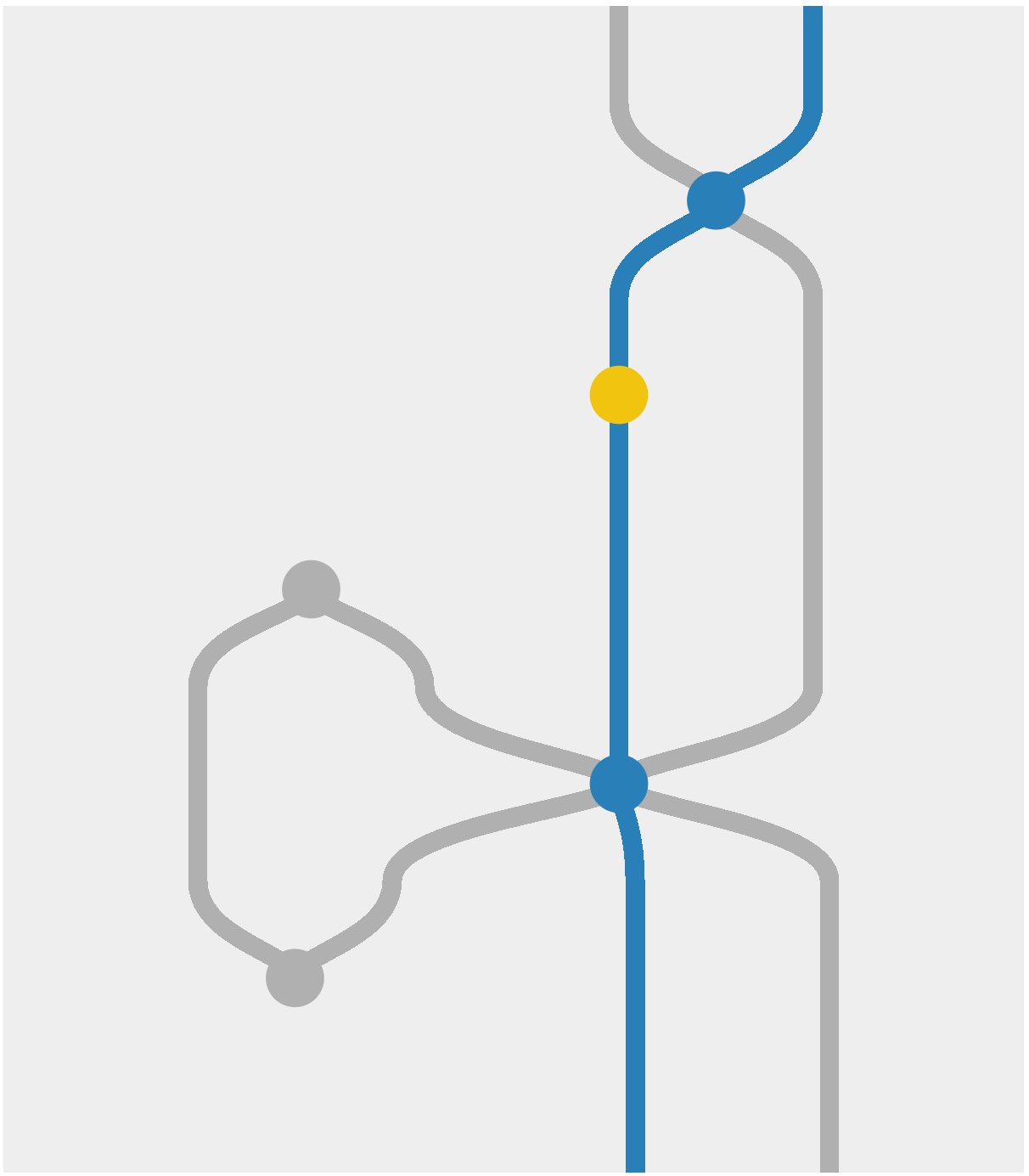}
\\\ldots\figgap\ot\figgap\alignedgraphics{width=\figwidth}{images/nn_5}
\figgap\to\figgap\alignedgraphics{width=\figwidth}{images/nn_5}
\figgap\ot\figgap\alignedgraphics{width=\figwidth}{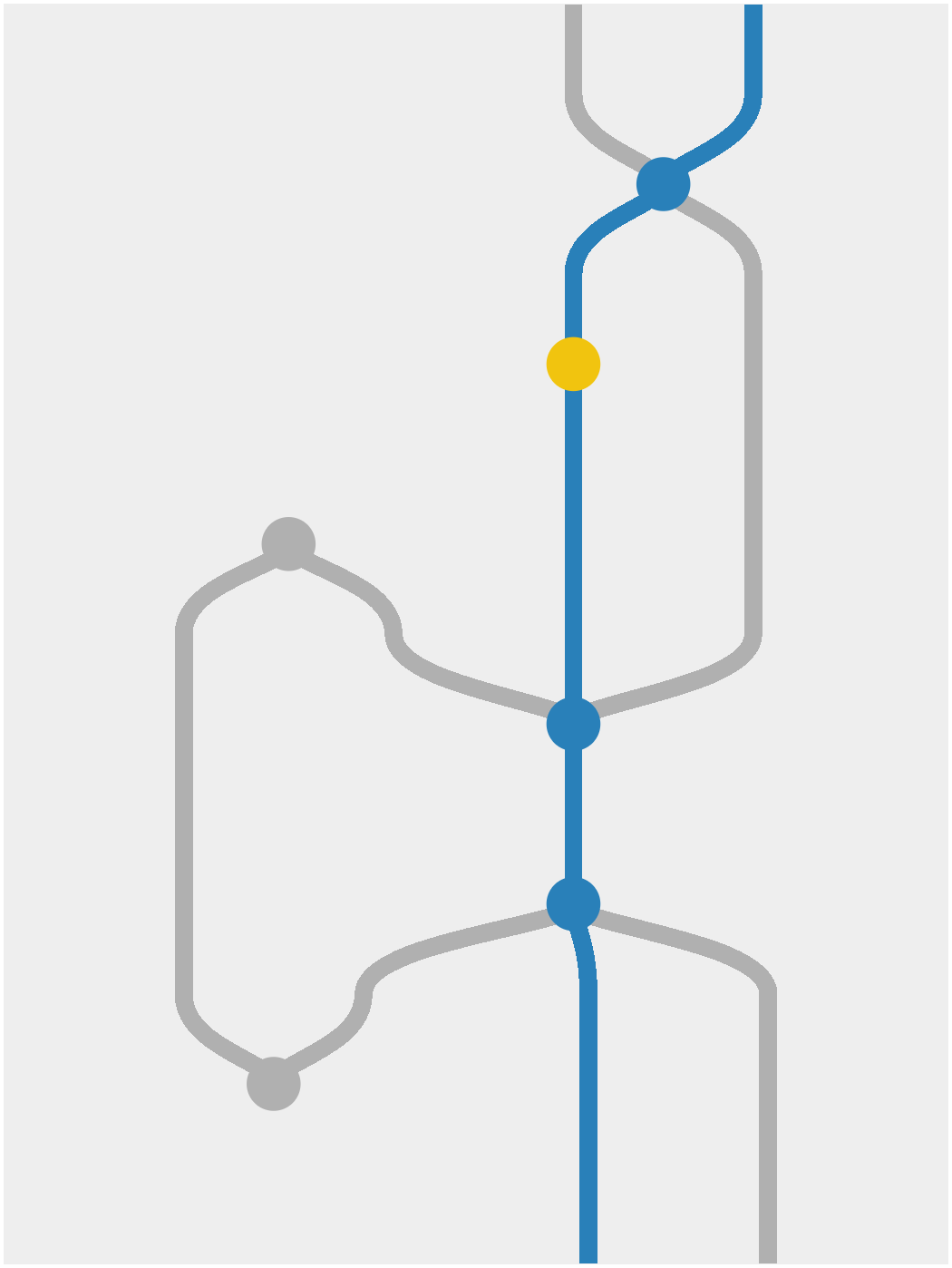}
\figgap\to\figgap\alignedgraphics{width=\figwidth}{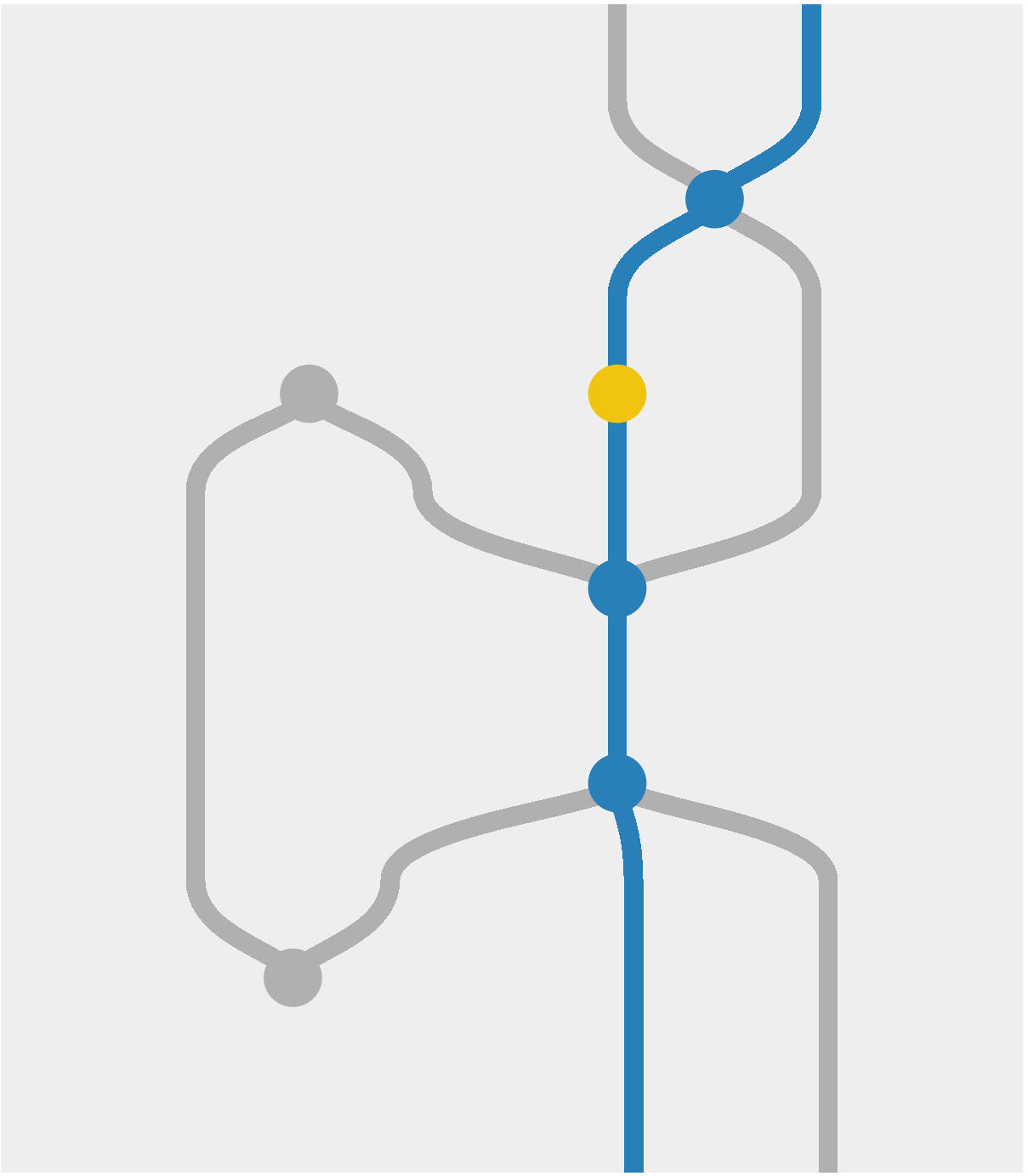}
\figgap\ot\figgap\alignedgraphics{width=\figwidth}{images/nn_7}
\figgap\to\figgap\alignedgraphics{width=\figwidth}{images/nn_7}
\figgap\ot\figgap\alignedgraphics{width=\figwidth}{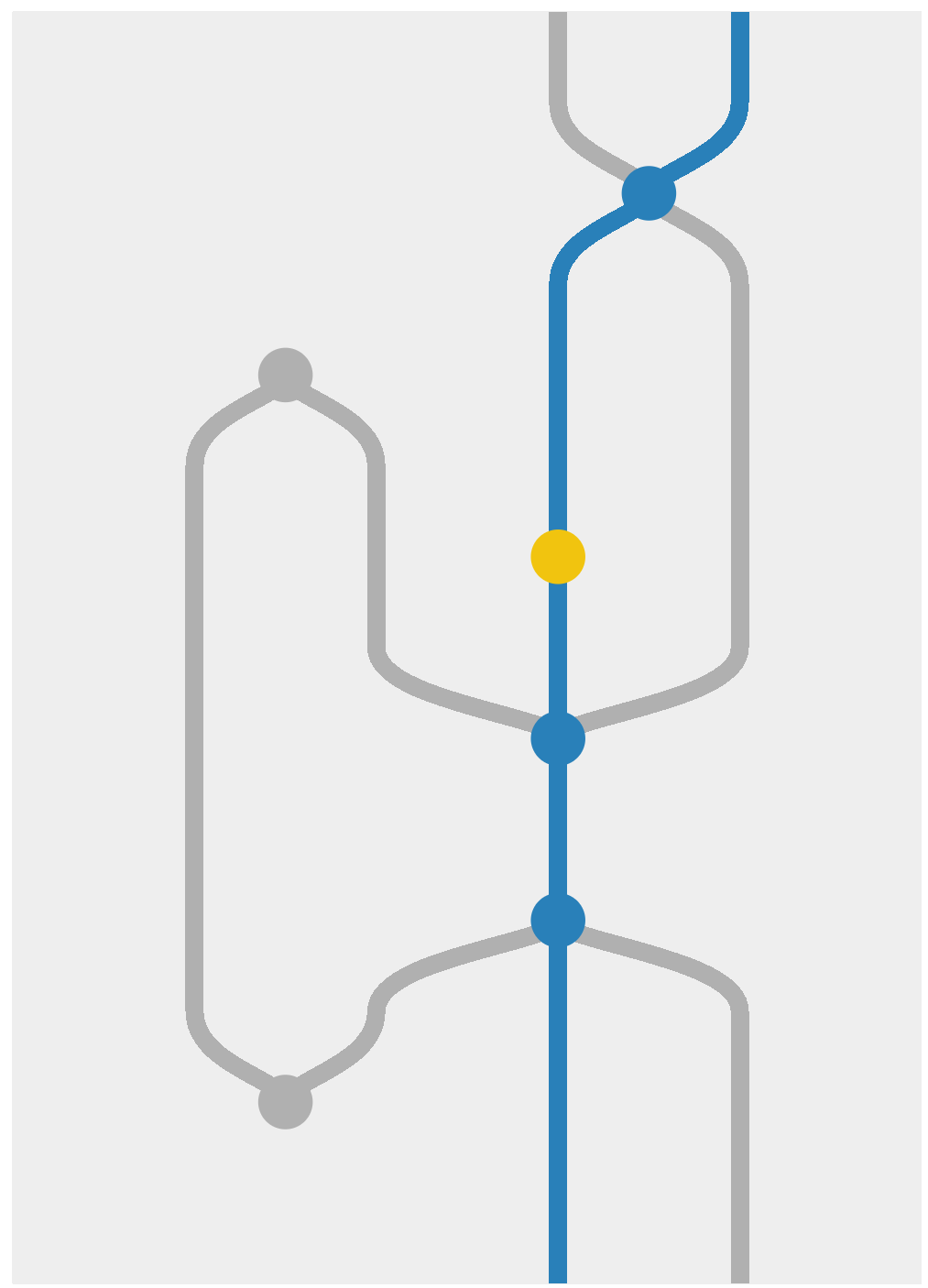}
\\\ldots\figgap\to\figgap\alignedgraphics{width=\figwidth}{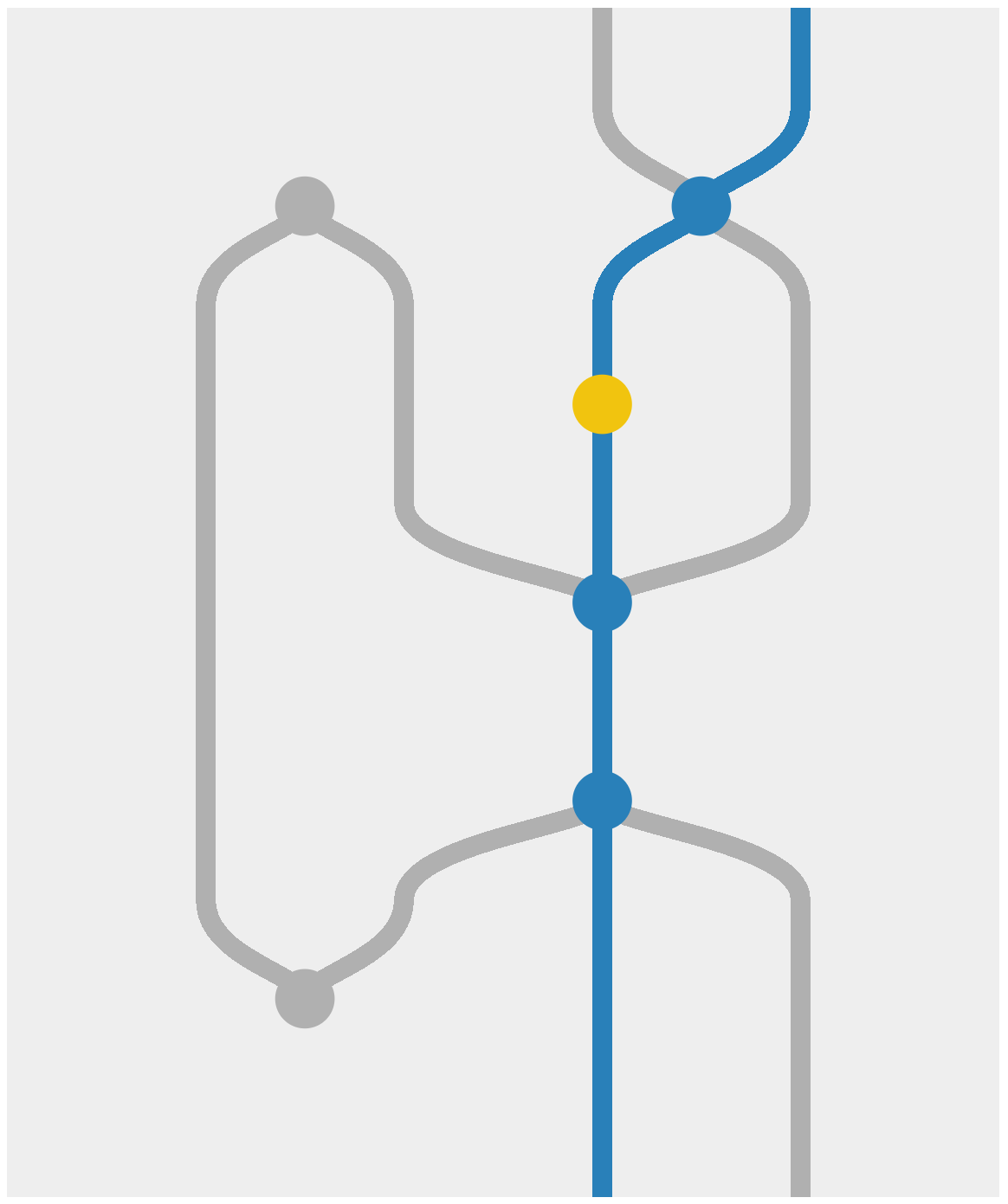}
\figgap\ot\figgap\alignedgraphics{width=\figwidth}{images/nn_9}
\figgap\to\figgap\alignedgraphics{width=\figwidth}{images/nn_9}
\figgap\ot\figgap\alignedgraphics{width=\figwidth}{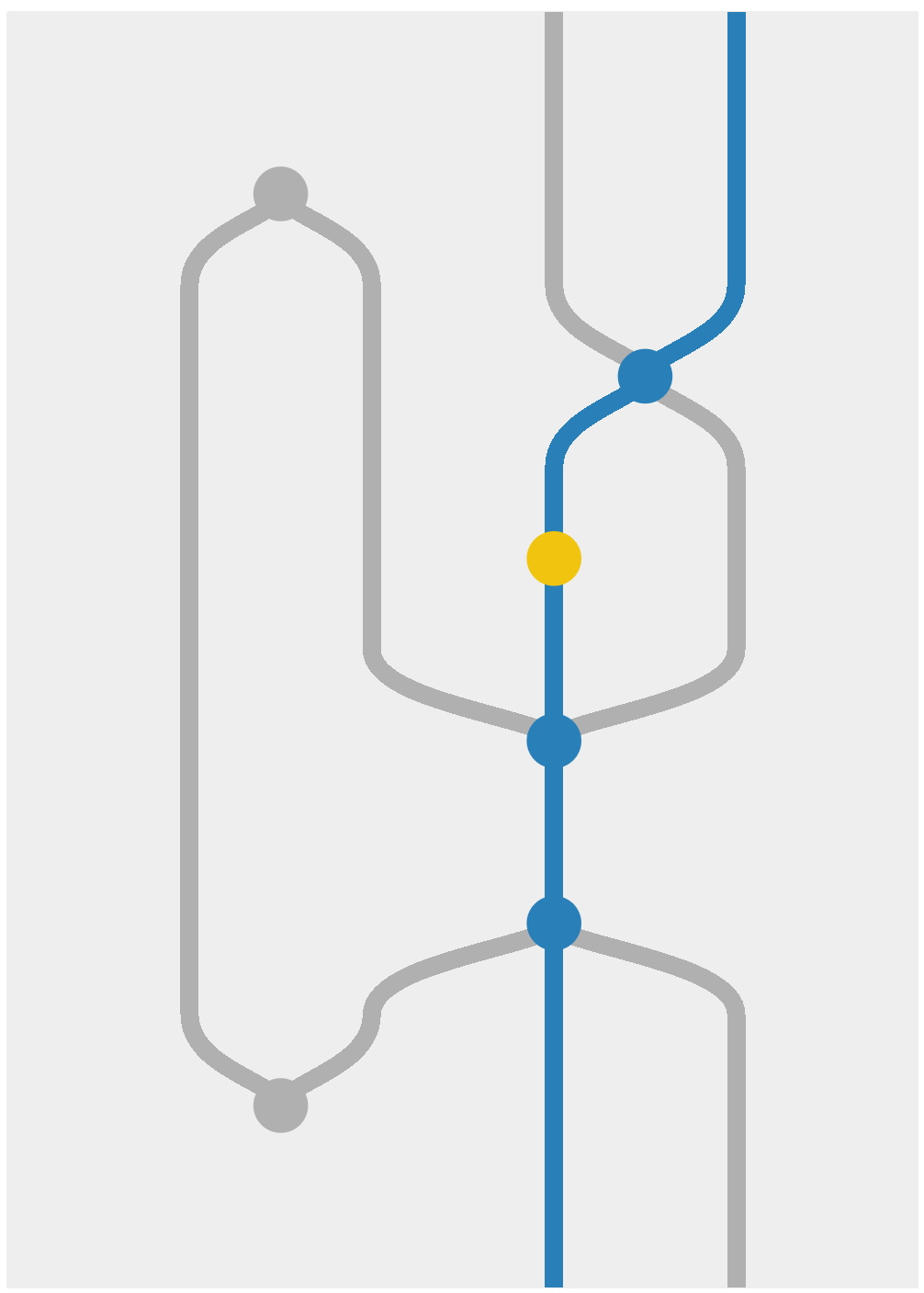}
\figgap\to\figgap\alignedgraphics{width=\figwidth}{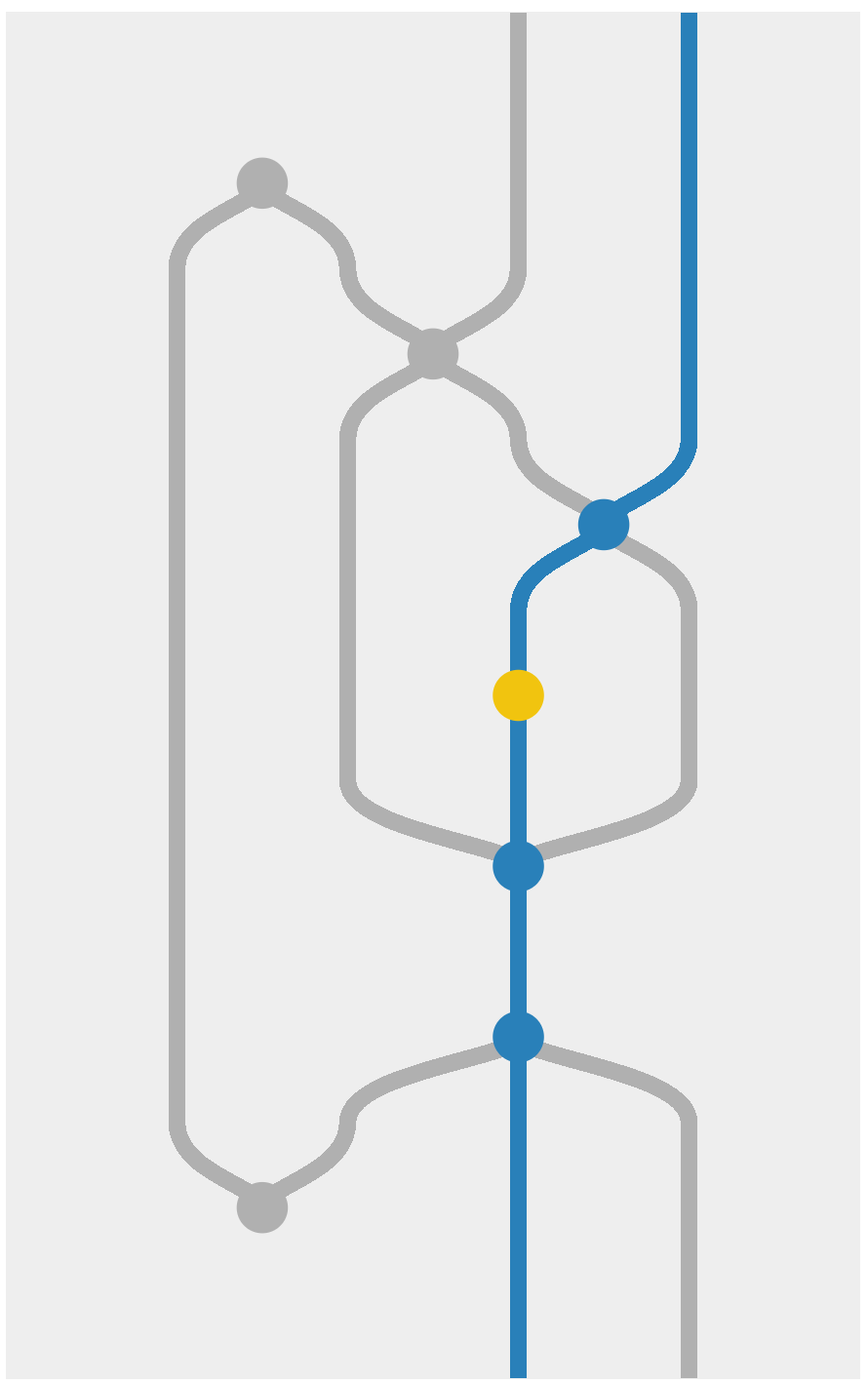}
\figgap\ot\figgap\alignedgraphics{width=\figwidth}{images/nn_11}
\figgap\to\figgap\alignedgraphics{width=\figwidth}{images/nn_11}
\\\ldots\figgap\ot\figgap\alignedgraphics{width=\figwidth}{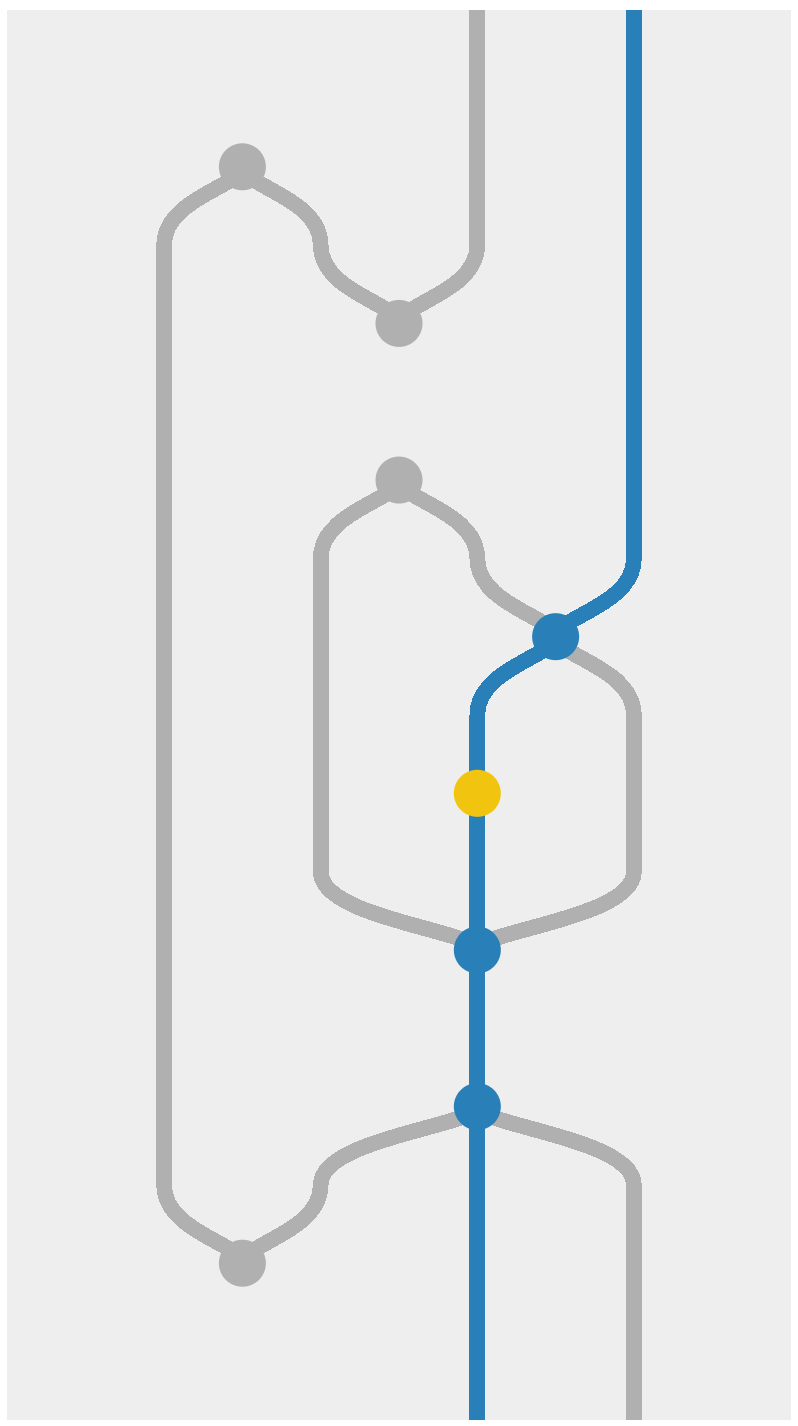}
\figgap\to\figgap\alignedgraphics{width=\figwidth}{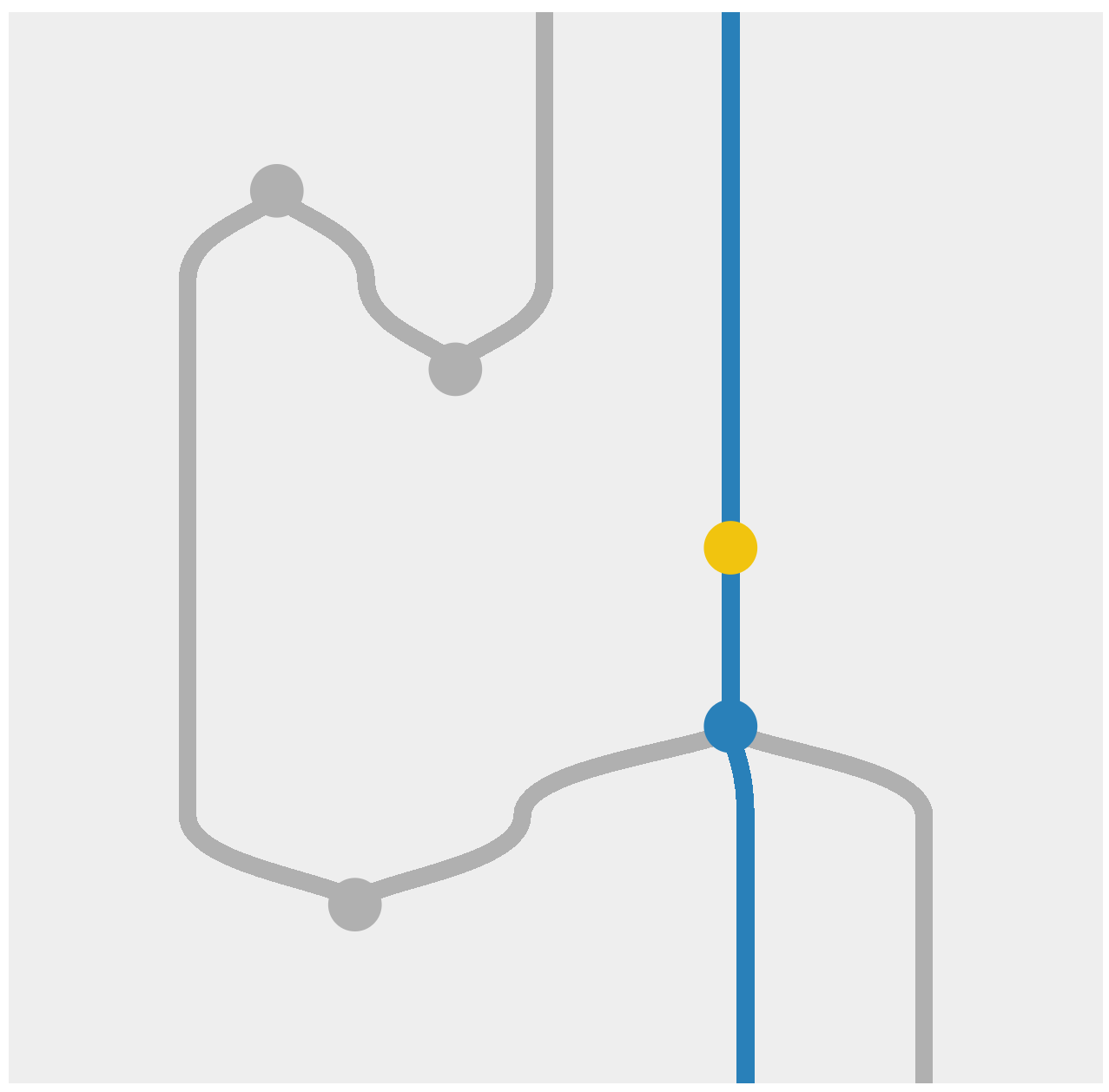}
\figgap\ot\figgap\alignedgraphics{width=\figwidth}{images/nn_13}
\figgap\to\figgap\alignedgraphics{width=\figwidth}{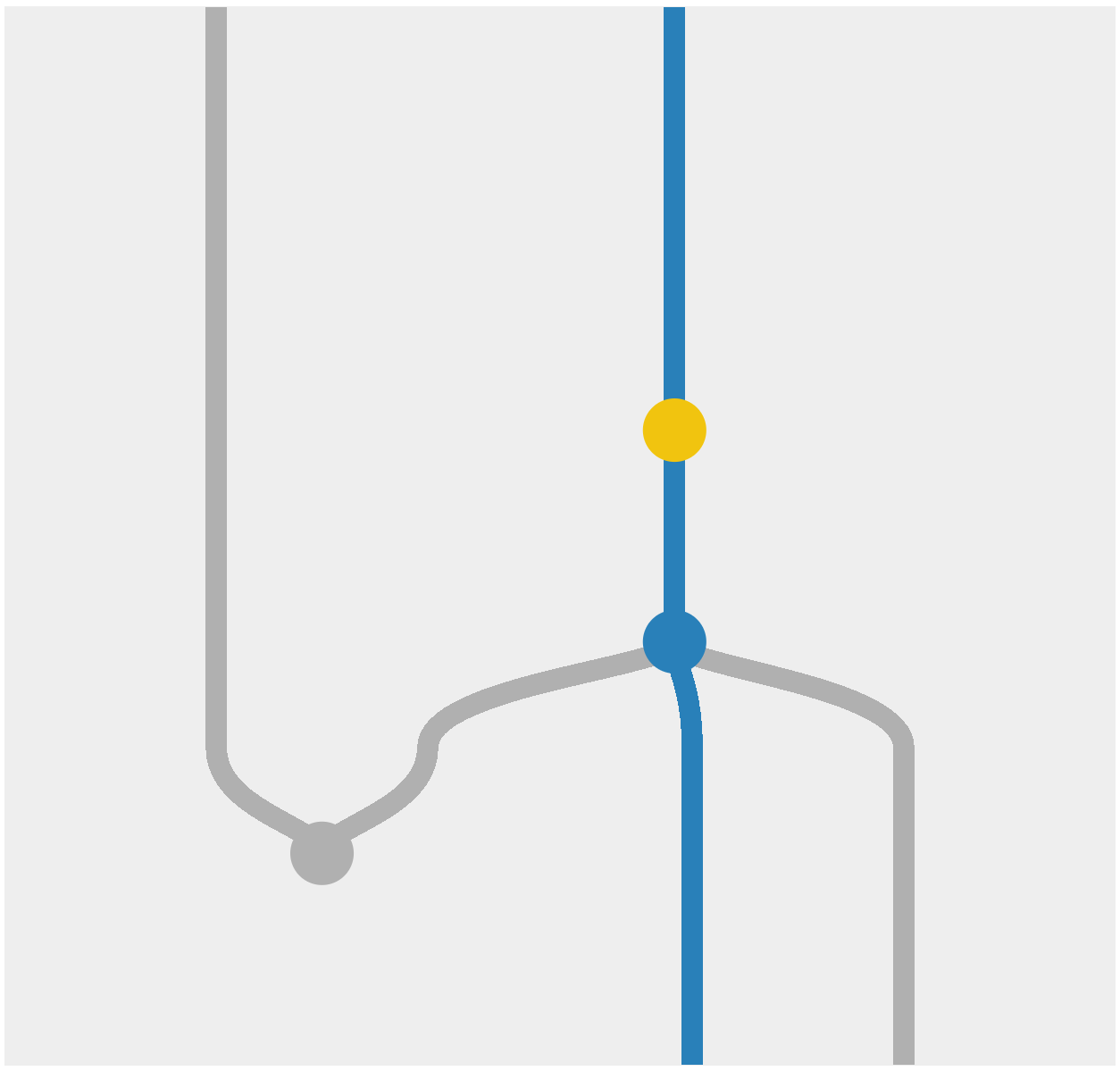}
\figgap\ot\figgap\alignedgraphics{width=\figwidth}{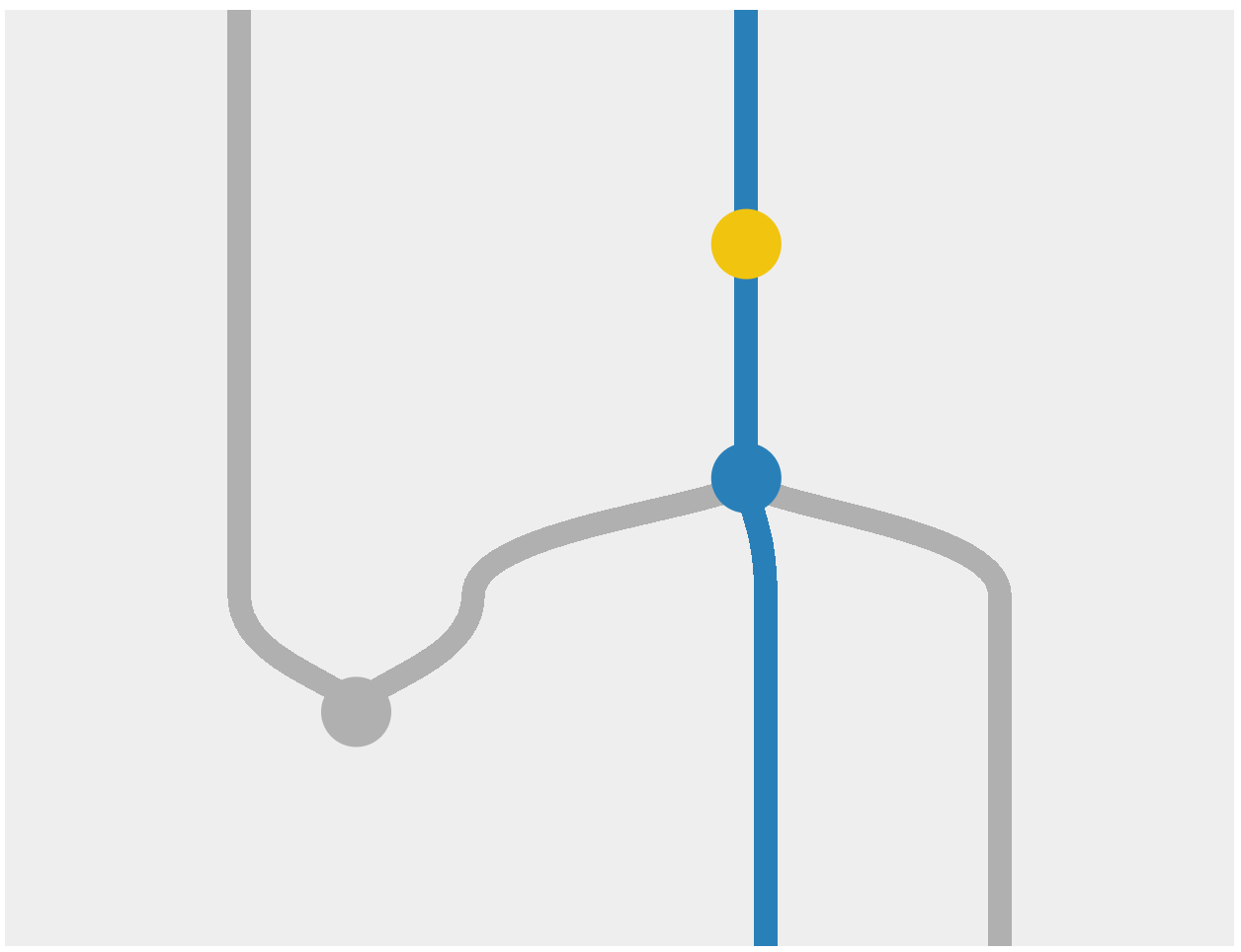}
\figgap\to\figgap\alignedgraphics{width=\figwidth}{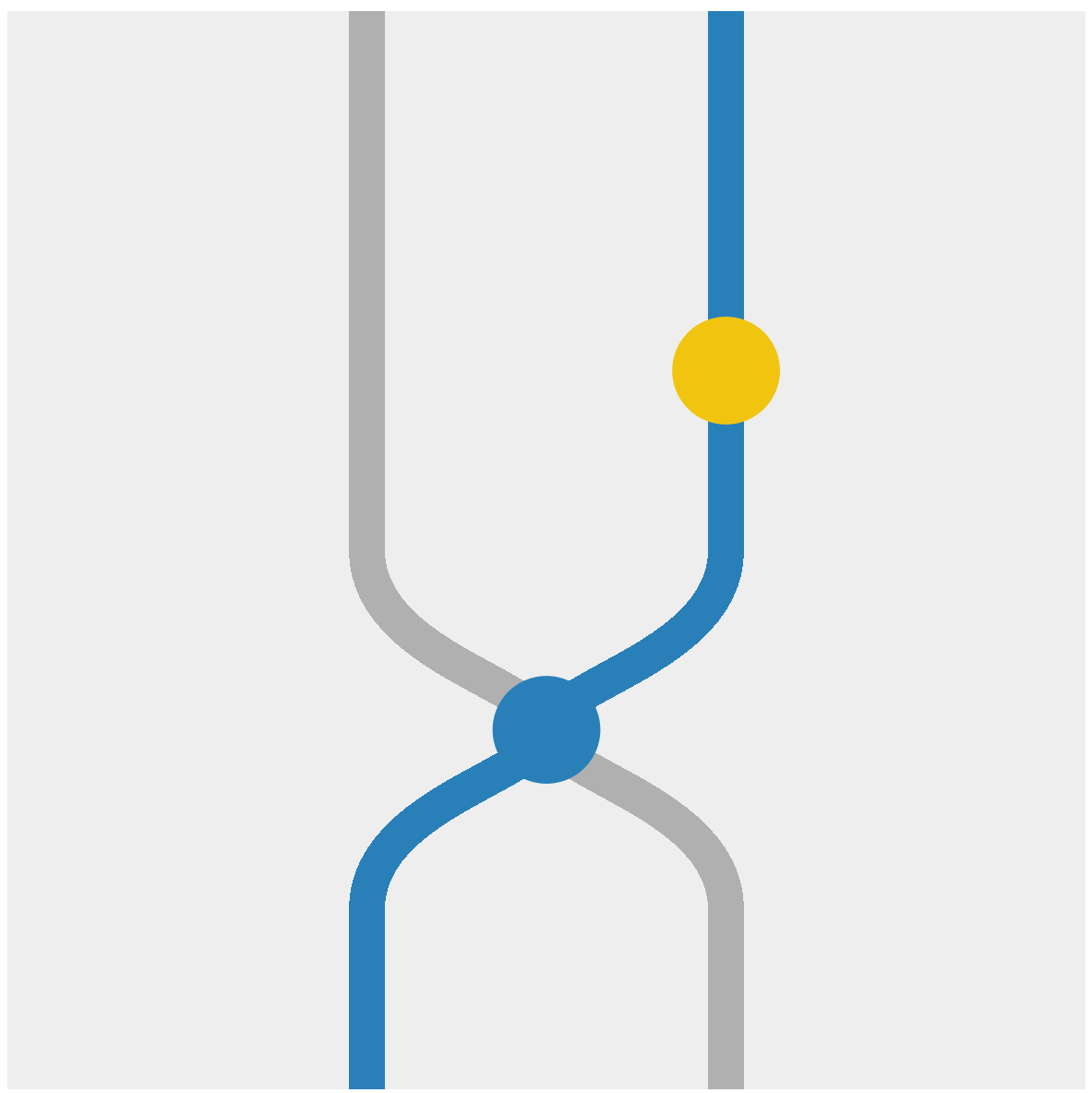}
\figgap\ot\figgap\alignedgraphics{width=\figwidth}{images/nn_16}
$
\caption{\label{fig:naturalitynaturality}Constructing the 5-dimensional naturality of naturality homotopy.}
\end{figure*}

The ``naturality of naturality'' homotopy exhibits that the composites~\eqref{eq:nnsource} and~\eqref{eq:nntarget}  are homotopic. The proof is constructed in \homotopyio as a zigzag of length 14 in $Z_\cat L^4$ (or alternatively, an object in $Z_\cat L^5$.) We build it by starting with the source 4\-diagram~\eqref{eq:nnsource}, and manipulating it using our contraction-based methods (that is, the contraction operation of Section~\ref{sec:contraction}, along with the additional recursive methods of Section~\ref{sec:additionalmethods}.) We give it as a movie (in which every frame shows a 4\-dimensional diagram, viewed in 2\-dimensional projection) in Figure~\ref{fig:naturalitynaturality}.

Note that the first frame of this movie is given by the projection of the contracted naturality homotopy (the second image in~\eqref{eq:naturalityproof}), composed with the yellow 4\-cell, and the last frame has these same components composed a different way; the yellow 4\-cell is indeed ``pulled through the naturality'' over the course of the proof, as we expect.

The proof as a whole has an interesting structure, which the movie of Figure~\ref{fig:naturalitynaturality} makes clear: we create a bubble, enlarge it, wrap it around the yellow 4\-cell, and then contract the remaining parts, with the result being that the yellow 4\-cell has moved to the other side of the naturality homotopy. Building this proof required repeated use of contraction, not only on the 4\-dimensional term being manipulated as the proof was being developed, but also in 3\-dimensional slices of those terms, which then propagated recursively to the entire 4\-diagram by the methods of Section~\ref{sec:additionalmethods}.

As before, the entire 5\-dimensional proof can be viewed as a single 2\-dimensional projected image, by ignoring the lowest 3 dimensions. We represent it in the first image here:
\begin{equation}
\alignedgraphics{height=5.5cm}{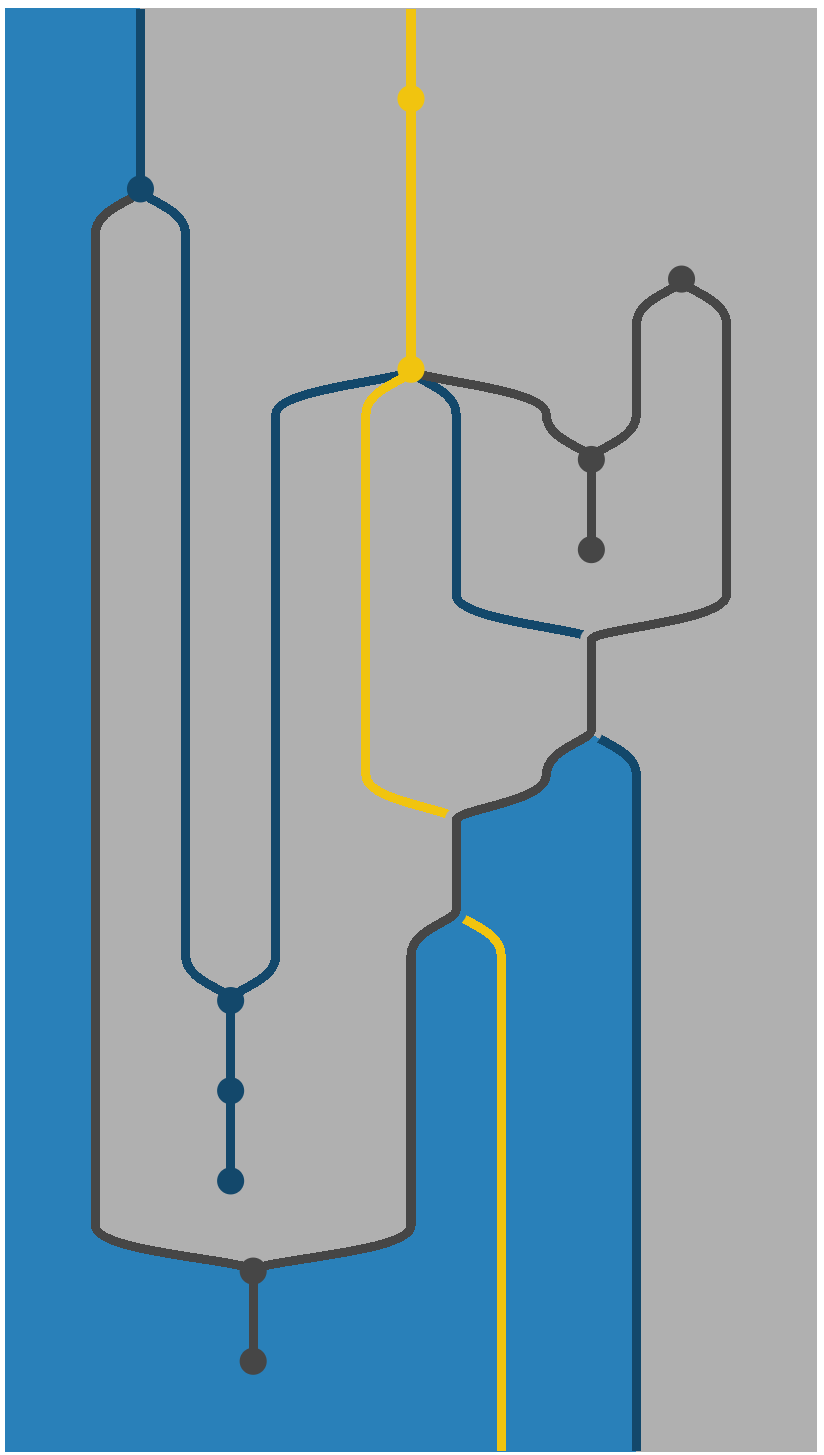}
\quad\leadsto\quad
\alignedgraphics{height=1.5cm}{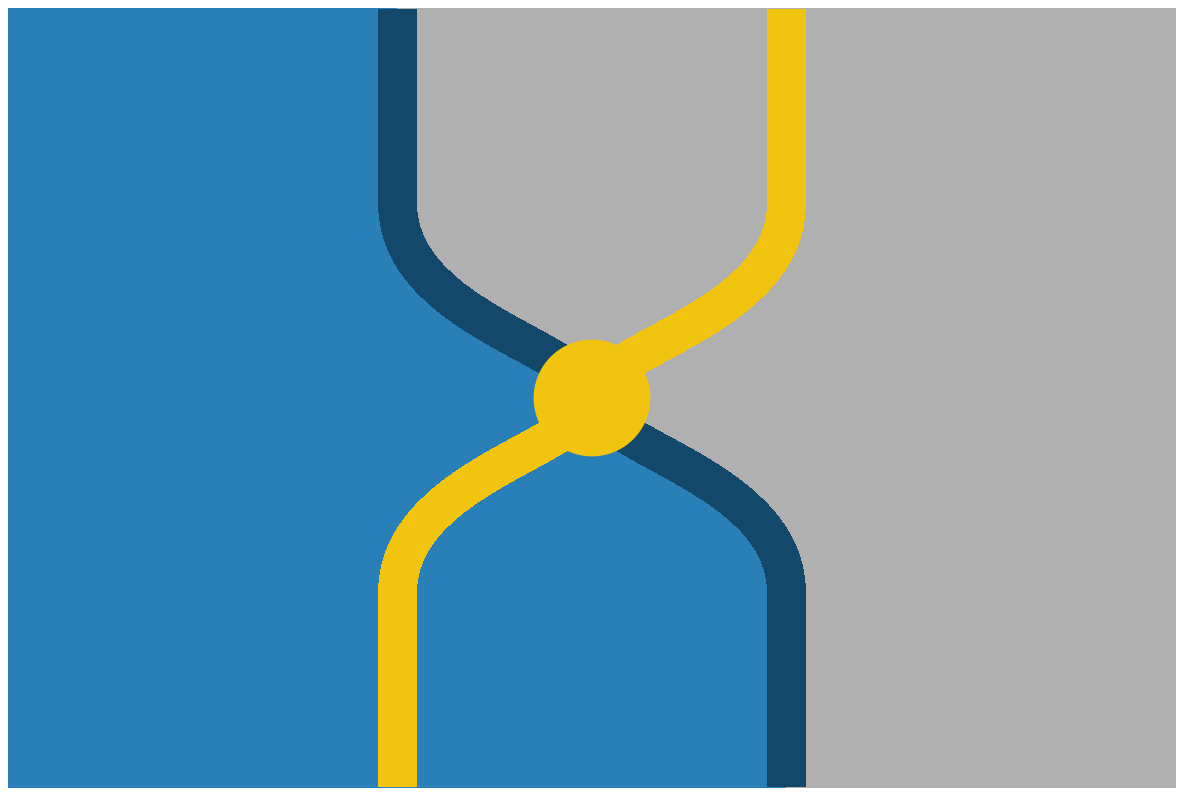}
\end{equation}
This entire proof contracts to a zigzag of length 1, and we give this contraction as the second image.

To summarize, we have shown  how contraction can be used as the main workhorse for manipulating terms in an associative $n$\-category, including the tasks of building an initial diagram, manipulating it (both at the top dimension and in lower dimensions) to obtain a proof object, and then contracting that proof object itself to yield a short witness for the logical statement being established.

%% file: sectionA-moremethods.tex

\section{Additional methods}

\label{sec:additionalmethods}

\noindent
As well as contraction, the implementation \homotopyio uses some additional methods to construct homotopies. We summarize these here. The main idea that as well as contractions, we have a dual notion of \emph{expansion}, and that we allow these to propagate recursively from regular or singular heights of a diagram, to the main diagram, wherever possible. These methods are far less complex than the basic contraction procedure described in Section~\ref{sec:contraction}, but extend its power considerably. 

\subsection{Recursive contraction}

\noindent
Here we define a \emph{generalized contraction} for a diagram $D$ as any morphism $D \to[c] C$ in the appropriate zigzag category. We show how these generalized contractions arise recursively, with the contractions arising from colimits, as described in Section~\ref{sec:contraction}, being the base case.

Here we give recursive schemes to describe the behaviour when a (generalized) contraction is triggered within a regular or singular height of the diagram. As with the base case contraction procedure in Section~\ref{sec:contraction}, the result of these recursive schemes will not necessarily type check; in the implementation, type checking is done as a final step before the result of the recursive contraction is returned to the user.

If a singular height of a diagram is contracted, this can be promoted to a generalized contraction of the diagram itself:
\[
\begin{tz}
\node (1) at (0,0) {$r_n$};
\node (2) at (1,.5) {$s_n$};
\node (3) at (2,0) {$r_{n+1}$};
\draw [->] (1) to node [above left] {$f_n$} (2);
\draw [->] (3) to node [above right] {$b_n$} (2);
\node (C) at (1,2) {$C$};
\draw [->] (2) to node [right] {$c$} (C);
\end{tz}
\quad\leadsto\quad
\begin{tz}[xscale=1.3]
\node (1) at (0,0) {$r_n$};
\node (2) at (1,.5) {$s_n$};
\node (3) at (2,0) {$r_{n+1}$};
\node (1') at (0,1.5) {$r_n$};
\node (3') at (2,1.5) {$r_{n+1}$};
\draw [->] (1) to node [above left, pos=0.9] {$f_n$} (2);
\draw [->] (3) to node [above right, pos=1] {$b_n$} (2);
\node (C) at (1,2) {$C$};
\draw [->] (2) to node [right] {$c$} (C);
\draw [->] (1') to node [above left] {$c \circ f_n$} (C);
\draw [->] (3') to node [above right] {$c \circ b_n$} (C);
\draw [double equal sign distance] (1) to (1');
\draw [double equal sign distance] (3) to (3');
\end{tz}
\]
On the left, we have a single zigzag, with a morphism out of a singular height; on the right, we have promoted this to a zigzag morphism. It is clear that the zigzag map conditions are satisfied. Here and throughout this section, we show only the central part of the zigzag local to the operation; this procedure is the identity on the rest of the zigzag.

If a regular height is contracted, this can also be promoted to a generalized contraction of the diagram itself, by ``bubbling'' the contraction at the regular height, increasing the size of the diagram by 1:
\[
\begin{tz}[every node/.style={inner sep=2pt}]
\node (1) at (0,.5) {$s_n$};
\node (2) at (1,0) {$r_{n+1}$};
\node (3) at (2,.5) {$s_{n+1}$};
\draw [->] (2) to node [below left] {$b_{n{+}1}$} (1);
\draw [->] (2) to node [below right] {$f_{n{+}1}$} (3);
\node (4) at (1,1.5) {$C$};
\draw [->] (2) to node[right] {$c$} (4);
\end{tz}
\quad\leadsto\quad
\begin{tz}[every node/.style={inner sep=2pt}]
\node (1) at (0,.5) {$s_n$};
\node (2) at (1,0) {$r_{n+1}$};
\node (3) at (2,.5) {$s_{n+1}$};
\draw [->] (2) to node [below left] {$b_{n{+}1}$} (1);
\draw [->] (2) to node [below right] {$f_{n{+}1}$} (3);
\node (4) at (1,2) {$C$};
\node (5) at (2,1.5) {$ r_{n+1}$};
\node (6) at (3,2) {$s_{n+1}$};
\node (7) at (0,1.5) {$r_{n+1}$};
\node (8) at (-1,2) {$s_{n}$};
\draw [double equal sign distance] (1) to (8);
\draw [double equal sign distance] (2) to (7);
\draw [double equal sign distance] (2) to (5);
\draw [double equal sign distance] (3) to (6);
\draw [->] (5) to node [below left, pos=0.2] {$c$} (4);
\draw [->] (7) to node [below right, pos=0.2] {$c$} (4);
\draw [->] (7) to node [above right] {$b_{n{+}1}$} (8);
\draw [->] (5) to node [above left] {$f_{n{+}1}$} (6);
\end{tz}
\]
Again, it is clear that this satisfies the zigzag map axioms.

\subsection{Recursive expansion}

\noindent
As well as generalized contractions, we have the dual notion of \emph{generalized expansions}. While the base case of generalized contraction is colimit, a rich notion explored in detail in Section~\ref{sec:contraction}, the base case of generalized expansion is much more trivial. It is simply the following \textit{ad hoc} process, where two vertices at the same height are perturbed, causing them to have distinct heights:
\begin{equation}
\alignedgraphics{width=2cm}{images/two_beads_singular}
\quad\leadsto\quad
\alignedgraphics{width=2cm}{images/two_beads_generic_1}
\end{equation}
This can be thought of as a ``reverse contraction''. More broadly, we define a \emph{generalized expansion} for a diagram $D$ to be any zigzag map $E \to [e] D$, dual to the notion $D \to [c] C$ of generalized contraction given above.

If we perform a generalized expansion on a regular slice of a diagram, we have the following situation:
\[
\begin{tz}[every node/.style={inner sep=2pt}]
\node (1) at (0,.5) {$s_n$};
\node (2) at (1,0) {$r_{n+1}$};
\node (3) at (2,.5) {$s_{n+1}$};
\draw [->] (2) to node [below left] {$b_{n{+}1}$} (1);
\draw [->] (2) to node [below right] {$f_{n{+}1}$} (3);
\node (4) at (1,-1.5) {$E$};
\draw [->] (4) to node[right] {$e$} (2);
\end{tz}
\]
We might expect that we can perform a recursive trick to extend this expansion to the entire target zigzag. However, this is not possible. The reason is the ``globular'' nature of the theory as we have implemented it, with zigzag maps acting as equalities on regular heights. The original theory of associative $n$\-categories is fundamentally cubical~\cite{Dorn2018, Douglas2019}, and an implementation of that broader theory would relieve this issue.

If we perform a generalized expansion on a singular slice $s_n$ of a diagram, the tool will attempt to find factorizing maps $f', b'$ with $f_n = e \circ f'$ and $b_n = e \circ b'$, as follows:
\[
\begin{tz}
\node (1) at (0,0) {$r_n$};
\node (2) at (1,.5) {$s_n$};
\node (3) at (2,0) {$r_{n+1}$};
\draw [->] (1) to node [above left] {$f_n$} (2);
\draw [->] (3) to node [above right] {$b_n$} (2);
\node (C) at (1,-1) {$E$};
\draw [->] (C) to node [right] {$e$} (2);
\end{tz}
\quad\leadsto\quad
\begin{tz}[xscale=1.3]
\node (1) at (0,0) {$r_n$};
\node (2) at (1,.5) {$E$};
\node (3) at (2,0) {$r_{n+1}$};
\node (1') at (0,1.5) {$r_n$};
\node (3') at (2,1.5) {$r_{n+1}$};
\draw [->] (1) to node [above left, pos=.8] {$f'$} (2);
\draw [->] (3) to node [above right, pos=.7] {$b'$} (2);
\node (C) at (1,2) {$s_n$};
\draw [->] (2) to node [right] {$e$} (C);
\draw [->] (1') to node [above left] {$f_n$} (C);
\draw [->] (3') to node [above right] {$b_n$} (C);
\draw [double equal sign distance] (1) to (1');
\draw [double equal sign distance] (3) to (3');
\end{tz}
\]
If it cannot find such factorizations, it will instead ``bubble'', analogous to the case of recursive contraction on a regular height described above:
\[
\begin{tz}
\node (1) at (0,0) {$r_n$};
\node (2) at (1,.5) {$s_n$};
\node (3) at (2,0) {$r_{n+1}$};
\draw [->] (1) to node [above left] {$f_n$} (2);
\draw [->] (3) to node [above right] {$b_n$} (2);
\node (C) at (1,-1) {$E$};
\draw [->] (C) to node [right] {$e$} (2);
\end{tz}
\quad\leadsto\quad
\begin{tz}[xscale=1]
\node (1) at (0,0.5) {$s_n$};
\node (2) at (1,0) {$E$};
\node (3) at (2,0.5) {$s_n$};
\node (1') at (0,1.5) {$r_n$};
\node (3') at (2,1.5) {$r_{n+1}$};
\draw [->] (2) to node [above right, pos=.7] {$e$} (1);
\draw [->] (2) to node [above left, pos=.7] {$e$} (3);
\node (C) at (1,2) {$s_n$};
\draw [->] (1') to node [above left] {$f_n$} (C);
\draw [->] (3') to node [above right] {$b_n$} (C);
\node (4) at (-1,0) {$r_n$};
\node (5) at (3,0) {$r_{n+1}$};
\draw [->] (4) to node [below right, pos=0] {$f_n$} (1);
\draw [->] (5) to node [below left, pos=0] {$b_n$} (3);
\draw [->] (3) to node [below left] {$\scriptstyle\id$} (C);
\draw [->] (1) to node [below right] {$\scriptstyle\id$} (C);
\draw [double equal sign distance] (1') to (4);
\draw [double equal sign distance] (3') to (5);
\end{tz}
\]
Again, it is clear that these yield zigzag maps.

\subsection{Biased cocones in $\Delta$}

\noindent
The colimit procedure in Section~\ref{sec:contraction} makes use of colimits in $\Delta$ as a subroutine. Since $\Delta$ does not have all colimits, this can sometimes fail, which is unlikely to be the outcome desired by the user. To overcome this, the implementation uses a fallback technique: when the colimit in $\Delta$ does not exist, it will instead supply a mere cocone for the diagram, which it computes by ``breaking the symmetry'' of the input colimit diagram. This method is \textit{ad hoc}, but successfully allows construction of non-symmetrical generalized contractions, such as that described in Example~\ref{ex:opposingunitcounit}, which could otherwise not be produced.

%% file: sectionB-moreproofs.tex

\section{Additional proofs}

\noindent
Here we prove two categorical results used in the main paper.

\begin{proposition}
\label{prop:opfibrationapp}
Let $F:\cA \to \cB$ be a Grothendieck opfibration and let $D:\cJ \to \cA$ be a diagram such that $F D$ has a colimit. If all fibres have $\cJ$-colimits and the base change functor $\sigma_*:F^{-1}(b) \to F^{-1}(b')$ preserves them for all $\sigma:b\to b'$ in $\cB$, then $D$ has a colimit and $F$ preserves it.
\end{proposition}
\begin{proof}
This proof is adapted from a StackExchange post of Pierre Cagne~\cite{Cagne2017}.
Lift the universal cocone morphisms $\lambda^j : FD_j \to x$ to opcartesian morphisms $\phi^j: D_j \to \lambda^j_*(D_j)$ and note that $F(\lambda^j_*(D_j))= x$. Note that every morphism $\sigma:j \to j'$ in $\cJ$ fulfills $\lambda_{j'} \circ FD\sigma  = \lambda_j$. Hence, by opcartesianity of the lifts $\phi^j$, this gives rise to morphisms $\lambda_\sigma: \lambda_*^j(D_j) \to \lambda^{j'}_*(D_{j'})$ with $F(\lambda_{\sigma}) = \id_{x}$. Uniqueness of these morphisms shows that they can be arranged into a functor  $D':\cJ \to F^{-1}(x)$. Let $\mu^j: \lambda^j_*(D_j) \to y$ be the colimit of this functor $D'$ and define $\rho^j = \mu^j\circ \phi^j$. We claim that $\rho^j: D_j \to y$ is a colimiting cocone of $D$. To show this, let $f^j: D_j \to z$ be another cocone. Since $F(f^j): FD_j \to Fz$ is a cocone for $FD$, there is a unique morphisms $\psi: x \to Fz$ such that $F(f^j)  = \psi \circ \lambda^j$. By opcartesianity of $\phi^j: D_j \to \lambda_*^j(D_j)$, this gives rise to morphisms $g^j:\lambda^j_*(D_j) \to z$ with $f^j = g^j \circ \phi^j$ and $F(g^j) = \psi$. Hence, a map of cocones $\{ \rho^j \} \to \{f^j\}$ in $\cC$ for the diagram $D$ amounts to a morphism $y \to z$ in $\cC$ such that:
\begin{equation}\label{eq:proof1} \lambda^j_*(D_j) \to[\mu^j] y \to z = \lambda^j_*(D_j) \to [g^j] z.
\end{equation}
However, since the colimiting cocone $\mu^j: \lambda^j_*(D_j) \to y$ is a colimiting cocone in the fibre $F^{-1}(x)$, whereas the cocone $g^j$ lives over $\psi$, we cannot yet conclude that such a map exists.

For every $j$, let $\overline{\psi}^j: \lambda^j_*D_j \to \psi_*\lambda^j_* D_j$ be an opcartesian lift of $\psi$. Opcartesianity implies that there are morphisms $\overline{g}^j: \psi_*\lambda^j_*D_j \to z$ in the fibre $F^{-1}(Fz)$ such that $g^j = \overline{g}^j \circ \overline{\psi}^j$. By definition of $\psi_*$, it follows that $\overline{g}^j$ forms a cocone for the diagram $\psi_*D': \cJ \to F^{-1}(Fz)$. Since $\psi_*$ preserves $\cJ$-colimits, it follows that $\psi_*\mu^j: \psi_* \lambda^j_*D_j \to \psi_*y$ is a colimiting cocone, and in particular that there is a morphism $\alpha: \psi_*y \to z$ with $\alpha \circ \psi_*\mu^j = \overline{g}^j$ and with $F(\alpha) = \id_{Fz}$. Precomposing $\alpha$ with the opcartesian lift of $\psi$ to a morphism $y \to \psi_*y$ (used for defining $\psi_*$) leads to a morphism $y \to z$ fulfilling equation~\eqref{eq:proof1} as required. Uniqueness of this morphism follows analogously.
\end{proof}

\begin{proposition}
\label{prop:commalimit}
Let $F:\cC\to \cD$ be a functor and let $X$ be an object of $\cD$. Then, the forgetful functor $F/X\to \cC$ from the comma category into $\cC$ preserves connected limits. 
\end{proposition}
\begin{proof} Given a diagram $\cJ \to F/X$, we denote the objects in the image of the diagram by $D_j = (d_j, \alpha_j:d_j \to X)$ for $j \in J$ and the morphisms by $d\sigma: d_j \to d_{j'}$ for $\sigma:j\to j'$, fulfilling $\alpha_{j'}\circ d\sigma =\alpha_{j'}$. Suppose this diagram has a limit $D=(d, \alpha:d\to X)$ with universal cone $\lambda_j : D\to D_j$. Then, $\lambda_j : d\to d_j$ is a universal cone for the diagram $\cJ \to F/X \to \cC$. Indeed, given another cone $\mu_j: r\to d_j$ we pick an object $j_0$ in $\cJ$ and define $\epsilon:= \alpha_{j_0} \circ \mu_{j_0}: r\to X$. Since $\cJ$ is connected and $\{\mu_j\}$ is a cone, the definition of $\epsilon$ does not depend on the chosen object $j_0$. Since $\mu_j: (r, \epsilon)  \to (d_j, \alpha_j)$ is a cone in $F/X$, there is a morphism $f:(r, \epsilon) \to (d, \alpha)$ making the relevant diagrams commute. The proposition follows since any morphism $f:r\to d$ in $\cC$ fulfilling $ \lambda_j \circ f = \mu_j$ is automatically a morphism $(r, \epsilon) \to (d, \alpha)$.
\end{proof}